\documentclass{article}

\usepackage{
	amsthm,
	amsmath,
	accents,
	amsfonts,
	bm,
	enumerate,
    graphicx,
    hyperref,
    shuffle,
    tikz
}

\usepackage[all,pdf]{xy}
\usepackage[shortlabels]{enumitem}

\usepackage[margin=2cm]{geometry}
\usepackage[notref,notcite,color]{showkeys} % show labels.

\setlist[enumerate]{label=\textup{(\arabic*)},ref=\arabic*}

\numberwithin{equation}{section}
\newtheorem{thm}{Theorem}[section]
\newtheorem{pro}[thm]{Proposition}
\newtheorem{crl}[thm]{Corollary}
\newtheorem{lem}[thm]{Lemma}

\theoremstyle{definition}

\newtheorem{rem}[thm]{Remark}

\def\red{\color{red}}
\addtocontents{toc}{\setcounter{tocdepth}{2}}

\title{Twisted algebras of invariant formal power series over monoids}

\author{
    Diego Arcis\footnote{Departamento de Matem\'aticas, Universidad de La Serena, Cisternas 1200 -- 1700000 La Serena, Chile (\texttt{diego.arcis@userena.cl}) - ORCID: \href{https://orcid.org/0000-0003-0027-9323}{0000-0003-0027-9323}.}
\and
    Camilo Gonz\'alez\footnote{Departamento de Matem\'atica, Universidad de Concepci\'on, Casilla 160-C -- 4030000 Concepci\'on, Chile (\texttt{camgonzalezp@udec.cl}) - ORCID: \href{https://orcid.org/0000-0003-0433-7180}{0000-0003-0433-7180}.}
\and
    Sebasti\'an M\'arquez\footnote{Departamento de Matem\'aticas, Universidad Aut\'onoma de Chile, Pedro de Valdivia 425 -- 7500000 Providencia, Chile (\texttt{sebastian.marquez@uautonoma.cl}) - ORCID: \href{https://orcid.org/0009-0000-8043-0112}{0009-0000-8043-0112}.}
}

\date{}
\newcommand{\aspdf}{1}
\newcommand{\asdiagpdf}{1}
\newcommand{\vcdraw}[1]{\vcenter{\hbox{#1}}}

\newcommand{\diagtwo}[1]{
	\centering
	\ifnum #1=1
		\ifnum\asdiagpdf=1
			\includegraphics{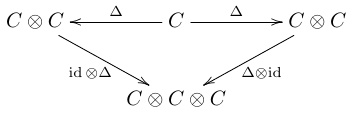}\qquad\qquad
			\includegraphics{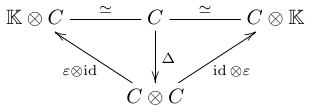}		
		\else
			\xymatrix{
				C\otimes C\ar[rd]_{\id\otimes\Delta}&C\ar[l]_{\kern+1em\Delta}\ar[r]^{\Delta\kern+1em}&C\otimes C\ar[ld]^{\Delta\otimes\id}\\
				&C\otimes C\otimes C
			}
			\qquad\qquad
			\xymatrix{
				\field\otimes C\ar@{-}[r]^{\kern+1em\simeq}&C\ar[d]^\Delta&C\otimes\field\ar@{-}[l]_{\kern-1em\simeq}\\
				&C\otimes C\ar[lu]^{\varepsilon\otimes\id}\ar[ru]_{\id\otimes\varepsilon}
			}
		\fi
	\else\ifnum #1=2
		\xymatrix{
			C\otimes C\ar[r]^{\id\otimes\Delta\kern+1em}&C\otimes C\otimes C\ar@{}[d]|\circlearrowleft&C\otimes C\ar[l]_{\kern+1em\Delta\otimes\id}\\
			&C\ar[lu]^\Delta\ar[ru]_\Delta
		}
		\qquad\qquad
		\xymatrix{
			C\otimes\field\ar@{-}[r]^{\kern+1em\simeq}\ar@{}[rd]|\circlearrowright&C&\field\otimes C\ar@{-}[l]_{\kern-1em\simeq}\ar@{}[ld]|\circlearrowleft\\
			C\otimes C\ar[u]^{\id\otimes\varepsilon}&C\ar[l]^\Delta\ar[u]_\id\ar[r]_\Delta&C\otimes C\ar[u]_{\varepsilon\otimes\id}
		}
	\fi\fi
}

\newcommand{\diagone}[1]{
	\centering
	\ifnum #1=1
		\ifnum\asdiagpdf=1
			\includegraphics{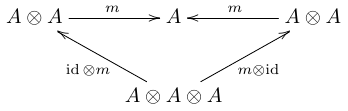}\qquad\qquad
			\includegraphics{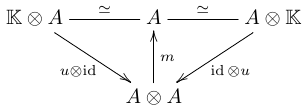}
		\else
			\xymatrix{
				A\otimes A\ar[r]^{\kern+1em m}&A&A\otimes A\ar[l]_m\\
				&A\otimes A\ar[lu]^{\id\otimes m}\otimes A\ar[ru]_{m\otimes\id}
			}
			\qquad\qquad
			\xymatrix{
				\field\otimes A\ar@{-}[r]^{\kern+1em\simeq}\ar[rd]_{u\otimes\id}&A&A\otimes\field\ar@{-}[l]_{\kern-1em\simeq}\ar[ld]^{\id\otimes u}\\
				&A\otimes A\ar[u]_m
			}
		\fi
	\else\ifnum #1=2
		\xymatrix{
			A\otimes A\ar[rd]_m&A\otimes A\otimes A\ar[l]_{\id\otimes m\kern+0.8em}\ar[r]^{\kern+0.8em m\otimes\id}\ar@{}[d]|\circlearrowleft&A\otimes A\ar[ld]^m\\
			&A
		}
		\qquad\qquad
		\xymatrix{
			A\otimes\field\ar@{-}[r]^{\kern+1em\simeq}\ar[d]_{\id\otimes u}\ar@{}[rd]|\circlearrowleft&A\ar[d]^\id&\field\otimes A\ar@{-}[l]_{\kern-1em\simeq}\ar[d]^{u\otimes\id}\ar@{}[ld]|\circlearrowright\\
			A\otimes A\ar[r]_m&A&A\otimes A\ar[l]^m
		}
	\fi\fi
}

\newcommand{\figtwo}{
	\centering
	\ifnum\aspdf=1
		\includegraphics{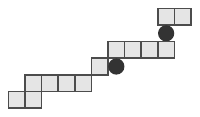}
	\else
		\tabpart{
			\&\&\&\&\&\&\&\&\&\snode\&\snode\\
			\&\&\&\&\&\&\&\&\&\dnode\\			
			\&\&\&\&\&\&\snode\&\snode\&\snode\&\snode\\
			\&\&\&\&\&\snode\&\dnode\\			
			\&\snode\&\snode\&\snode\&\snode\\
			\snode\&\snode
		}
	\fi
}

\newcommand{\figone}{
	\centering
	\ifnum\aspdf=1
		\includegraphics{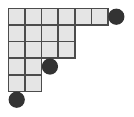}
	\else
		\tabpart{
			\snode\&\snode\&\snode\&\snode\&\snode\&\snode\&\dnode\\
			\snode\&\snode\&\snode\&\snode\\
			\snode\&\snode\&\snode\&\snode\\
			\snode\&\snode\&\dnode\\
			\snode\&\snode\\
			\dnode
		}
	\fi
}

\renewcommand{\operatorname}{\mathsf}
\newcommand{\condition}{{\red$\xi_{\arabic*}$}}

\newcommand{\ubar}[1]{\underaccent{\bar}{#1}}
\newcommand{\refeq}[1]{\stackrel{\eqref{#1}}{=}}
\renewcommand{\operatorname}{\mathsf}

\newcommand{\pols}[2]{#1\langle #2\rangle}
\newcommand{\funcs}[2]{#1\langle\kern-0.2em\langle #2\rangle\kern-0.2em\rangle}
\newcommand{\dbrack}[1]{[\kern-0.2em[#1]\kern-0.2em]}

\newcommand{\Sym}{\operatorname{Sym}}
\newcommand{\QSym}{\operatorname{QSym}}
\newcommand{\NCSym}{\operatorname{NCSym}}
\newcommand{\NCQSym}{\operatorname{NCQSym}}
\newcommand{\sSym}{\operatorname{sSym}}
\newcommand{\sQSym}{\operatorname{sQSym}}
\newcommand{\sNCSym}{\operatorname{sNCSym}}
\newcommand{\sNCQSym}{\operatorname{sNCQSym}}
\newcommand{\MNCSym}{\operatorname{MNCSym}}
\newcommand{\MSym}{\operatorname{MSym}}

\newcommand{\G}{\mathcal{G}}
\renewcommand{\P}{\mathcal{P}}

\newcommand{\N}{\mathbb{N}}
\newcommand{\Q}{\mathbb{Q}}
\newcommand{\Z}{\mathbb{Z}}
\newcommand{\field}{\mathbb{K}}

\newcommand{\sym}{\mathfrak{S}}

\newcommand{\id}{\operatorname{id}}
\newcommand{\GL}{\operatorname{GL}}
\newcommand{\Aut}{\operatorname{Aut}}
\newcommand{\bdg}{\operatorname{bdg}}
\newcommand{\ind}{\operatorname{ind}}
\newcommand{\inv}{\operatorname{inv}}
\newcommand{\sgn}{\operatorname{sgn}}
\newcommand{\pack}{\operatorname{pack}}
\newcommand{\sort}{\operatorname{sort}}
\newcommand{\supp}{\operatorname{supp}}
\renewcommand{\ker}{\operatorname{ker}}
\renewcommand{\alph}{\operatorname{alph}}

\newcommand{\hsf}{\mathsf{h}}
\newcommand{\psf}{\mathsf{p}}
\newcommand{\qsf}{\mathsf{q}}
\newcommand{\ssf}{\mathsf{s}}

\newcommand{\conditions}[1]{
	\ifnum#1=0
		$M$ is a monoid equipped with a twisting $\tau:M\times M\to\field$, and $\G$ is a commutative monoid endowed with a bicharacter $\phi:\G\times\G\to\field^\times$.
	\else\ifnum#1=1
		$M$ is locally finite with respect to $\tau$ and is equipped with a $\G$-grading $\{\ubar{M}_x\mid x\in\G\}$ for the trivial twisting.
	\else\ifnum#1=2
		$M$ is graded with respect to $\tau$ with degree map $\deg:M\to\N$ and grading $\{M^d\mid d\in\N\}$, where each homogeneous component $M^d$ intersects only finitely many fibres $\ubar{M}_x$.
	\else\ifnum#1=3
		$G$ is a group equipped with a quasi-multiplicative bigraded action $\rho:G\to\sym(M)$ together with a weight function $\omega:G\times M\to\field^\times$ such that $(\tau,\omega)$ is a multiplicative pair.
	\else\ifnum#1=4
		$T$ denotes a system of representatives of the $\omega$-admissible $G$-orbits, and each fibre $M^d$ intersects only finitely many of these orbits for all $d\in\N$.
	\else\ifnum#1=5
		$M$ is equipped with a coproduct structure $(\,\psf:M\to M\times M,\,\hsf:G\times G\to G\,)$.
	\else\ifnum#1=6
		$M$ is equipped with an invariant-retraction pair $(S\subseteq M\times M,\,U\subseteq M)$.
	\else\ifnum#1=7
	\fi\fi\fi\fi\fi\fi\fi\fi
}

\begin{document}

\maketitle

\begin{abstract}
We develop a general monoid-theoretic framework to construct and study graded and color Hopf algebras of formal power series of bounded degree invariant under group actions. Rather than working with external alphabets of commuting, noncommuting, or anti-commuting variables, we start from an abstract monoid $M$ equipped with a $2$-cocycle twisting $\tau:M\times M\to\field$, which is allowed to vanish, a bigrading over a commutative monoid $\G\times\N$ endowed with a bicharacter $\phi$, and a group $G$ acting on $M$ via a bigraded action $\rho$ that is quasi-multiplicative with respect to $\tau$, coupled with a $1$-cocycle weight function $\omega:G\times M\to\field^\times$. By formulating local finiteness, gradings, quasi-multiplicativity, and ideals relative to the non-vanishing of $\tau$, we prove that the invariant space $\funcs{\field}{M}^G$ forms a well-defined bigraded subalgebra of bounded-degree formal power series with a canonical monomial basis indexed by $\omega$-admissible orbits. Furthermore, we establish four functorial constructions on $\funcs{\field}{M}^G$: embeddings into larger invariant algebras via restricted actions, color tensor product decompositions for direct products and diagonal actions, connected color Hopf algebra structures induced by internal monoid duplication, and quotient color Hopf algebras governed by invariant-retraction pairs with an explicit description of their kernel bi-ideals. As applications, we uniformly recover and extend the classical and $(\pm1)$-twisted Hopf algebras of symmetric and quasisymmetric functions in commuting and noncommuting variables ($\Sym_q$, $\QSym_q$, $\NCSym_q$, $\NCQSym_q$), MacMahon symmetric functions, the full square of symmetric and quasisymmetric Hopf algebras in superspace ($\sNCSym_q$, $\sNCQSym_q$, $\sSym_q$, $\sQSym_q$), symmetric and quasisymmetric functions in partially commutative variables, and signed-symmetric functions associated with the hyperoctahedral group.
\end{abstract}

\medskip{\small
}

\tableofcontents

\section{Introduction}

The theory of symmetric functions stands as a central pillar of modern algebraic combinatorics, with deep connections to representation theory, algebraic geometry, and mathematical physics. Over the past decades, the classical ring of symmetric functions $\Sym$~\cite{Mac99,St99} has evolved into a rich hierarchy of algebraic structures, most notably quasisymmetric functions $\QSym$~\cite{Ges84,LuMyWi13}, symmetric functions in noncommuting variables $\NCSym$~\cite{Wol36,RoSa06,BeReRoZa08}, quasisymmetric functions in noncommuting variables $\NCQSym$~\cite{BeZa09,NoThWi10}, as well as their quantum, multi-parameter, and multi-variable counterparts~\cite{Thi96,Ro01,RoRoSt02}. More recently, motivated by integrable systems in mathematical physics, this landscape has expanded to superspace through the introduction of symmetric and quasisymmetric functions in commuting, noncommuting and anti-commuting, namely $\sSym$~\cite{DeLaMa03,DeLaMa04}, $\sQSym$~\cite{FiLaPi19}, $\sNCSym$~\cite{ArGoMa25}, and $\sNCQSym$~\cite{PreArGoMa26A}.

Despite the diversity of these algebras, their traditional constructions rely almost exclusively on a variable-centric paradigm. One typically starts with an infinite alphabet of commuting, noncommuting, or anti-commuting variables, forms a ring of formal power series of bounded degree, and extracts invariants under the action of a group of permutations, or local order-preserving shifts, acting on the variable indices. While this perspective has been remarkably fruitful, it depends on case-by-case heuristics that frequently obscure the underlying algebraic mechanisms and encounter fundamental obstacles when signs, quantum parameters, or nilpotent variables enter the picture.

A first obstacle arises from deformed relations and nilpotence. In superspace or quantum settings, variables are subject to relations such as $\theta_i\theta_j=-\theta_j\theta_i$, $\theta_i^2=0$, or $xy=q\,yx$. Because neither minus signs, scalar multiples, nor zero-divisors of the form $\theta_i^2=0$ can exist inside a cancellative monoid of words, variable-centric approaches are forced to abandon pure monoid combinatorics and work inside quotient algebras of formal power series, where the notion of a monomial orbit becomes delicate due to sign cancellations. A second issue concerns non-multiplicative group actions. Several fundamental algebras of invariants, including $\QSym$, $\NCQSym$, $\sQSym$, and $\sNCQSym$, are defined as invariant subspaces under Hivert's quasi-symmetrizing action~\cite[Section~3]{Hi00}. Unlike the standard permutation action, this action is not a monoid homomorphism, raising the conceptual question of what structural property of a non-multiplicative group action guarantees that the product of two invariants remains invariant. A third difficulty appears in the construction of Hopf algebra structures, which are almost universally introduced through the heuristic of alphabet duplication, that is, evaluating a series on two ordered, disjoint copies of the variable set. While intuitively appealing, this procedure is rarely formalized at the level of the underlying monoids, making rigorous proofs of coassociativity cumbersome, especially in braided or super-graded settings where cross-factors commute up to signs or powers of $q$. Finally, canonical projections from noncommuting to commuting variables, such as $\NCSym\to\Sym$ or $\sNCSym\to\sSym$, are informally understood as forcing variables to commute or anti-commute via relations like $xy=yx$ or $xy=-yx$. However, as observed in~\cite[Subsection~3.6]{HiNoTh08}, the commutator $xy-yx$ is not a symmetric function, so this relation does not generate an ideal inside the algebra of invariants. Even more critically, in superspace or twisted settings, anti-commutation relations of the form $xy=-yx$ do not arise from any monoid congruence at all, and certain noncommutative monomial invariants project to zero because their symmetrized orbits collapse by sign cancellation.

In this paper, we resolve these tensions simultaneously by shifting the focus from external alphabets of variables to the internal algebra of an abstract monoid. Our guiding philosophy is that the monoid provides everything: every deformed product, sign convention, non-multiplicative symmetry, canonical basis, tensor decomposition, Hopf coproduct, and algebraic projection can be encoded directly on a monoid $M$ and lifted functorially to its space of formal power series.

A central conceptual innovation of our work is the systematic relativization of monoid-theoretic notions with respect to a $2$-cocycle twisting $\tau:M\times M\to\field$. Rather than forcing relations like $uv=-vu$, $uv=q\,vu$, or $\theta_i^2=0$ into the monoid $M$ itself, we keep $M$ purely combinatorial and encode all deformations, signs, and nilpotencies into $\tau$, which we allow to take the value zero. For instance, fermionic variables in superspace can be modeled over an elementary abelian group $\Theta^\wedge$ where $\theta_r^2=1$ in the monoid, while the twisting $\tau_\wedge(\theta_I,\theta_J)$ vanishes whenever $I\cap J\neq\emptyset$ and records the Koszul sign $(-1)^{\inv(\theta_I,\theta_J)}$ when $I\cap J=\emptyset$. Because $\tau$ may vanish, the monoid $M$ can possess factorizations or relations that are annihilated in the twisted convolution product
\[(g\cdot h)(u)=\sum_{vw=u}\tau(v,w)\,g(v)\,h(w).\]
Consequently, classical monoid hypotheses are too rigid and would exclude superspace and nilpotent settings. By requiring algebraic and combinatorial properties of $M$ to hold only on pairs $(u,v)$ with $\tau(u,v)\neq0$, we obtain a flexible yet completely rigorous foundation.

We develop this framework progressively across Sections~\ref{167} to~\ref{169}, introducing each structural condition precisely when the algebra demands it. In Section~\ref{167}, we lay the foundations for twisted power series and their invariant subalgebras. We begin in Subsection~\ref{007} by considering a monoid $M$ equipped with a $2$-cocycle twisting $\tau:M\times M\to\field$, together with a commutative grading monoid $\G$ and a bicharacter $\phi:\G\times\G\to\field^\times$ (Condition~\eqref{015}). Working in the category of $\G$-graded vector spaces endowed with the $\phi$-colored tensor product $A\otimes_\phi B$ allows us to incorporate classical, quantum, and super-braided signs directly into the ambient linear algebra. For the twisted convolution product on $\field^M$ to be well defined, $M$ need not be locally finite in the classical sense; indeed, the group $\Theta^\wedge$ is not locally finite since the relation $\theta_r^2=1$ produces infinitely many factorizations of any element. Instead, it suffices that every $u\in M$ admits only finitely many factorizations $vw=u$ with $\tau(v,w)\neq0$, a property we call local finiteness with respect to $\tau$. Pairing this requirement with a $\G$-grading $\{M_x\mid x\in\G\}$ for the trivial twisting yields Condition~\eqref{016}, under which $\field^M$ becomes a well-defined algebra containing the twisted monoid algebra $\pols{\field}{M}$ as a $\G$-graded subalgebra.

To handle the infinite sums inherent in symmetric functions while preserving the $\G$-grading on tensor products, we introduce in Subsection~\ref{188} an $\N$-grading $\{M^d\mid d\in\N\}$ with respect to $\tau$, meaning that the degree map $\deg:M\to\N$ is additive whenever $\tau(u,v)\neq0$, such that each homogeneous fibre $M^d$ intersects only finitely many $\G$-components $M_x$ (Condition~\eqref{017}). Here again, the relativization with respect to $\tau$ is decisive: in the fermionic monoid $\Theta^\wedge$, the degree map $\deg(\theta_I)=|I|$ fails to be additive on products where $\theta_r^2=1$, but it is strictly additive whenever $\tau_\wedge(\theta_I,\theta_J)\neq0$, that is, when $I\cap J=\emptyset$. Condition~\eqref{017} guarantees that the space $\funcs{\field}{M}$ of formal power series of bounded $\N$-degree is closed under the twisted convolution product and forms a $\G$-graded subalgebra of $\field^M$ (Proposition~\ref{005}).

In Subsection~\ref{189}, we determine when a group $G$ acting on $M$ via $\rho:G\to\sym(M)$ induces a subalgebra of invariants $\funcs{\field}{M}^G$. To accommodate signed and quantum symmetries, we twist the induced action on $\field^M$ by a $1$-cocycle weight function $\omega:G\times M\to\field^\times$, setting $\bm{\rho}_\sigma(g)(u)=\omega(\sigma,\rho_{\sigma^{-1}}(u))\,g(\rho_{\sigma^{-1}}(u))$. Crucially, to capture quasisymmetric-type actions, we do not require $\rho_\sigma$ to be a monoid homomorphism. Instead, we define $\rho$ to be quasi-multiplicative with respect to $\tau$ if, for every $\sigma\in G$ and every target product $w\in M$, there exist $\alpha,\varsigma\in G$, depending only on $\sigma$ and $w$, such that $\rho_\sigma(uv)=\rho_\alpha(u)\rho_\varsigma(v)$ for all factorizations $uv=w$ with $\tau(u,v)\neq0$. When combined with the compatibility relation $\tau(u,v)\,\omega(\sigma,uv)=\tau(\rho_\alpha(u),\rho_\varsigma(v))\,\omega(\alpha,u)\,\omega(\varsigma,v)$ whenever $\tau(u,v)\neq0$, which defines $(\tau,\omega)$ as a multiplicative pair with respect to $\tau$, the map $(u,v)\mapsto(\rho_\alpha(u),\rho_\varsigma(v))$ restricts to a bijection between the finite $\tau$-non-zero factorization fibres of $\rho_{\sigma^{-1}}(w)$ and $w$. This bijection allows the convolution sum to be re-indexed cleanly, proving that if $G$ acts via a bigraded quasi-multiplicative action with a multiplicative pair $(\tau,\omega)$ (Condition~\eqref{018}), then the invariant space $\funcs{\field}{M}^G$ is a $\G$-graded subalgebra of $\funcs{\field}{M}$ (Proposition~\ref{029}).

In Subsection~\ref{190}, we construct canonical linear bases for $\funcs{\field}{M}^G$. Due to the presence of the weight function $\omega$, not every $G$-orbit in $M$ can support a non-zero invariant: an orbit $Gt$ supports a non-zero monomial invariant $m_t$ if and only if $\omega(\sigma,t)=1$ for all $\sigma$ in the stabilizer $G_{(t)}$, in which case we say that $t$ is $\omega$-admissible. Under the assumption that each degree component $M^d$ intersects only finitely many $\omega$-admissible $G$-orbits (Condition~\eqref{051}), we prove that the orbit sums $\{m_t\mid t\in T\}$ indexed by a system of admissible representatives $T$ form a canonical $\G$-homogeneous linear basis of the algebra $\funcs{\field}{M}^G$ (Theorem~\ref{050}).

In Section~\ref{168}, we develop four general constructions that explain the structural relationships connecting different families of invariant algebras. We begin in Subsection~\ref{171} with a general restriction mechanism that accounts for the embeddings of symmetric algebras into quasisymmetric ones. Given a symmetric generating set $S$ of $G$, we assign to each $s\in S$ a subset $M_s\subseteq M$ that is stable and multiplicative with respect to $\tau$ under $\rho_s$, and whose complement $M\setminus M_s$ is an ideal with respect to $\tau$, meaning that if $uv\in M_s$ and $\tau(u,v)\neq0$, then both $u$ and $v$ belong to $M_s$. Defining $\bar{\rho}_s$ to act as $\rho_s$ on $M_s$ and as the identity on $M\setminus M_s$ generates a restricted group $\bar{G}$ acting on $M$ with a restricted weight $\bar{\omega}$. Once again, formulating the ideal condition with respect to $\tau$ is essential for superspace: in the monoid $X^\theta=\Theta^\wedge\times X^*$, the complement $X^\theta\setminus X_s^\theta$ is not a classical monoid ideal because the relation $\theta_r^2=1$ allows fermionic letters to cancel and modify the support of an element, yet whenever such a cancellation occurs in a product $uv$ one has $\tau^\times(u,v)=0$, ensuring that $X^\theta\setminus X_s^\theta$ is an ideal with respect to $\tau^\times$. We prove that $(\bar{G},\bar{\rho},\bar{\omega})$ automatically satisfies Conditions~\eqref{018} and~\eqref{051}, so that $\funcs{\field}{M}^{\bar{G}}$ is an algebra of invariants containing $\funcs{\field}{M}^G$ as a subalgebra, and each $G$-monomial basis element $m_t$ expands explicitly as a $\bar{\omega}$-weighted sum of $\bar{G}$-monomial basis elements $\bar{m}_{t'}$ (Theorem~\ref{057}).

Next, in Subsections~\ref{120} and~\ref{183}, we examine direct products and diagonal actions. Given monoids $M_1,\dots,M_k$ satisfying Conditions~\eqref{015} through~\eqref{051} under groups $G_1,\dots,G_k$, we show that the direct product $M_1\times\cdots\times M_k$ inherits a combined twisting $\tau_\phi=\tau^\times\tau^\phi$ incorporating the bicharacter $\phi$, and that the invariant algebra under the product group $G_1\times\cdots\times G_k$ is canonically isomorphic to the $\phi$-colored tensor product $\funcs{\field}{M_1}^{G_1}\otimes_\phi\cdots\otimes_\phi\funcs{\field}{M_k}^{G_k}$ (Theorem~\ref{088} and Corollary~\ref{097}). Complementarily, when a single group $G$ acts diagonally on a direct product $M_1\times\cdots\times M_k$ or on a free product $M_1*\cdots*M_k$, we establish how twistings, gradings, and weights combine (Propositions~\ref{076} through~\ref{068}), providing a systematic method to build multi-alphabet and superspace invariant algebras from simpler factors.

In Subsection~\ref{191}, we combine the product isomorphism $\Phi_2:\funcs{\field}{M\times M}^{G\times G}\to\funcs{\field}{M}^G\otimes_\phi\funcs{\field}{M}^G$ with duality to formalize alphabet duplication entirely inside the monoid as a coproduct structure $(\psf:M\to M\times M,\hsf:G\times G\to G)$ (Condition~\eqref{127}). Here $\psf$ is a bigraded, finite-fibre, twisting-preserving monoid homomorphism intertwined with $\hsf$ via $\psf\circ\rho_{\hsf(\sigma,\varsigma)}=(\rho_\sigma\times\rho_\varsigma)\circ\psf$, satisfying fibrewise coassociativity and counitality on admissible orbits. We prove that Condition~\eqref{127} induces a connected $\G$-graded color Hopf algebra structure on $\funcs{\field}{M}^G$ with coproduct $\Delta=\Phi_2\circ\nabla$ (Theorem~\ref{130}). Furthermore, we establish compatibility criteria under which a coproduct structure transfers to a restricted action $\bar{G}$, making the inclusion $\funcs{\field}{M}^G\hookrightarrow\funcs{\field}{M}^{\bar{G}}$ an embedding of connected color Hopf algebras (Proposition~\ref{170}), and we provide a canonical construction of $\psf$ (Proposition~\ref{036}).

In Subsection~\ref{176}, we address algebraic projections between invariant algebras by introducing the notion of an invariant-retraction pair $(S,U)$ on $M$ (Condition~\eqref{161}), consisting of a $G$-stable, finite-class monoid congruence $S$ and a section $U\subseteq M$ of class representatives compatible with $\tau$, the gradings, $\omega$, and $(\psf,\hsf)$. Because the induced quotient twisting $\tau^{_S}$ and quotient weight $\omega^{_S}$ on $M/S$ are defined by pulling back to the section $U$, a congruence $S$ on $M$ that merely commutes generators in the monoid automatically generates inversion signs or $q$-weights in the quotient algebra $\funcs{\field}{M/S}^G$, bypassing the impossibility of quotienting a monoid by signed relations like $xy=-yx$. We prove that the linear projection $[\ ]_G:\funcs{\field}{M}^G\to\funcs{\field}{M/S}^G$ is a surjective homomorphism of connected color Hopf algebras (Theorem~\ref{153}). Moreover, we show that each monomial invariant $m_t$ maps either to $c_t\,m_{[t]}$ with $c_t=|[t]\cap Gt|\geq1$ when $[t]$ is $\omega^{_S}$-admissible, or to $0$ when $[t]$ is not $\omega^{_S}$-admissible (Lemma~\ref{160}), and we give an explicit generating set for the kernel bi-ideal $J=\ker([\ ]_G)$ (Proposition~\ref{166}).

Finally, in Section~\ref{169}, we apply our framework across three broad settings to demonstrate both its unifying power and its ability to generate new structures. In Subsection~\ref{179}, we apply our results to the free monoid $X^*$ and its commutative quotient $X^\oplus$ over a rational-indexed alphabet $X=\{x_r\mid r\in\Q\}$ equipped with the inversion twisting $\tau_q$ and weight $\omega_q$. We simultaneously recover the classical Hopf algebras $\NCSym$, $\NCQSym$, $\Sym$, and $\QSym$ at $q=1$ together with their $(\pm1)$-twisted counterparts $\NCSym_q$, $\NCQSym_q$, $\Sym_q$, and $\QSym_q$, assembling them into a commutative square of connected color Hopf algebras, and we derive the noncommutative and commutative MacMahon symmetric functions via $k$-fold free and direct products. In Subsection~\ref{192}, we construct the complete square of $(\pm1)$-twisted symmetric and quasisymmetric Hopf algebras in superspace, namely $\sNCSym_q$, $\sNCQSym_q$, $\sSym_q$, and $\sQSym_q$, by taking the diagonal product $X^\theta=\Theta^\wedge\times X^*$ of the fermionic group $\Theta^\wedge$ with $X^*$, illustrating how $\tau$-relative local finiteness, $\tau$-relative ideals, and invariant-retraction pairs naturally govern Grassmann variables, super-admissibility, and projections to superpartitions and supercompositions. We conclude in Subsection~\ref{193} by constructing symmetric and quasisymmetric functions in partially commutative variables as well as signed-symmetric functions in noncommuting and commuting variables associated with the hyperoctahedral group.\medskip

Throughout this paper, $\field$ denotes a field of characteristic zero, and all linear algebraic structures are defined over $\field$. We write $\N$ for the commutative monoid of nonnegative integers under addition; for $n\in\N$, we write $[n]_0$ for $\{0,\dots,n\}$, and for $n\geq1$, we write $[n]$ for $\{1,\ldots,n\}$. The group of permutations of a set $X$ is denoted by $\sym(X)$.

%\begin{itemize}
%\item[\eqref{015}]\conditions{0}
%\item[\eqref{016}]\conditions{1}
%\item[\eqref{017}]\conditions{2}
%\item[\eqref{018}]\conditions{3}
%\item[\eqref{051}]\conditions{4}
%\item[\eqref{127}]\conditions{5}
%\item[\eqref{161}]\conditions{6}
%\end{itemize}

\section{Background}

To keep this paper self-contained and establish a unified terminology, this section collects the essential algebraic preliminaries on monoids, gradings, and color Hopf algebras that form the foundation of our work.

\subsection{Monoids}

We begin by recalling the fundamental constructions in monoid theory, specifically congruences, free products, and gradings, that will serve as the combinatorial ambient spaces for our formal power series. For further details on semigroup theory, we refer the reader to~\cite{How95}.

Given a monoid $M$, we write $1$, or $1_{_M}$ when necessary, for its identity element. Throughout this subsection, $M$ and $N$ denote two monoids.

\subsubsection{Congruences}\label{002}

If $R\subseteq M\times M$ is a relation on $M$, it is common to write $u=v$ for any pair $(u,v)\in R$, referring to these pairs as \emph{relations}. When $R$ is an equivalence relation on $M$, the class of an element $u\in M$ is denoted by $[u]_R$, or simply $[u]$. Given a system of representatives $T\subseteq M$, we typically identify the class $[u]$ with its representative $u\in T$.

A \emph{congruence} on $M$ is an equivalence relation $R\subseteq M\times M$ compatible with multiplication, meaning that $[u][v]=[uv]$ for all $u,v\in M$. In this case, the quotient set $M/R$ inherits a natural monoid structure. For any relation $S\subseteq M\times M$ on $M$, we denote by $\bar{S}$ the \emph{congruence closure} generated by $S$, defined as the smallest congruence on $M$ containing $S$. Given a monoid homomorphism $\psi:M\to N$, a congruence $R$ on $M$ is said to be \emph{stable} under $\psi$ if $\psi(u)=\psi(v)$ for all $(u,v)\in R$. When $R$ is stable under $\psi$, the universal property of quotients guarantees the existence of a unique monoid homomorphism $\hat{\psi}:M/R\to N$ defined by $\hat{\psi}([u])=\psi(u)$~\cite[Theorem~1.5.3]{How95}. If $S$ is another congruence on $M$ such that $R\subseteq S$, the set $S/R:=\{([u]_R,[v]_R)\mid (u,v)\in S\}$ defines a congruence on $M/R$ such that $(M/R)/(S/R)\simeq M/S$~\cite[Theorem~1.5.4]{How95}. Consequently, we naturally identify $[[u]_R]_{S/R}$ with $[u]_S$.

\subsubsection{Free products}

For a set $X$, we denote by $X^*$ the \emph{free monoid} generated by $X$, that is, the collection of finite words with \emph{letters} in $X$, under the product given by \emph{concatenation}. The \emph{free product} of $M$ and $N$ is the quotient monoid~\cite[Subsection~8.2]{How95}:\[M*N=(M\sqcup N)^*/\bar{R}\,,\quad\text{where}\quad R=\{(1_{_M},\varepsilon),(1_{_N},\varepsilon)\}\cup\{(u*v,uv)\mid u,v\in M\text{ or }u,v\in N\}.\]More precisely, every non-identity element of $M*N$ can be uniquely written as an alternating product $u_1v_1\cdots u_kv_k$, where $u_1,\ldots,u_k\in M\setminus\{1\}$ and $v_1,\ldots,v_k\in N\setminus\{1\}$, except possibly that $u_1$ or $v_k$ may be the identity element. This construction preserves the internal multiplicative structure of each monoid, while introducing no relations between elements of $M$ and $N$. In particular, in $M*M$, the two copies of $M$ interact freely. On the other hand, if $[M,N]$ denotes the congruence generated by $\{(uv,vu)\mid(u,v)\in M\times N\}$, then $M\times N\simeq(M*N)/[M,N]$. Given monoid homomorphisms $\phi:M\to L$ and $\psi:N\to L$, there exists a unique homomorphism $\phi*\psi:M*N\to L$ such that $(\phi*\psi)(u)=\phi(u)$ and $(\phi*\psi)(v)=\psi(v)$ for all $u\in M$ and $v\in N$~\cite[Proposition~8.2.1]{How95}.

This construction naturally extends to arbitrary multiple free products of monoids~\cite[Subsection~8.2]{How95}. In particular, if $M_1,\ldots,M_k$ are monoids, every non-identity element in their free product $M_1*\cdots*M_k$ is uniquely written as an alternating product $u_1\cdots u_q$, with $u_1,\ldots,u_q\in(M_1\sqcup\cdots\sqcup M_k)\setminus\{1\}$, such that $u_i$ and $u_{i+1}$ do not belong to the same monoid for all $i\in[q-1]$. Given $u,v\in M_1*\cdots*M_k$, the pair formed by the last factor in the alternating product of $u$ and the first factor in the alternating product of $v$ is called the \emph{central factors} of $u$ and $v$.

\subsection{Color bialgebras}

We now turn to color bialgebras and color Hopf algebras, which provide the algebraic framework to handle both superspace sign conventions and deformed product structures. For a comprehensive treatment of color algebraic structures, we refer the reader to~\cite{AndAngBag14,ChePetOys06,DaNaRa01}.

\subsubsection{Graded algebras}

Recall that a vector space $V$ is said to be \emph{$\P$-graded}, for some commutative monoid $\P$, if there exists a family of subspaces $\{V_x\mid x\in\P\}$, referred to as a \emph{$\P$-grading}, such that $V$ admits the direct sum decomposition:\[V=\bigoplus_{x\in\P}V_x.\]In this context, each subspace $V_x$ is called a \emph{homogeneous component} of $V$, and every nonzero element $v\in V_x$ is said to be \emph{homogeneous} of \emph{$\P$-degree} $x$. Furthermore, the map\[
\nu:\bigsqcup_{x\in\P}(V_x\setminus\{0\})\to\P\]defined by $\nu(v)=x$ if and only if $v\in V_x$, is called the \emph{degree map} of the $\P$-grading. For any $v\in V$, we denote by $v_x$ its projection onto the homogeneous component $V_x$. A subspace $U\subseteq V$ is itself \emph{$\P$-graded} if $v_x\in U$ for all $x\in\P$ and all $v\in U$. Given two $\P$-graded vector spaces $V$ and $W$, a linear map $\psi:V\to W$ is said to be \emph{$\P$-graded} if $\psi(V_x)\subseteq W_x$ for all $x\in\P$. In the specific case where $\P=\N$, the prefix "$\N$-" is typically omitted from the related terminology. Note that the ground field $\field$ is naturally $\P$-graded by setting $\field_1=\field$ and $\field_x=\{0\}$ for all $x\in\P\setminus\{1\}$. More generally, if $V_1,\ldots,V_k$ are $\P$-graded vector spaces with respective degree maps $\nu_1,\ldots,\nu_k$, their tensor product $V_1\otimes\cdots\otimes V_k$ is $\P$-graded with the degree map defined on homogeneous simple tensors by $v_1\otimes\cdots\otimes v_k\mapsto\nu_1(v_1)\cdots\nu_k(v_k)$. Hence,\begin{equation}\label{096}
V_1\otimes\cdots\otimes V_k=\bigoplus_{x\in\P}(V_1\otimes\cdots\otimes V_k)_x,\quad\text{where}\quad(V_1\otimes\cdots\otimes V_k)_x=\bigoplus_{x_1\cdots x_k=x}(V_1)_{x_1}\otimes\cdots\otimes(V_k)_{x_k}.\end{equation}

By the universal property of tensor products, an algebra $(A,\cdot\,,1)$ is equivalently a vector space $A$ equipped with two linear maps: the \emph{product} $m:A\otimes A\to A$, defined by $m(a\otimes b)=ab$, and the \emph{unit} $u:\field\to A$, defined by $u(\alpha)=\alpha 1$, satisfying the following associativity and unitality diagrams:\begin{equation}\label{008}\diagone{1}\end{equation}We say that $A$ is \emph{$\P$-graded} if its underlying vector space is endowed with a $\P$-grading $\{A_x\mid x\in\P\}$ such that the product and unit are $\P$-graded linear maps, where $A\otimes A$ is equipped with the $\P$-grading defined in \eqref{096}. Equivalently, $A_x\,A_y\subseteq A_{xy}$ for all $x,y\in\P$, and the identity element satisfies $1\in A_1$. In this context, an algebra homomorphism between $\P$-graded algebras is called \emph{$\P$-graded} if it is a $\P$-graded linear map. A \emph{$\P$-graded subalgebra} of a $\P$-graded algebra $A$ is a subalgebra $B \subseteq A$ such that $B=\bigoplus_{x\in\P}(B\cap A_x)$; that is, for every element $b\in B$, all its homogeneous components also belong to $B$. A \emph{graded ideal} of $A$ is an ideal $I\subseteq A$ that is also a $\P$-graded subspace, meaning that $I_x=I\cap A_x$ for all $x\in\P$. In this case, the quotient algebra $A/I$ inherits a natural $\P$-grading with homogeneous components given by $(A/I)_x=A_x/I_x$.

\subsubsection{Color Hopf algebras}\label{012}

In the remainder of this subsection, $\G$ denotes a commutative monoid equipped with a \emph{bicharacter}~\cite[Subsection~1.1]{AndAngBag14}; that is, a map $\phi:\G\times\G\to\field^\times$ satisfying the following conditions:\begin{equation}\label{082}\phi(xy,z)=\phi(x,z)\phi(y,z),\qquad\phi(x,yz)=\phi(x,y)\phi(x,z).\end{equation}From \eqref{082}, we immediately deduce that $\phi(1,x)=1=\phi(x,1)$ for all $x\in\G$. We further assume that every vector space $V$ is endowed with a $\G$-grading $\{\ubar{V}_x\mid x\in\G\}$. Formally, we work within the category of $\G$-graded vector spaces, where all linear maps and algebraic structures are implicitly assumed to be $\G$-graded. The $\G$-degree of a homogeneous vector $v$ is denoted by $|v|$, so that for each $x\in\G$, we have $\ubar{V}_x\setminus\{0\}=\{v\in V\mid|v|=x\}$.

Since $\phi$ is a bicharacter, if $A$ and $B$ are $\G$-graded algebras, their tensor product $A\otimes B$ can also be endowed with an algebra structure, denoted $A\otimes_\phi B$, by linearly extending the product defined on homogeneous tensors by:\begin{equation}\label{080}(a\otimes b)(c\otimes d)=\phi(|b|,|c|)\,ac\otimes bd.\end{equation}Indeed, the $\G$-grading is compatible with this multiplication since $\G$ is commutative and $\phi$ takes values in $\field^\times$; specifically, $|(a\otimes b)(c\otimes d)|=|ac\otimes bd|=|ac||bd|=|ab||cd|=|a\otimes b||c\otimes d|$.

\begin{rem}\label{081}
The multiplicative structure in \eqref{080} extends naturally to arbitrary tensor products. Specifically, if $A_1,\ldots,A_k$ are $\G$-graded algebras, the tensor product $A_1\otimes\cdots\otimes A_k$ can be endowed with the algebra structure $A_1\otimes_\phi\cdots\otimes_\phi A_k$ by linearly extending the product defined on homogeneous tensors as:\[(a_1\otimes\cdots\otimes a_k)\cdot(b_1\otimes\cdots\otimes b_k)=\Bigg(\prod_{i\,<\,j}\phi(|a_j|,|b_i|)\Bigg)a_1b_1\otimes\cdots\otimes a_kb_k.\]Moreover, $A_1\otimes_\phi\cdots\otimes_\phi A_k$ inherits a natural $\P$-graded algebra structure with the $\P$-grading given in~\eqref{096}.
\end{rem}

Dually, a coalgebra is defined by reversing the arrows in the diagrams defining an algebra in \eqref{008}. A \emph{coalgebra} is a vector space $C$ equipped with two linear maps $\Delta:C\to C\otimes C$ and $\varepsilon:C\to\field$, called \emph{coproduct} and \emph{counit}, respectively, satisfying the following \emph{coassociativity} and \emph{counitality} diagrams:\[\diagtwo{1}\]Since the coproduct and counit are $\G$-graded linear maps, they satisfy $\Delta(\ubar{C}_x)\subseteq\bigoplus_{yz=x}\ubar{C}_y\otimes\ubar{C}_z$ for each $x\in\G$, while the counit satisfies $\varepsilon(\ubar{C}_x)=\{0\}$ for all $x\in\G\setminus\{1\}$. 

As is customary, for each homogeneous element $c\in C$, we employ \emph{Sweedler's notation} \cite[Section~1.2]{Sw69} to write the coproduct as $\Delta(c)=c_{_{(1)}}\otimes c_{_{(2)}}$, where $c_{_{(1)}}\otimes c_{_{(2)}}$ represents the finite linear combination of simple tensors defining $\Delta(c)$. In this context, it is understood that $c_{_{(1)}}$ and $c_{_{(2)}}$ are homogeneous, with $|c_{_{(1)}}|\,|c_{_{(2)}}|=|c|$.

A \emph{bialgebra} is a vector space $B$ that is both an algebra and a coalgebra, such that its coproduct and counit are unital algebra homomorphisms; that is, for all homogeneous elements $b,c\in B$, it holds\[\Delta(bc)=\phi(|b_{_{(2)}}|,|c_{_{(1)}}|)\,b_{_{(1)}}c_{_{(1)}}\otimes b_{_{(2)}}c_{_{(2)}},\quad\Delta(1)=1\otimes1,\quad\varepsilon(bc)=\varepsilon(b)\varepsilon(c),\quad\varepsilon(1)=1.\]

A bialgebra $H$ is a \emph{Hopf algebra} if it is equipped with a linear map $S:H\to H$, called the \emph{antipode}, such that for every homogeneous element $h\in H$, it holds\begin{gather}\label{011}h_{_{(1)}}S(h_{_{(2)}})=\varepsilon(h)1=S(h_{_{(1)}})h_{_{(2)}}.\end{gather}

Algebraic structures equipped with a $\G$-grading and a bicharacter $\phi$ are usually termed \emph{color algebraic structures} \cite[Section 10.5]{Mont93}, including \emph{color algebras}, \emph{color coalgebras}, \emph{color bialgebras}, and \emph{color Hopf algebras}~\cite{ChePetOys06,AndAngBag14,Sci24}. In particular, when $\G=\{1\}$ and $\phi$ is identically $1$, these recover the classical algebraic structures. For the sake of brevity, and since the bicharacter $\phi$ remains fixed throughout our discussion, we shall refer to them simply as \emph{bialgebras} and \emph{Hopf algebras}, implicitly assuming the underlying $\G$-grading.

\subsubsection{Connected bialgebras}

A bialgebra $B$ is said to be \emph{connected} if it is equipped with an additional grading $\{B_n\mid n\in\N\}$, that is compatible with the $\G$-grading, such that all its operations are graded maps:\[\ubar{(B_n)}_x=B_n\cap\ubar{B}_x,\quad B_m\,B_n\subseteq B_{m+n},\quad B_0=\field,\quad\,\,\Delta(B_n)\subseteq\bigoplus_{r+s=n}B_r\otimes B_s,\quad\,\,\varepsilon\biggl(\bigoplus_{\,n\geq1}B_n\biggr)=\{0\}.\]

Consequently, the family of subspaces $\{B_{n,x}:=\ubar{(B_n)}_x\mid(n,x)\in\N\times\G\}$ defines a $(\N\times\G)$-grading on $B$, which we shall hereafter refer to as the \emph{bigrading} of the algebra, which satisfies the following properties:\begin{equation}\label{013}
\!\!B_{m,x}\,B_{n,y}\subseteq B_{m+n,xy},\quad\!B_{0,1}=\field,\quad\Delta(B_{n,x})\subseteq\bigoplus_{\substack{r+s=n\\yz=x}}B_{r,y}\otimes B_{s,z},\quad\varepsilon(B_{n,x})=\{0\}\,\,\text{ if }\,\,(n,x)\neq(0,1).\end{equation}

As shown in \cite[Exercise~1.3.20]{GrRe20}, if $B$ is a connected bialgebra, then every $g\in B_n$ with $n\geq1$ satisfies $\Delta(g)=1\otimes g+\bar{\Delta}(g)+g\otimes1$, where $\bar{\Delta}(g)\in\bigoplus_{k=1}^{n-1}B_k\otimes B_{n-k}$. This property ensures that $B$ can be equipped with a unique antipode satisfying the conditions in \eqref{011} \cite[Proposition~1.4.16]{GrRe20}. Specifically, the antipode is defined by $S(1)=1$ and, for $g\in B_n$, by Takeuchi's recursive formula~\cite{Ta70}:\begin{gather}\label{010}S(g)=-g-S(g_{_{(1)}})g_{_{(2)}},\quad\text{where}\quad\bar{\Delta}(g)=g_{_{(1)}}\otimes g_{_{(2)}}.\end{gather}

\begin{pro}\label{105}
The antipode defined in \eqref{010} is $\G$-graded, granting $B$ the structure of a Hopf algebra.
\end{pro}
\begin{proof}
Let $g\in B_{n,x}$ for some $x\in\G$ and $n\in\N$. We proceed by induction on the degree $n$. If $n=0$, then $g\in B_0\cap\ubar{B}_x$. Since $B$ is connected, then $g=\alpha1$ for some $\alpha\in\field$. As we established that $S(1)=1$, it follows that $S(g)=\alpha S(1)=\alpha1=g$. Thus, $S(g)\in B_{0,x}$. If $n\geq1$, assume that $S(B_{m,y})\subseteq\ubar{B}_y$ for all $m<n$ and $y\in\G$. By the recursive formula \eqref{010}, we have $S(g)=-g-S(g_{_{(1)}})g_{_{(2)}}$. By \eqref{013}, each term in the reduced coproduct $\bar{\Delta}(g)=g_{_{(1)}}\otimes g_{_{(2)}}$ consists of homogeneous elements $g_{_{(1)}}\in B_{r,y}$ and $g_{_{(2)}}\in B_{s,z}$, where $r+s=n$ and $yz=x$. Since $r<n$, the inductive hypothesis implies that $S(g_{_{(1)}})\in\ubar{B}_y$. Furthermore, because the product in $B$ is $\G$-graded, we have $S(g_{_{(1)}})g_{_{(2)}}\in\ubar{B}_y\ubar{B}_z\subseteq\ubar{B}_{yz}=\ubar{B}_x$. Finally, since $g\in\ubar{B}_x$ and every term in the summation also belongs to $\ubar{B}_x$, their sum $S(g)$ must be in $\ubar{B}_x$. This completes the induction.
\end{proof}

A \emph{graded bi-ideal} of a connected bialgebra $B$ is a graded subspace $I=\bigoplus_{n\geq0}I_n$ satisfying\[I_{n,x}:=\ubar{(I_n)}_x=I\cap B_{n,x},\quad I_0=\{0\},\quad IB\subseteq I,\quad BI\subseteq I,\quad\Delta(I)\subseteq I\otimes B+B\otimes I,\quad\varepsilon(I)=\{0\}.\]Note that the condition $S(I)\subseteq I$ is automatically satisfied for the antipode $S$ defined in \eqref{010}. The quotient space $B/I$ inherits a natural structure of a connected bialgebra, whose homogeneous components are given by $\ubar{(B/I)}_x=\ubar{B}_x/\ubar{I}_x$, $B_n/I_n$ and $B_{n,x}/I_{n,x}$, and whose operations are induced by those of $B$ \cite[Proposition 4.2.13]{DaNaRa01} \cite[Remark~1]{ChePetOys06} \cite[Subsection~3.2]{Sci24}; specifically,\begin{gather*}
(g+I)(h+I)=gh+I,\qquad\tilde{u}(\alpha)=u(\alpha)+I,\\\tilde{\Delta}(g+I)=\Delta(g)+(I\otimes B+B\otimes I),\quad\tilde{\varepsilon}(g+I)=\varepsilon(g),\quad\tilde{S}(g+I)=S(g)+I.
\end{gather*}

A \emph{bigraded sub-bialgebra} of a connected bialgebra $B$ is a graded subspace $A=\bigoplus_{n\geq0}A_n$ satisfying $A_{n,x}:=\ubar{(A_n)}_x=A\cap B_{n,x}$, $1\in A$, $A^2\subseteq A$, and $\Delta(A)\subseteq A\otimes A$. Note that any bigraded sub-bialgebra is itself a connected bialgebra. Given connected bialgebras $B$ and $B'$ with coproducts $\Delta,\Delta'$ and counits $\varepsilon,\varepsilon'$, respectively, a \emph{connected bialgebra homomorphism} is a bigraded algebra homomorphism $\phi:B\to B'$ satisfying $(\phi\otimes\phi)\circ\Delta=\Delta'\circ\phi$ and $\varepsilon'\circ\phi=\varepsilon$. An \emph{isomorphism} is a bijective bialgebra homomorphism.

\begin{pro}\label{165}
Let $\psi:B\to B'$ be a connected bialgebra homomorphism. Then, the kernel $\ker(\psi)=\psi^{-1}(\{0\})$ is a graded bi-ideal of $B$, and $\psi(B)$ is a connected bialgebra that is isomorphic to $B/\ker(\psi)$.
\end{pro}
\begin{proof}
By the first isomorphism theorem for algebras, $\ker(\psi)$ is an algebra ideal of $B$, $\psi(B)$ is a subalgebra of $B'$, and the induced map $\bar{\psi}:B/\ker(\psi)\to\psi(B)$, defined by $\bar{\psi}(g+\ker(\psi))=\psi(g)$, is an algebra isomorphism.

Let $g\in\ker(\psi)$ and let $g_{n,x}$ be its projection onto the homogeneous component $B_{n,x}$. Since $\psi(g)=0$ and each $\psi(g_{n,x})\in B'_{n,x}$, the linear independence of homogeneous components implies $\psi(g_{n,x})=0$; hence $g_{n,x}\in\ker(\psi)$, showing that $\ker(\psi)$ is a bigraded subspace with $\ker(\psi)_{n,x}=\ker(\psi)\cap B_{n,x}$. In degree zero, since $B$ is connected, we have $B_0=\field 1$. Because $\psi(1)=1\neq0$, it follows that $\ker(\psi)_0=\ker(\psi)\cap B_0=\{0\}$. Consequently, $(B/\ker(\psi))_0\simeq B_0\simeq\field$, ensuring that the quotient inherits the connectedness condition.

Regarding the coalgebra structure, since $\psi$ preserves the coproduct and counit, it is a well-known property that $\Delta(\ker(\psi))\subseteq\ker(\psi)\otimes B+B\otimes\ker(\psi)$ and $\varepsilon(\ker(\psi))=\{0\}$ \cite[Proposition~1.4.9]{DaNaRa01}. Therefore, $\ker(\psi)$ is a graded bi-ideal of $B$, ensuring that the quotient $B/\ker(\psi)$ inherits a connected bialgebra structure. Similarly, $\psi(B)$ satisfies $\Delta'(\psi(B))\subseteq\psi(B)\otimes\psi(B)$ and $\varepsilon'(\psi(B))\subseteq\field$ \cite[Proposition~1.4.9]{DaNaRa01}, and by the fundamental isomorphism theorem for coalgebras \cite[Theorem~1.4.10 and Corollary~1.4.11]{DaNaRa01}, the algebra isomorphism $\bar{\psi}$ naturally preserves the coproduct and counit, that is, $(\bar{\psi}\otimes\bar{\psi})\circ\tilde{\Delta}=\Delta'\circ\bar{\psi}$ and $\varepsilon'\circ\bar{\psi}=\tilde{\varepsilon}$.

Finally, since $\psi$ is a bigraded map, $\psi(B)$ inherits the bigrading $\psi(B)_{n,x}=\psi(B_{n,x})$. In degree zero, we have $\psi(B)_0=\psi(B_0)=\psi(\field 1)=\field 1$, confirming that $\psi(B)$ is a connected bialgebra. Since for any homogeneous element $g\in B_{n,x}$, its canonical class in the quotient satisfies $\bar{\psi}(g+\ker(\psi))=\psi(g)\in\psi(B)_{n,x}$, the isomorphism $\bar{\psi}$ is bigraded, establishing that $\bar{\psi}$ is an isomorphism of connected bialgebras.
\end{proof}

\section{Invariants}\label{167}

Recall that if $G$ is a group acting on a vector space $V$ by linear maps via a group homomorphism $\rho:G\to\GL(V)$, an \emph{invariant} with respect to this action, or simply a $G$-\emph{invariant}, is an element $g\in V$ fixed by every group element; that is, $\rho_\sigma(g)=g$ for all $\sigma\in G$. The $G$-invariants form a vector subspace of $V$ denoted by $V^G$. If $V$ is an algebra and $G$ acts by algebra homomorphisms, then $V^G$ is a subalgebra of $V$. If $H$ is a subgroup of $G$, then, by restricting the action, each $G$-invariant is also $H$-invariant; hence $V^G\subseteq V^H$.

In what follows, we are interested in the study of invariants within the framework of formal power series. Specifically, we focus on the case where the vector space $V$ consists of functions from a monoid to a field, and the group $G$ acts on these functions by composing them with a given action on the monoid, possibly modified by a scalar weight. To this end, we begin by establishing the necessary structure for these spaces.\medskip

Throughout this section, $M$ denotes a monoid, whose elements are referred to as \emph{simple monomials}.

\subsection{Power series}\label{007}

A \emph{formal power series} over $M$ is a map $g:M\to\field$ \cite[Section 1.3]{DrKu09}. The elements in the range of $g$ are referred to as its \emph{coefficients}. A series as above is typically expressed as a generalized linear combination:\[g=\sum_{u\in M}g(u)u,\]where each summand is called a \emph{monomial} of $g$. The \emph{support} of $g$, denoted by $\supp(g)$, is the set of elements in $M$ with nonzero coefficients. A series with finite support is called a \emph{polynomial}. A simple monomial $u$ is said to \emph{occur} in $g$ if $u\in\supp(g)$. More generally, since $g(u)=0$ for all $u\notin\supp(g)$, we may write\[g=\sum_{u\in S}g(u)u,\]where $S\subseteq M$ is any subset containing the support of $g$. If $\supp(g)=\{u\}$, we abbreviate the series as $g=g(u)u$, and if $g(u)=1$, we simply write $g=u$. Conversely, if the support of $g$ is empty, we write $g=0$.

The collection of all formal power series over $M$ is denoted by $\field^M$. This set is naturally endowed with a binary operation, \emph{addition}, and a scalar action of $\field$, which grants it the structure of a vector space:\[g+h=\sum_{u\in M}\left(\,g(u)+h(u)\,\right)\,u\,;\qquad ag=\sum_{u\in M}\,\left(\,ag(u)\,\right)\,u.\]The set of polynomials over $M$, denoted by $\pols{\field}{M}$, constitutes a vector subspace of $\field^M$; indeed, this follows from the property that for any $g,h\in\field^M$ and any scalar $a\in\field$, the support satisfies\begin{equation}\label{000}\supp(g+h)\subseteq\supp(g)\cup\supp(h)\quad\text{and}\quad\supp(ag)\subseteq\supp(g).\end{equation}

Now, we establish the condition under which the product of the monoid $M$ can be extended to the space of all formal power series, thereby endowing $\field^M$ with an algebra structure.

A \emph{twisting} on $M$ \cite[Definition 3]{Wi07} is a map $\tau:M\times M\to\field$ that satisfies the following 2-cocycle conditions (cf. \cite[Section 2.1]{KarV293}):\begin{equation}\label{001}\tau(u,v)\,\tau(uv,w)=\tau(u,vw)\,\tau(v,w)\quad\text{and}\quad\tau(u,1)=1=\tau(1,u)\quad\text{for all}\quad u,v,w\in M.\end{equation}We say that $\tau$ is \emph{trivial} if $\tau(u,v)=1$ for all $u,v\in M$.

Given another monoid $N$ equipped with a twisting $\gamma:N\times N\to\field$, a monoid homomorphism $\psi:M\to N$ is said to be \emph{twisting-preserving} if it satisfies $\gamma(\psi(u),\psi(v))=\tau(u,v)$ for all $u,v\in M$.\medskip

Throughout the remainder of the section we assume the following condition:
\newcounter{condition}
\begin{enumerate}[start=\thecondition,label={(\condition)},ref={\condition}]
\item\label{015}\conditions{0}
\end{enumerate}

\begin{pro}\label{014}
Assume that $M$ is \emph{locally finite with respect to $\tau$}; that is, for each $u\in M$, the set of pairs $(v,w)\in M\times M$ such that $vw=u$ and $\tau(v,w)\neq0$ is finite. Then, the \emph{twisted convolution product} defined by\begin{equation}\label{023}(g\cdot h)(u)=\sum_{vw=u}\tau(v,w)g(v)h(w),\quad\text{or equivalently,}\quad g\cdot h=\sum_{u\in M}\Biggl(\,\sum_{vw=u}\tau(v,w)g(v)h(w)\Biggr)u,\end{equation}is well-defined for all $g,h\in\field^M$. Furthermore, this operation endows $\field^M$ with the structure of an associative algebra over $\field$, with identity element $1_{_M}$.
\end{pro}
\begin{proof}
The local finiteness condition ensures that the sum defining $(g\cdot h)(u)$ contains only finitely many non-zero terms for each $u\in M$; thus, the product is well-defined as a series from $M$ to $\field$. The associativity of this product follows directly from the first condition in \eqref{001}, while the second condition in \eqref{001} ensures that the identity element of the monoid $M$ acts as the unit of the algebra.
\end{proof}
When the twisting $\tau$ is understood from the context, we will often omit the phrase ``with respect to $\tau$'' and simply say that $M$ is locally finite. 

Next, we show that if the monoid $M$ is equipped with a grading, this structure naturally induces a direct product decomposition on the space of all formal power series.

We say that $M$ is \emph{$\P$-graded with respect to $\tau$}, for some commutative monoid $\P$, if it is equipped with a map $\nu:M\to\P$, called a \emph{degree map}, such that $\nu(1)=1$ and $\nu(uv)=\nu(u)\nu(v)$ whenever $\tau(u,v)\neq0$. For each $u\in M$, $\nu(u)$ is referred to as the \emph{$\P$-degree} of $u$. For any $x\in\P$, the fibre $\nu^{-1}(\{x\})$ is denoted by $M_x$, and the collection $\{M_x\mid x\in\P\}$ is called a \emph{$\P$-grading}. By definition, these fibres satisfy $uv\in M_{xy}$ for all $u\in M_x$ and $v\in M_y$ such that $\tau(u,v)\neq0$, and the monoid admits a decomposition:\[M=\bigsqcup_{x\in\P}M_x.\]We say that $M$ is \emph{graded with respect to $\tau$} if it is $\N$-graded by a degree map $\nu:M\to\N$ satisfying the condition $\nu^{-1}(\{0\})=\{1\}$. For this specific class of monoids, the prefix ``$\N$-'' is typically omitted from the related terminology.

When the twisting $\tau$ is clear from the context, we will often omit the qualifying phrase ``with respect to $\tau$'' and simply refer to $M$ as being $\P$-graded or graded. Furthermore, throughout this paper, any $\G$-grading on a monoid is always considered with respect to the trivial twisting, meaning that the degree map is a monoid homomorphism.

If $M$ and $N$ are $\P$-graded monoids, a map $\psi:M\to N$ is said to be \emph{$\P$-graded} if $\psi(M_x)\subseteq N_x$ for all $x\in\P$.

\begin{pro}\label{005}
If $M$ is a $\P$-graded monoid with $\P$-grading $\{M_x\mid x\in\P\}$ and degree map $\nu:M\to\P$, then the vector space $\field^M$ admits a direct product decomposition\[\field^M=\prod_{x\in\P}\field_x^M,\]where each $\field_x^M$ is the subspace consisting of all series $g\in\field^M$ such that $\supp(g)\subseteq M_x$. Consequently, the subspace $\pols{\field}{M}$ is $\P$-graded with homogeneous components $\pols{\field}{M_x}:=\{g\in\pols{\field}{M}\mid\supp(g)\subseteq M_x\}$. 

Moreover, if $M$ is locally finite, the twisted convolution product respects the direct product decomposition; that is, $\field_x^M\cdot\field_y^M\subseteq\field_{xy}^M$ for all $x,y\in\P$.
\end{pro}
\begin{proof}
It follows from \eqref{000} that each $\field_x^M$ is a vector subspace of $\field^M$. Furthermore, for any distinct $x,y\in\P$, the disjointness of fibres implies $\field_x^M\cap\field_y^M=\{0\}$. Since $\{M_x\mid x\in\P\}$ is a partition of $M$, any series $g\in\field^M$ admits a unique decomposition\begin{equation}\label{004}g=\sum_{x\in\P}g_x,\quad\text{where}\quad g_x=\sum_{u\in M_x}g(u)u.\end{equation}Since each $g_x$ is clearly contained in $\field_x^M$, we obtain the required direct product decomposition. Finally, if $g\in\pols{\field}{M}$, then $\supp(g)$ is finite. This implies that $g_x=0$ for all but finitely many $x\in\P$, so the decomposition lies in the direct sum. Therefore, $\pols{\field}{M}$ is $\P$-graded.

Given $g\in\field_x^M$ and $h\in\field_y^M$, for any $u\in\supp(g\cdot h)$, there must exist $v\in\supp(g)$ and $w\in\supp(h)$ such that $vw=u$ and $\tau(v,w)\neq 0$. This implies $\nu(v)=x$ and $\nu(w)=y$. Since $\tau(v,w)\neq 0$, we have $\nu(u)=\nu(v)\nu(w)=xy$, which implies that $\supp(g\cdot h)\subseteq M_{xy}$. This confirms that $g\cdot h\in\field_{xy}^M$ and, consequently, $\field_x^M\cdot\field_y^M\subseteq\field_{xy}^M$.
\end{proof}

\begin{crl}\label{006}
Assume that $M$ is locally finite and equipped with a $\P$-grading $\{M_x\mid x\in\P\}$. Let $A$ be a subalgebra of $\field^M$ that admits the direct sum decomposition\[A=\bigoplus_{x\in\P}\left(A\cap\field_x^M\right).\]Then $A$ is a $\P$-graded algebra with homogeneous components $A\cap\field_x^M=\{g\in A\mid\supp(g)\subseteq M_x\}$.
\end{crl}
\begin{proof}
By hypothesis, $A$ is a $\P$-graded vector space with homogeneous components $A\cap\field_x^M$. To verify the compatibility of the product with the grading, let $g\in A\cap\field_x^M$ and $h\in A\cap\field_y^M$. Since $A$ is a subalgebra, we have $g\cdot h\in A$. On the other hand, Proposition~\ref{005} implies that $g\cdot h\in\field_{xy}^M$. Consequently, $g\cdot h\in A\cap\field_{xy}^M$, which confirms that $A$ is a $\P$-graded algebra.
\end{proof}

In the remainder of this section, we assume the following condition:
\stepcounter{condition}
\begin{enumerate}[start=\thecondition,label={(\condition)},ref={\condition}]
\item\label{016}\conditions{1}
\end{enumerate}
For each simple monomial $u\in M$, we denote its $\G$-degree by $|u|$, so that $\ubar{M}_x=\{u\in M\mid|u|=x\}$ for each $x\in\G$. Consequently, $M$ admits a decomposition into fibres:\[M=\bigsqcup_{x\in\G}\ubar{M}_x.\]Under condition \eqref{016}, Proposition~\ref{014} and Proposition~\ref{005} imply that $\field^M$ is an algebra admitting a direct product decomposition into subspaces $\ubar{\field}_x^M:=\{g\in\field^M\mid\supp(g)\subseteq\ubar{M}_x\}$. We refer to $\field^M$, equipped with this structure, as the \emph{twisted algebra of formal power series} on $M$. For $g\in\field^M\setminus\{0\}$, we write $|g|=x$ whenever $\supp(g)\subseteq\ubar{M}_x$, thereby extending the degree map from $M$ to the non-zero elements in these subspaces. Corollary~\ref{006} implies that $\pols{\field}{M}$ is a $\G$-graded subalgebra called the \emph{twisted algebra of polynomials}, also known as the \emph{twisted monoid algebra} \cite[Section 3]{Wi07}. The homogeneous components of $\pols{\field}{M}$ are denoted by $\pols{\field}{\ubar{M}_x}$. Note that $\pols{\field}{\ubar{M}_x}$  is precisely the vector subspace spanned by the set $\ubar{M}_x$. In particular, if $\tau$ is identically equal to $1$, then $\pols{\field}{M}$ coincides with the standard monoid algebra of $M$.

\subsection{Functions}\label{188}

Throughout the remainder of the section, we also assume the following condition:
\stepcounter{condition}
\begin{enumerate}[start=\thecondition,label={(\condition)},ref={\condition}]
\item\label{017}\conditions{2}
\end{enumerate}
Under condition \eqref{017}, Proposition~\ref{014} and Proposition~\ref{005} imply that $\field^M$ admits a direct product decomposition into subspaces $\{g\in\field^M\mid\supp(g)\subseteq M^d\}$. For $g\in\field^M\setminus\{0\}$, we write $\deg(g)=d$ whenever $\supp(g)\subseteq M^d$, thereby extending the degree map from $M$ to the non-zero elements in these subspaces. Corollary~\ref{006} implies that $\pols{\field}{M}$ defines a graded subalgebra, whose homogeneous components are denoted by $\pols{\field}{M^d}$.

A series $g\in\field^M$ is said to be \emph{of bounded degree} if there exists $m\in\N$ such that $\deg(u)\leq m$ for all $u\in\supp(g)$. We denote the set of all such series by $\funcs{\field}{M}$. By Proposition~\ref{005}, the twisted
convolution product respects the grading, so $\funcs{\field}{M}$ defines a subalgebra of $\field^M$ that contains $\pols{\field}{M}$. We call $\funcs{\field}{M}$ the \emph{twisted algebra of formal power series of bounded degree} on $M$, and we shall refer to its elements simply as \emph{functions}. By definition, the support of any function is contained in a union of finitely many fibres $M^d$. Consequently, by Corollary~\ref{006}, $\funcs{\field}{M}$ constitutes a graded subalgebra, whose homogeneous components are denoted by $\funcs{\field}{M^d}$.

Since $M$ is simultaneously $\G$-graded and $\N$-graded, it carries a natural $(\G\times\N)$-grading with respect to $\tau$ given by the fibres $M^d_x:=M^d\cap\ubar{M}_x$, which we call the \emph{bigrading}. Consequently, $M$ admits the following decomposition:\[M=\bigsqcup_{(d,x)\in\N\times\G}M_x^d.\]Furthermore, condition \eqref{017} implies that the support of any function is also contained in a union of finitely many $\G$-fibres $\ubar{M}_x$. Consequently, by Corollary~\ref{006}, $\funcs{\field}{M}$ is both a $\G$-graded subalgebra, whose homogeneous components are denoted by $\funcs{\field}{\ubar{M}_x}$, and a $(\G\times\N)$-graded subalgebra, whose homogeneous components are denoted by $\funcs{\field}{M_x^d}$. In summary, we obtain the following decompositions:\[\funcs{\field}{M}\,\,\,=\,\,\,\bigoplus_{d\in\N}\funcs{\field}{M^d}\,\,\,=\,\,\,\bigoplus_{x\in\G}\funcs{\field}{\ubar{M}_x}\,\,\,=\bigoplus_{(x,d)\in\G\times\N}\funcs{\field}{M_x^d}.\]

For $g\in\funcs{\field}{M}\setminus\{0\}$, we write $\bdg(g)=(x,d)$ if and only if $g\in\funcs{\field}{M_x^d}$; that is, $\bdg(g)=(|g|,\deg(g))$. We call $\bdg(g)$ the \emph{bidegree} of $g$. Note that $M$ also inherits this terminology; for any $u\in M$, we have $\bdg(u)=(|u|,\deg(u))$. The $(\G\times\N)$-grading of any subalgebra of functions is referred to as its \emph{bigrading}. Accordingly, a $(\G\times\N)$-graded map $\Psi$ between two such subalgebras, meaning $\Psi(\funcs{\field}{M_x^d})\subseteq\funcs{\field}{N_x^d}$ for all $(x,d)\in\G\times\N$, is also called \emph{bigraded}. Such a map is necessarily $\G$-graded and $\N$-graded.

We conclude this subsection by establishing the conditions under which a monoid homomorphism can be extended to an algebra homomorphism between spaces of functions. In what follows, $N$ denotes a monoid satisfying the same structural assumptions as $M$. Specifically, $N$ is locally finite with respect to a twisting $\gamma:N\times N\to\field$ and is endowed with a $\G$-grading $\{\ubar{N}_x\mid x\in\G\}$ for the trivial twisting. Furthermore, $N$ is equipped with a grading $\{N^d\mid d\in\N\}$ with respect to $\gamma$, where each $N^d$ intersects only finitely many $\G$-fibres $\ubar{N}_x$.

A \emph{bigraded homomorphism} from $M$ to $N$ is a $(\G\times\N)$-graded twisting-preserving monoid homomorphism $\psi:M\to N$ such that $\psi^{-1}(\{u\})$ is finite for all $u\in N$. The condition $\psi(M_x^d)\subseteq N_x^d$ implies that $\psi$ preserves both gradings; specifically, $|\psi(u)|=|u|$, $\deg(\psi(u))=\deg(u)$, and $\bdg(\psi(u))=\bdg(u)$ for all $u \in M$. If $\psi$ is a bijection, the condition on the fibres is automatically satisfied and we call it a \emph{bigraded isomorphism}.

\begin{pro}\label{019}
Let $\psi:M\to N$ be a bigraded homomorphism. Then:
\begin{enumerate}
\item\label{020}There exists an algebra homomorphism $\Psi:\funcs{\field}{M}\to\funcs{\field}{N}$ defined for each $g\in\funcs{\field}{M}$ by\[\Psi(g)=\sum_{u\in N}\Biggl(\,\sum_{\psi(v)=u}g(v)\Biggr)u.\]
\item\label{021}The homomorphism $\Psi$ is bigraded, that is, $\Psi(\funcs{\field}{M_x^d})\subseteq\funcs{\field}{N_x^d}$ for all $x\in\G$ and $d\in\N$.
\item\label{022}If $\psi$ is surjective, then $\Psi$ is surjective. If $\psi$ is bijective, then $\Psi$ is an isomorphism with $\Psi(g)=g\circ\psi^{-1}$.
\end{enumerate}
\end{pro}
\begin{proof}
\eqref{020} Since $\psi$ is a bigraded homomorphism, each fibre $\psi^{-1}(\{u\})$ is finite, ensuring the inner sum in $\Psi(g)$ is well-defined. Moreover, since $\supp(\Psi(g))\subseteq\psi(\supp(g))$ and $\psi$ preserves the degree, $\Psi$ maps $\funcs{\field}{M}$ to $\funcs{\field}{N}$. The linearity of $\Psi$ follows directly from the linearity of the finite sums. To prove that $\Psi$ is multiplicative, let $g,h \in \funcs{\field}{M}$ and $u \in N$. Since $\psi$ is twisting-preserving, we have\[\underbrace{
	\sum_{\psi(v)=u}\Biggl(\,\sum_{ab=v}\tau(a,b)g(a)h(b)\!\Biggr)
}_{\Psi(g\cdot h)(u)}
\,\,\,=\,\,\,
\sum_{\psi(a)\psi(b)=u}\tau(a,b)g(a)h(b)
\,\,\,=\,\,\,
\underbrace{
	\sum_{cd=u}\gamma(c,d)\Biggl(\,\sum_{\psi(a)=c}g(a)\!\Biggr)\!\Biggl(\,\sum_{\psi(b)=d}h(b)\!\Biggr)
}_{\left(\Psi(g)\Psi(h)\right)(u)}.\]Hence, $\Psi$ is an algebra homomorphism.

\eqref{021} Let $g\in\funcs{\field}{M_x^d}$ and $u\in\supp(\Psi(g))$. By definition, there must exist $v\in\supp(g)\subseteq M_x^d$ such that $\psi(v)=u$. Since $\psi$ is a bigraded homomorphism, it preserves the bidegree, and thus $\bdg(u)=\bdg(\psi(v))=\bdg(v)=(x,d)$. Consequently, $\supp(\Psi(g))\subseteq N_x^d$, proving that $\Psi(\funcs{\field}{M_x^d})\subseteq\funcs{\field}{N_x^d}$.

\eqref{022} Suppose $\psi$ is surjective, and for each $u\in N$, choose a preimage $v_u\in\psi^{-1}(\{u\})$. For any $h\in\funcs{\field}{N}$, define $g:M\to\field$ by $g(v_u)=h(u)$ for each $u\in\supp(h)$, and $g(v)=0$ otherwise. Then $\Psi(g)(u)=g(v_u)=h(u)$. Since $\psi$ preserves degree, $g$ is of bounded degree whenever $h$ is; thus $\Psi$ is surjective. 

Finally, if $\psi$ is bijective, then for each $u\in N$, the fibre $\psi^{-1}(\{u\})$ contains a unique element $\psi^{-1}(u)$. The definition of $\Psi$ then simplifies to $\Psi(g)(u)=g(\psi^{-1}(u))=(g\circ\psi^{-1})(u)$. Thus $\Psi(g)=g\circ\psi^{-1}$, which is clearly bijective with inverse $\Psi^{-1}(h)=h\circ\psi$.
\end{proof}

\subsection{Invariant functions}\label{189}

In this subsection, we formalize the notion of invariant formal power series and study the conditions under which they form a well-behaved algebra.

Given a group action $\rho:G\to\sym(M)$, for $u\in M$, we denote its orbit by $Gu=\{\,\rho_\sigma(u)\mid\sigma\in G\,\}$ and its stabilizer by $G_{(u)}=\{\,\sigma\in G\mid\rho_\sigma(u)=u\,\}$. These orbits form the equivalence classes of the relation $R_G$ on $M$ consisting of the pairs $(u,v)\in M\times M$ such that $v=\rho_\sigma(u)$ for some $\sigma\in G$. If $M$ is $\P$-graded, we say that $\rho$ is \emph{$\P$-graded} if each bijection $\rho_\sigma$ is $\P$-graded. We say that $\rho$ is \emph{bigraded} if it is $(\G\times\N)$-graded; that is, $\rho_\sigma(M_x^d)=M_x^d$, or equivalently, $\bdg(\rho_\sigma(u))=\bdg(u)$ for all $\sigma\in G$ and $u\in M$. Note that such an action is both $\G$-graded and $\N$-graded. We say that $\rho$ is \emph{multiplicative with respect to $\tau$} if $\rho_\sigma(uv)=\rho_\sigma(u)\rho_\sigma(v)$ for all $\sigma\in G$ and $u,v\in M$ with $\tau(u,v)\neq0$. More generally, $\rho$ is called \emph{quasi-multiplicative with respect to $\tau$} if for each $\sigma\in G$ and $w\in M$, there exist $\alpha,\varsigma\in G$ such that $\rho_\sigma(uv)=\rho_\alpha(u)\rho_\varsigma(v)$ for all $u,v\in M$ satisfying $uv=w$ and $\tau(u,v)\neq0$.

A \emph{weight function} for the action $\rho:G\to\sym(M)$ is a map $\omega:G\times M\to\field$ that satisfies the following 1-cocycle conditions (cf. \cite[Subsection~8.1]{Aar97}, \cite[Section 3.1]{KarV293}):\begin{gather}\label{003}\omega(\sigma\varsigma,u)=\omega(\sigma,\rho_\varsigma(u))\,\omega(\varsigma,u)\quad\text{and}\quad\omega(1,u)=1=\omega(\sigma,1)\quad\text{for all}\quad\sigma,\varsigma\in G\quad\text{and}\quad u\in M.\end{gather}It follows from \eqref{003} that $\omega$ takes values in $\field^\times$, with the inverse given by $\omega(\sigma,u)^{-1}=\omega(\sigma^{-1},\rho_\sigma(u))$. If $\rho$ is quasi-multiplicative with respect to $\tau$, we say that $(\tau,\omega)$ is a \emph{multiplicative pair} for the action $\rho$ if for all $\sigma\in G$ and $w\in M$, there exist $\alpha,\varsigma\in G$ with $\rho_\sigma(uv)=\rho_\alpha(u)\rho_\varsigma(v)$ such that\begin{equation}\label{028}\tau(u,v)\,\omega(\sigma,uv)=\tau(\rho_\alpha(u),\rho_\varsigma(v))\,\omega(\alpha,u)\,\omega(\varsigma,
v),\end{equation}for all $u,v\in M$ with $uv=w$ and $\tau(u,v)\neq0$; whenever $\rho$ is multiplicative with respect to $\tau$, we require~\eqref{028} to hold with $\alpha=\varsigma=\sigma$.

\begin{pro}\label{029}
Let $\rho:G\to\sym(M)$ be a $\G$-graded group action together with a weight function $\omega:G\times M\to\field^\times$. Then:
\begin{enumerate}
\item\label{030}The map $\bm{\rho}:G\to\GL(\field^M)$ defined by $\bm{\rho}_\sigma(g)(u)=\omega(\sigma,\rho_{\sigma^{-1}}(u))\,g(\rho_{\sigma^{-1}}(u))$ constitutes a group action on the space $\field^M$ such that $\supp(\bm{\rho}_\sigma(g))=\rho_\sigma(\supp(g))$ and $\bm{\rho}_\sigma(\ubar{\field}_x^M)=\ubar{\field}_x^M$ for all $\sigma\in G$ and $x\in\G$.
\item\label{027}If $\rho$ is quasi-multiplicative, then $\bm{\rho}$ is quasi-multiplicative on the non-zero products of monomials of $\field^M$ if and only if $(\tau,\omega)$ is a multiplicative pair for $\rho$.
\item\label{181}If $\rho$ is quasi-multiplicative and $(\tau,\omega)$ is a multiplicative pair, for every $u\in M$ the map $(a,b)\mapsto(\rho_\alpha(a),\rho_\varsigma(b))$ is a bijection between the sets of factorizations with non-zero twisting of $\rho_{\sigma^{-1}}(u)$ and $u$, where $\alpha,\varsigma\in G$ are the elements associated with $(\sigma,\rho_{\sigma^{-1}}(u))$ by the quasi-multiplicativity of $\rho$.
\item\label{031}If $\rho$ is multiplicative and $(\tau,\omega)$ is a multiplicative pair, then $\bm{\rho}$ is also multiplicative.
\item\label{032}The restriction $\sigma\mapsto\bm{\rho}_\sigma|_{\pols{\field}{M}}$ defines a $\G$-graded group action on the vector space $\pols{\field}{M}$.
\item\label{033}If $\rho$ is a bigraded action, the restriction $\sigma\mapsto\bm{\rho}_\sigma|_{\funcs{\field}{M}}$ defines a bigraded group action on the vector space $\funcs{\field}{M}$, that is, $\bm{\rho}_\sigma(\funcs{\field}{M_x^d})=\funcs{\field}{M_x^d}$ for all $\sigma\in G$.
\item\label{034}If $\rho$ is a quasi-multiplicative bigraded action and $(\tau,\omega)$ is a multiplicative pair, then $\funcs{\field}{M}^G$ is a bigraded algebra that contains $\pols{\field}{M}^G$ as a bigraded subalgebra.
\end{enumerate}
\end{pro}
\begin{proof}
\eqref{030}~Since $\rho$ is a group action and $\omega$ satisfies \eqref{003}, we have $\bm{\rho}_1(g)(u)=\omega(1,\rho_{1^{-1}}(u))g(\rho_{1^{-1}}(u))=g(u)$, and for any $\sigma,\varsigma\in G$,\[\begin{array}{rcl}
\bm{\rho}_\sigma(\bm{\rho}_\varsigma(g))(u)
&=&\omega(\sigma,\rho_{\sigma^{-1}}(u))\,\bm{\rho}_\varsigma(g)(\rho_{\sigma^{-1}}(u))\\
&=&\omega(\sigma,\rho_{\sigma^{-1}}(u))\,\omega(\varsigma,\rho_{\varsigma^{-1}}(\rho_{\sigma^{-1}}(u)))\,g(\rho_{\varsigma^{-1}}(\rho_{\sigma^{-1}}(u)))\\
%&=&\omega(\sigma,\rho_{\sigma^{-1}}(u))\,\omega(\varsigma,\rho_{(\sigma\varsigma)^{-1}}(u))\,g(\rho_{(\sigma\varsigma)^{-1}}(u))\\
&=&\omega(\sigma,\rho_\varsigma(\rho_{(\sigma\varsigma)^{-1}}(u)))\,\omega(\varsigma,\rho_{(\sigma\varsigma)^{-1}}(u))\,g(\rho_{(\sigma\varsigma)^{-1}}(u))\\
&=&\omega(\sigma\varsigma,\rho_{(\sigma\varsigma)^{-1}}(u))\,g(\rho_{(\sigma\varsigma)^{-1}}(u))\kern+0.5em=\kern+0.5em\bm{\rho}_{\sigma\varsigma}(g)(u).
\end{array}\]Thus, $\bm{\rho}$ is a well-defined group action on $\field^M$. Since $\omega(\sigma,\rho_{\sigma^{-1}}(u))\neq0$ for all $\sigma\in G$ and $u\in M$, we have\[\supp(\bm{\rho}_\sigma(g))=\{u\in M\mid g(\rho_{\sigma^{-1}}(u))\neq0\}=\{u\in M\mid \rho_{\sigma^{-1}}(u)\in\supp(g)\}=\rho_\sigma(\supp(g)).\]Furthermore, since $\rho$ is $\G$-graded, we have $\rho_{\sigma^{-1}}(\ubar{M}_x)=\ubar{M}_x$ for all $\sigma\in G$, and thus, for any $g\in\ubar{\field}_x^M$, $g(\rho_{\sigma^{-1}}(M\setminus\ubar{M}_x))=\{0\}$. Hence, $\bm{\rho}_\sigma(g)(M\setminus\ubar{M}_x)=\{0\}$, meaning that $\bm{\rho}_\sigma(g)\in\ubar{\field}_x^M$.

\eqref{027}~Assume now that $\rho$ is quasi-multiplicative, and let $u,v\in M$ with $\tau(u,v)\neq0$. Let $\alpha,\varsigma\in G$ be the elements given by the quasi-multiplicativity of $\rho_\sigma$ on $uv$. On the one hand, we have
\[\bm{\rho}_\sigma(au\cdot bv)=\bm{\rho}_\sigma(ab\,\tau(u,v)\,uv)=ab\,\tau(u,v)\,\omega(\sigma,uv)\,\rho_\sigma(uv)=ab\,\tau(u,v)\,\omega(\sigma,uv)\,\rho_\alpha(u)\rho_\varsigma(v).\]On the other hand, multiplying the respective images under $\alpha$ and $\varsigma$, we obtain\[\bm{\rho}_\alpha(au)\cdot\bm{\rho}_\varsigma(bv)=(a\,\omega(\alpha,u)\,\rho_\alpha(u))\cdot(b\,\omega(\varsigma,v)\,\rho_\varsigma(v))=ab\,\omega(\alpha,u)\,\omega(\varsigma,v)\,\tau(\rho_\alpha(u),\rho_\varsigma(v))\,\rho_\alpha(u)\rho_\varsigma(v).\]
Hence, $\bm{\rho}_\sigma(au\cdot bv)=\bm{\rho}_\alpha(au)\cdot\bm{\rho}_\varsigma(bv)$ if and only if $(\tau,\omega)$ is a multiplicative pair.

\eqref{181}~By~\eqref{016}, for each $x\in M$ the set of factorizations with non-zero twisting, $F(x):=\{(v,w)\in M\times M\mid vw=x\text{ and }\tau(v,w)\neq0\}$, is finite. Given $\sigma\in G$, set $z:=\rho_{\sigma^{-1}}(u)$. For any $(v,w)\in F(z)$, we have $\rho_\alpha(v)\rho_\varsigma(w)=\rho_\sigma(z)=u$, and since $\omega$ takes values in $\field^\times$, \eqref{028} implies that $\tau(\rho_\alpha(v),\rho_\varsigma(w))\neq0$. Because $\rho_\alpha$ and $\rho_\varsigma$ are bijections on $M$, the assignment $(v,w)\mapsto(\rho_\alpha(v),\rho_\varsigma(w))$ defines an injective map from $F(z)$ to $F(u)$, whence $|F(z)|\leq|F(u)|$. Applying the symmetric argument to $\sigma^{-1}$ and $u$ yields $|F(u)|\leq|F(z)|$. Since these sets are finite, this assignment is a bijection from $F(z)$ onto $F(u)$.

\eqref{031}~Assume now that $\rho$ is multiplicative and $(\tau,\omega)$ is a multiplicative pair. By~\eqref{181}, for any $\sigma\in G$ and $u\in M$, the map $(a,b)\mapsto(\rho_\sigma(a),\rho_\sigma(b))$ is a bijection between the sets of factorizations with non-zero twisting of $\rho_{\sigma^{-1}}(u)$ and $u$. Thus,\[\begin{array}{rcl}(\bm{\rho}_\sigma(g)\cdot\bm{\rho}_\sigma(h))(u)
&=&{\displaystyle\sum_{vw=u}\tau(v,w)\bm{\rho}_\sigma(g)(v)\bm{\rho}_\sigma(h)(w)}\\
&=&{\displaystyle\sum_{vw=u}\tau(v,w)\,\omega(\sigma,\rho_{\sigma^{-1}}(v))\,\omega(\sigma,\rho_{\sigma^{-1}}(w))\,g(\rho_{\sigma^{-1}}(v))\,h(\rho_{\sigma^{-1}}(w))}\\
&\stackrel{\eqref{028}}{=}&{\displaystyle\sum_{vw=u}\tau(\rho_{\sigma^{-1}}(v),\rho_{\sigma^{-1}}(w))\,\omega(\sigma,\rho_{\sigma^{-1}}(u))\,g(\rho_{\sigma^{-1}}(v))\,h(\rho_{\sigma^{-1}}(w))}\\
&=&{\displaystyle\omega(\sigma,\rho_{\sigma^{-1}}(u))\!\!\!\!\!\sum_{ab=\rho_{\sigma^{-1}}(u)}\!\!\!\!\tau(a,b)\,g(a)\,h(b)}\\
&=&{\displaystyle\omega(\sigma,\rho_{\sigma^{-1}}(u))\,(g\cdot h)(\rho_{\sigma^{-1}}(u))}\quad=\quad\bm{\rho}_\sigma(g\cdot h)(u).
\end{array}\]Hence, $\bm{\rho}_\sigma$ is an algebra automorphism.

\eqref{032}~By \eqref{030}, we have $|\supp(\bm{\rho}_\sigma(g))|=|\rho_\sigma(\supp(g))|=|\supp(g)|$ and $\bm{\rho}_\sigma(\ubar{\field}_x^M)=\ubar{\field}_x^M$ for all $\sigma\in G$, which implies that $\bm{\rho}_\sigma(\pols{\field}{\ubar{M}_x})=\pols{\field}{\ubar{M}_x}$. Thus, the map $\sigma\mapsto\bm{\rho}_\sigma|_{\pols{\field}{M}}$ defines a $\G$-graded group action on $\pols{\field}{M}$.

\eqref{033}~If $\rho$ is a bigraded action, then by \eqref{030} we have $\supp(\bm{\rho}_\sigma(g))=\rho_\sigma(\supp(g))\subseteq\rho_\sigma(M_x^d)=M_x^d$ for all $\sigma\in G$ and $g\in\funcs{\field}{M_x^d}$, which implies that $\bm{\rho}_\sigma(\funcs{\field}{M_x^d})\subseteq\funcs{\field}{M_x^d}$. Thus, the map $\sigma\mapsto\bm{\rho}_\sigma|_{\funcs{\field}{M}}$ defines a bigraded group action on the space $\funcs{\field}{M}$.

\eqref{034}~Let $g,h\in\funcs{\field}{M}^G$. For any $\sigma\in G$ and $u\in M$, we have\[\bm{\rho}_\sigma(g\cdot h)(u)=\omega(\sigma,\rho_{\sigma^{-1}}(u))\sum_{vw=\rho_{\sigma^{-1}}(u)}\tau(v,w)\,g(v)\,h(w).\]Let $\alpha,\varsigma\in G$ be the elements associated with $(\sigma,\rho_{\sigma^{-1}}(u))$ by the quasi-multiplicativity of $\rho$. Furthermore, the $G$-invariance of $g$ and $h$ implies that $\bm{\rho}_\alpha(g)=g$ and $\bm{\rho}_\varsigma(h)=h$. Evaluating these at $\rho_\alpha(v)$ and $\rho_\varsigma(w)$, respectively, we obtain\[g(\rho_\alpha(v))=\bm{\rho}_\alpha(g)(\rho_\alpha(v))=\omega(\alpha,v)g(v)\quad\text{and}\quad h(\rho_\varsigma(w))=\bm{\rho}_\varsigma(h)(\rho_\varsigma(w))=\omega(\varsigma,w)h(w).\]Substituting $g(v)=\omega(\alpha,v)^{-1}g(\rho_\alpha(v))$ and $h(w)=\omega(\varsigma,w)^{-1}h(\rho_\varsigma(w))$ into the sum above, and noting that $\omega(\sigma,\rho_{\sigma^{-1}}(u))=\omega(\sigma,vw)$, it follows that\[\begin{array}{rcl}
\bm{\rho}_\sigma(g\cdot h)(u)&=&{\displaystyle\sum_{\rho_\sigma(vw)=u}\tau(v,w)\,\omega(\sigma,vw)\,\omega(\alpha,v)^{-1}\,\omega(\varsigma,w)^{-1}\,g(\rho_\alpha(v))\,h(\rho_\varsigma(w))}\\
&\stackrel{\eqref{028}}{=}&{\displaystyle\sum_{\rho_\alpha(v)\rho_\varsigma(w)=u}\tau(\rho_\alpha(v),\rho_\varsigma(w))\,g(\rho_\alpha(v))\,h(\rho_\varsigma(w))}\kern+0.6em=\kern+0.6em{\displaystyle\sum_{ab=u}\tau(a,b)\,g(a)\,h(b)}\kern+0.6em=\kern+0.6em(g\cdot h)(u),
\end{array}\]because the assignment $(a,b)=(\rho_\alpha(v),\rho_\varsigma(w))$ defines the bijection between $F(\rho_{\sigma^{-1}}(u))$ and $F(u)$ established in~\eqref{181}. This proves that $g\cdot h$ is $G$-invariant. Moreover, since $\rho$ is bigraded and $M^0=\{1\}$, we have $\rho_\sigma(1)=1$ and thus $\bm{\rho}_\sigma(1)=\omega(\sigma,1)\,1=1$ for all $\sigma\in G$, so $1\in\pols{\field}{M}^G\subseteq\funcs{\field}{M}^G$. Since $\rho$ is a bigraded action, the invariant subspaces preserve the bigrading, making $\funcs{\field}{M}^G$ a bigraded algebra and $\pols{\field}{M}^G$ a bigraded subalgebra.
\end{proof}

Given the results in Proposition~\ref{029}, for the remainder of the section, we assume the following condition:
\stepcounter{condition}
\begin{enumerate}[start=\thecondition,label={(\condition)},ref={\condition}]
\item\label{018}\conditions{3}
\end{enumerate}
Consequently, the space $\funcs{\field}{M}$ is equipped with a bigraded action $\bm{\rho}:G\to\GL(\funcs{\field}{M})$. This action defines the subalgebra of invariants $\funcs{\field}{M}^G$, referred to as the \emph{twisted algebra of invariant functions} on $M$. Furthermore, the bigrading of $\funcs{\field}{M}^G$ induces natural decompositions of the invariant subspace:\[\funcs{\field}{M}^G\,\,\,=\,\,\,\bigoplus_{d\in\N}\funcs{\field}{M^d}^G\,\,\,=\,\,\,\bigoplus_{x\in\G}\funcs{\field}{\ubar{M}_x}^G\,\,\,=\bigoplus_{(x,d)\in\G\times\N}\funcs{\field}{M_x^d}^G.\]

\begin{rem}[$G$-invariance]\label{044}
A function $g\in\funcs{\field}{M}$ is $G$-invariant if and only if $\bm{\rho}_{\sigma^{-1}}(g)=g$ for all $\sigma\in G$, which is equivalent to having\begin{equation}\label{043}\bm{\rho}_{\sigma^{-1}}(g)(u)=\omega(\sigma^{-1},\rho_\sigma(u))\,g(\rho_\sigma(u))=g(u)\quad\stackrel{\eqref{003}}{\Leftrightarrow}\quad g(\rho_\sigma(u))=\omega(\sigma,u)\,g(u).\end{equation}
\end{rem}

\begin{rem}[Restriction]\label{049}
For any subset $X\subseteq M$ and any function $g\in\funcs{\field}{M}$, define the \emph{restricted function} $g_{_X}\in\funcs{\field}{M}$ by setting $g_{_X}(u)=g(u)$ if $u\in X$, and $g_{_X}(u)=0$ otherwise. If $g$ is $G$-invariant, and the subset $X$ is \emph{$G$-stable}, that is, $\rho_\sigma(X)\subseteq X$ for all $\sigma\in G$, then $g_{_X}$ is also $G$-invariant. To see this, let $\sigma\in G$ and $u\in M$. If $u\in X$, the $G$-stability of $X$ guarantees that $\rho_\sigma(u)\in X$. Thus, the invariance of $g$ implies that $g_{_X}(\rho_\sigma(u))=g(\rho_\sigma(u))=\omega(\sigma,u)g(u)=\omega(\sigma,u)g_{_X}(u)$. On the other hand, if $u\notin X$, we must also have $\rho_\sigma(u)\notin X$ because $X$ is $G$-stable. Therefore, both $g_{_X}(\rho_\sigma(u))$ and $g_{_X}(u)$ evaluate to zero, and the equality $0=\omega(\sigma,u)\cdot0$ holds trivially. Remark~\ref{044} implies then that $g_{_X}$ is $G$-invariant. In particular, since any orbit $Gu$ is a $G$-stable subset, the function $g_{_{Gu}}$ is a $G$-invariant function supported on $Gu$.
\end{rem}

\subsection{Linear bases}\label{190}

In this subsection, we establish a linear basis for the twisted algebra of $G$-invariant functions $\funcs{\field}{M}^G$ and describe its multiplication.

\begin{pro}\label{045}
Let $u\in M$ be a simple monomial. There exists a nonzero $G$-invariant function $g\in \funcs{\field}{M}^G$ with $\supp(g)\subseteq Gu$ if and only if $\omega(\varsigma,u)=1$ for all $\varsigma\in G_{(u)}$.
\end{pro}
\begin{proof}
Suppose that $g$ is a nonzero $G$-invariant function with $\supp(g)\subseteq Gu$. As $g$ is nonzero, there exists some $v\in\supp(g)$ such that $g(v)\neq0$. Let $\sigma\in G$ be an element satisfying $\rho_\sigma(u)=v$. By the $G$-invariance of $g$, Remark~\ref{044} implies that $g(v)=g(\rho_\sigma(u))=\omega(\sigma,u)g(u)$. The fact that $g(v)\neq0$ forces $g(u)\neq0$. Now, let $\varsigma\in G_{(u)}$. Applying Remark~\ref{044} again and noting that $\rho_\varsigma(u)=u$, we obtain $g(u)=g(\rho_\varsigma(u))=\omega(\varsigma,u)g(u)$. Because $g(u)\neq0$, we conclude that $\omega(\varsigma,u)=1$ for all $\varsigma\in G_{(u)}$.

Conversely, assume that $\omega(\varsigma,u)=1$ for all $\varsigma\in G_{(u)}$. We construct a nonzero formal power series $m_u$ supported on $Gu$, which, since the action $\rho$ is bigraded, will be homogeneous and therefore of bounded degree. For any $v\in Gu$, choose $\sigma\in G$ such that $\rho_\sigma(u)=v$ and define $m_u(v)=\omega(\sigma,u)$. For $v \notin Gu$, we set $m_u(v)=0$. Note that if $\varsigma\in G$ is another element satisfying $\rho_\varsigma(u)=v$, then $\varsigma^{-1}\sigma\in G_{(u)}$, and condition \eqref{003} implies that $\omega(\sigma,u)=\omega(\varsigma(\varsigma^{-1}\sigma),u)=\omega(\varsigma,\rho_{\varsigma^{-1}\sigma}(u))\omega(\varsigma^{-1}\sigma,u)=\omega(\varsigma,\rho_{\varsigma^{-1}\sigma}(u))=\omega(\varsigma,u)$. Thus, $m_u(v)$ is independent of the choice of $\sigma$. By construction, $m_u(u)=\omega(1,u)=1\neq0$. Finally, for any $\varsigma\in G$ and $v \in Gu$, we have $m_u(\rho_\varsigma(v))=m_u(\rho_{\varsigma\sigma}(u))=\omega(\varsigma\sigma,u)=\omega(\varsigma,\rho_\sigma(u))\omega(\sigma,u)=\omega(\varsigma,v)m_u(v)$, which proves that $g:=m_u$ is a nonzero $G$-invariant function.
\end{proof}

Motivated by the previous proposition, we introduce the following terminology. A simple monomial $u$ is called \emph{$\omega$-admissible} if $\omega(\sigma,u)=1$ for all $\sigma\in G_{(u)}$. In this case, there is a nonzero $G$-invariant function $m_u$ constructed in the proof of Proposition~\ref{045}, defined by $m_u(v)=\omega(\sigma,u)$ for $v=\rho_\sigma(u)\in Gu$, and $0$ elsewhere, which we call the \emph{monomial $G$-invariant function} associated with $u$. Note that if $\omega$ is identically $1$, then every simple monomial is $\omega$-admissible, and therefore $m_u$ coincides with the classical indicator function $\delta_{Gu}$ which takes value $1$ on every element of $Gu$ and $0$ elsewhere.

From the proof of Proposition~\ref{045} we obtain the following result.
\begin{crl}\label{066}
Let $u\in M$ be an $\omega$-admissible simple monomial, and let $v\in Gu$. Then $\omega(\sigma,u)=\omega(\varsigma,u)$ for all $\sigma,\varsigma\in G$ satisfying $\rho_\sigma(u)=v=\rho_\varsigma(u)$.
\end{crl}

\begin{rem}
We say that $\omega$ is a \emph{coboundary} if there exists a map $\delta:M\to\field^\times$ such that $\omega(\sigma,u)=\delta(\rho_\sigma(u))\delta(u)^{-1}$ for all $\sigma\in G$ and $u\in M$ (cf. \cite[Subsection~8.1]{Aar97}). In this case, every simple monomial is $\omega$-admissible; indeed, $\omega(\sigma,u)=\delta(\rho_\sigma(u))\delta(u)^{-1}=\delta(u)\delta(u)^{-1}=1$ for all $\sigma\in G_{(u)}$. Consequently, the simple monomials take the explicit form $m_u(v)=\omega(\sigma,u)=\delta(v)\delta(u)^{-1}$ for all $v\in Gu$ with $v=\rho_\sigma(u)$.
\end{rem}

\begin{lem}\label{046}
Let $u\in M$ be $\omega$-admissible, and let $v\in Gu$. Then, $v$ is also $\omega$-admissible. Furthermore, if $\sigma\in G$ is such that $\rho_\sigma(u)=v$, then $m_u=\omega(\sigma,u)m_v$. In particular, $m_u$ and $m_v$ are linearly dependent.
\end{lem}
\begin{proof}
Since $u$ is $\omega$-admissible, the monomial $G$-invariant function $m_u$ is well-defined, is non-zero, and is supported exactly on $Gu$. Because $v \in Gu$, we have $Gv = Gu$, meaning $m_u$ is a non-zero $G$-invariant function with support contained in $Gv$. By Proposition~\ref{045}, this immediately implies that $v$ is $\omega$-admissible, guaranteeing the existence of its own monomial $G$-invariant function $m_v$.

Now, let $\sigma\in G$ such that $\rho_\sigma(u)=v$, and let $w\in Gu$. Since $w\in Gv$, there exists $\varsigma\in G$ such that $w=\rho_\varsigma(v)=\rho_\varsigma(\rho_\sigma(u))=\rho_{\varsigma\sigma}(u)$. Applying condition \eqref{003}, we obtain $m_u(w)=\omega(\varsigma\sigma,u)=\omega(\varsigma, \rho_\sigma(u))\omega(\sigma,u)=\omega(\varsigma,v)\omega(\sigma,u)$. On the other hand, the definition of $m_v$ implies that $m_v(w)=\omega(\varsigma,v)$. Substituting this into the previous identity yields $m_u(w)=m_v(w)\omega(\sigma,u)$. Since this equality holds for all $w\in Gu$ and both functions vanish outside this orbit, we conclude that $m_u=\omega(\sigma, u)m_v$.
\end{proof}

The previous lemma allows us to unambiguously extend the property of being $\omega$-admissible to entire orbits. An orbit in $M/R_G$ is called \emph{$\omega$-admissible} if any of its elements is $\omega$-admissible.

\begin{pro}\label{048}
Let $T\subseteq M$ be a system of representatives for the $\omega$-admissible orbits in $M/R_G$. The set $\{m_t\mid t\in T\}$ is linearly independent. Furthermore, every $G$-invariant function $g\in\funcs{\field}{M}^G$ can be uniquely written as a pointwise well-defined sum\[g=\sum_{t\in T}g(t)m_t.\]
\end{pro}
\begin{proof}
Let $t_1,\ldots,t_k\in T$ be distinct elements, and suppose that $a_1m_{t_1}+\cdots+a_km_{t_k}=0$ for some scalars. Since the functions $m_{t_1},\ldots,m_{t_k}$ have disjoint supports, evaluating this linear combination at $t_i$ yields $0=(a_1m_{t_1}+\cdots+a_km_{t_k})(t_i)=a_im_{t_i}(t_i)=a_i$ for each $i$. This proves linear independence.

Now, let $g\in\funcs{\field}{M}^G$ and let $u\in M$. Remark~\ref{049} implies that the restricted function $g_{_{Gu}}$ is a $G$-invariant function supported on the orbit $Gu$. If $Gu$ is not $\omega$-admissible, Proposition~\ref{045} implies that the only such invariant function is the zero function. Consequently, $g_{_{Gu}}=0$, and in particular $g(Gu)=\{0\}$. 

If $Gu$ is $\omega$-admissible, let $t\in T\cap Gu$ be its unique representative. Consider the function $h=g_{_{Gu}}-g(t)m_t$. Since both $g_{_{Gu}}$ and $m_t$ are $G$-invariant functions supported on $Gu$, $h$ is also a $G$-invariant function supported on $Gu$. Evaluating $h$ at the representative $t$ yields $h(t)=g_{_{Gu}}(t)-g(t)m_t(t)=g(t)-g(t)(1)=0$. Because $h$ is $G$-invariant, Remark~\ref{044} implies that for any $w=\rho_\sigma(t) \in Gu$, we have $h(w)=h(\rho_\sigma(t))=\omega(\sigma,t)h(t)=0$. This means that $h$ vanishes on the entire orbit $Gu$, concluding that $g_{_{Gu}}=g(t)m_t$. Evaluating at $u$, we obtain $g(u)=g_{_{Gu}}(u)=g(t)m_t(u)$. Because the orbits form a partition of $M$, for any $u\in M$ there is at most one $t\in T$ such that $m_t(u)\neq0$. Thus, the sum $\sum g_{_{Gt}}=\sum g(t)m_t$ is pointwise well-defined and equals $g$. The uniqueness of the coefficients $g(t)$ follows directly from the disjointness of the supports.
\end{proof}

\begin{thm}\label{050}
Let $T\subseteq M$ be a system of representatives for the $\omega$-admissible $G$-orbits in $M/R_G$. If each homogeneous component $M^d$ intersects only finitely many $\omega$-admissible $G$-orbits, then the set $\{m_t\mid t\in T\}$ defines a linear basis for the vector space $\funcs{\field}{M}^G$.
\end{thm}
\begin{proof}
By Proposition~\ref{048}, the set $\{m_t\mid t\in T\}$ is linearly independent, and any $G$-invariant function $g$ can be written as the pointwise well-defined sum $g=\sum g(t)m_t$. Since $g$ is of bounded degree, its support is contained in a finite union of homogeneous components $M^d$. By our assumption, each $M^d$ in this finite union intersects only finitely many $\omega$-admissible $G$-orbits. Consequently, the set of representatives $t \in T$ whose orbit $Gt$ intersects the support of $g$ is finite. Thus, the coefficients $g(t)$ are non-zero for only finitely many $t\in T$. The sum reduces to a finite linear combination, proving that $\{m_t\mid t\in T\}$ spans $\funcs{\field}{M}^G$.
\end{proof}

Given the result of Theorem~\ref{050}, to ensure that the space of invariant functions possesses a well-behaved linear basis, we assume the following condition for the remainder of the section.
\stepcounter{condition}
\begin{enumerate}[start=\thecondition,label={(\condition)},ref={\condition}]
\item\label{051}\conditions{4}
\end{enumerate}
Under Condition~\eqref{051}, the space $\funcs{\field}{M}^G$ possesses a linear basis $\{m_t\mid t\in T\}$, which we shall refer to as the \emph{monomial basis} of $\funcs{\field}{M}^G$ with respect to $T$.

\begin{rem}
Since $\rho$ is bigraded, it follows that $\supp(m_u)\subseteq M_x^d$ for every simple monomial $u\in M_x^d$. Consequently, the degrees of the invariant function match those of the monomial, meaning $\deg(m_u)=\deg(u)$ and $|m_u|=|u|$. In particular, this explicitly yields $\phi(|m_u|,|m_v|)=\phi(|u|,|v|)$ for any $\omega$-admissible $u,v\in M$. Furthermore, this implies that the sets $\{m_t\mid t\in T\cap M^d\}$, $\{m_t\mid t\in T\cap\ubar{M}_x\}$, and $\{m_t\mid t\in T\cap M_x^d\}$ constitute linear bases of the invariant subspaces $\funcs{\field}{M^d}^G$, $\funcs{\field}{\ubar{M}_x}^G$, and $\funcs{\field}{M_x^d}^G$, respectively.
\end{rem}

We finish this subsection by determining the decomposition in the monomial basis of the product of two monomial $G$-invariant functions.
\begin{pro}\label{052}
For all $\omega$-admissible simple monomials $u,v\in M$, the twisted convolution product $m_u\cdot m_v$ decomposes as\[m_u\cdot m_v=\sum_{t\in T}c_{u,v}^t\,m_t,\]where the structure constants $c_{u,v}^t\in\field$ are explicitly given by\begin{equation}\label{053}c_{u,v}^t\,=\!\!\!\sum_{\substack{(x,y)\in Gu\times Gv\\xy=t}}\!\!\!\!\tau(x,y)\,m_u(x)\,m_v(y).\end{equation}
\end{pro}
\begin{proof}
Since $\funcs{\field}{M}^G$ is an algebra, the product $m_u\cdot m_v$ is a $G$-invariant function, and, by Proposition~\ref{048}, it admits a unique linear expansion in terms of the monomial basis $m_u\cdot m_v=\sum c_{u,v}^t\,m_t$, where $c_{u,v}^t=(m_u\cdot m_v)(t)$. Applying the definition of the twisted convolution product immediately yields \eqref{053}, where the sum is naturally restricted to $x\in\supp(m_u)=Gu$ and $y\in\supp(m_v)=Gv$.
\end{proof}

%In what follows, the sequence formed by the monoid $M$, the twisting $\tau:M\times M\to\field$, the group $\G$, the bicharacter $\phi:\G\times\G\to\field^\times$, the degree maps $|\,|:M\to\G$ and $\deg:M\to\N$, the group $G$, the action $\rho:G\to\sym(M)$, the weight function $\omega:G\times M\to\field^\times$, and the system of representatives $T$, satisfying conditions \eqref{015} to \eqref{051}, will be called an \emph{invariant algebraic datum}:\[(M,\tau,\G,\phi,|\,|,\deg,G,\rho,\omega,T).\]

\section{Additional structures}\label{168}

In this section, we explore several algebraic constructions that naturally arise from our previous development.

As usual, $M$ denotes a monoid that is locally finite with respect to a twisting $\tau:M\times M\to\field$, and $G$ denotes a group equipped with a bigraded group action $\rho:G\to\sym(M)$ and a weight function $\omega:G\times M\to\field^\times$, satisfying conditions \eqref{015}~to~\eqref{051}. Under these assumptions, the twisted algebra of invariants $\funcs{\field}{M}^G$ is well-defined and admits the monomial basis $\{m_t\mid t\in T\}$, where $T$ is a system of representatives of the $\omega$-admissible $G$-orbits.

We begin by analysing the behaviour of $G$-invariant functions under restricted actions and examining the interaction with direct products of monoids. Subsequently, we determine the conditions under which $\funcs{\field}{M}^G$ admits a Hopf algebra structure and conclude by studying the induced structures on quotient monoids.

\subsection{Restricted actions}\label{171}

An \emph{ideal with respect to $\tau$} of $M$ is a subset $I\subseteq M$ such that for any $u,v\in M$ with $\tau(u,v)\neq0$, if $u\in I$ or $v\in I$, then $uv\in I$~(cf.~\cite[Subsection~1.1]{How95}). Equivalently, $I$ is an ideal with respect to $\tau$ if its complement $M\setminus I$ satisfies that for any $u,v\in M$ with $\tau(u,v)\neq0$, if $uv\in M\setminus I$, then $u,v\in M\setminus I$.

Let $S$ be a set of generators of $G$ such that $S^{-1}=S$. For each $s\in S$, we consider a subset $M_s\subseteq M$ that is stable and multiplicative with respect to $\tau$ under the action of $s$, and such that its complement $M\setminus M_s$ is an ideal of $M$ with respect to $\tau$; that is, $\rho_s(M_s)=M_s$, and if $uv\in M_s$ with $\tau(u,v)\neq0$, then $u,v\in M_s$, $\rho_s(uv)=\rho_s(u)\rho_s(v)$, and~\eqref{028} holds for $\sigma=s$ with $\alpha=\varsigma=s$. We define $\bar{\rho}_s:M\to M$ by $\bar{\rho}_s(u)=\rho_s(u)$ if $u\in M_s$, and $\bar{\rho}_s(u)=u$ otherwise.

\begin{lem}\label{056}
For each $s \in S$, the map $\bar{\rho}_s: M \to M$ is a bigraded bijection.
\end{lem}
\begin{proof}
By hypothesis, the original action $\rho_s$ is bigraded. Since $\bar{\rho}_s(u)$ is defined piece-wise as either $\rho_s(u)$ or $u$, both of which have the exact same bidegree as $u$, it follows immediately that $\bar{\rho}_s$ is a bigraded map.

Since $\rho_s$ is a bijection on $M$, it suffices to show that $\bar{\rho}_s$ maps $M_s$ onto $M_s$ and $M \setminus M_s$ onto $M \setminus M_s$ bijectively. By the stability condition, $\rho_s(M_s)=M_s$, and since $\rho_s$ is bijective, it restricts to a bijection on $M_s$. On the other hand, $\bar{\rho}_s$ acts as the identity on $M\setminus M_s$, which is clearly a bijection. This proves that $\bar{\rho}_s$ is bijective, completing the proof.
\end{proof}

Lemma~\ref{056} guarantees that the assignment $s\mapsto\bar{\rho}_s$ extends uniquely to a homomorphism $\bar{\rho}^*:F(S)\to\sym(M)$, where $F(S)$ denotes the free group generated by $S$. If $R$ is a relation on $F(S)$ such that $\bar{\rho}^*(r)=\bar{\rho}^*(t)$ for all relations $(r,t)\in R$, and $\bar{G}=F(S)/\bar{R}$, then there exists a unique group homomorphism $\bar{\rho}:\bar{G}\to\sym(M)$ satisfying $\bar{\rho}(s)=\bar{\rho}_s$ for all $s\in S$. We call $\bar{\rho}$ the \emph{restricted action} of $G$ on $M$ with respect to $\{M_s\mid s\in S\}$.

\begin{pro}\label{054}
\begin{enumerate}
\item\label{060}The restricted action only produces elements that are reached by the original action; that is, for any $\sigma\in\bar{G}$ and $u\in M$ there exists $\varsigma\in G$ such that $\bar{\rho}_\sigma(u)=\rho_\varsigma(u)$.
\item\label{061}For any $u\in M$, $\bar{G}u\subseteq Gu$. Consequently, the partition of $M$ into $\bar{G}$-orbits is a refinement of the partition into $G$-orbits.
\item\label{062}The restricted action $\bar{\rho}$ is a bigraded quasi-multiplicative action.
\end{enumerate}
\end{pro}
\begin{proof}
\eqref{060}~We proceed by induction on the length of a word representing $\sigma$. If $\sigma=1$, the claim holds trivially with $\varsigma=1$. If $\sigma=s$ for some generator $s\in S$, then by definition $\bar{\rho}_s(u)$ is either $\rho_s(u)$ or $u=\rho_1(u)$; in either case, the claim holds with $\varsigma=s$ or $\varsigma=1$. Suppose the claim holds for some $\sigma'\in\bar{G}$, where $\sigma=s\sigma'$ and $\bar{\rho}_{\sigma'}(u)=\rho_{\varsigma'}(u)$ for some $\varsigma'\in G$. Then $\bar{\rho}_\sigma(u)=\bar{\rho}_s(\bar{\rho}_{\sigma'}(u))=\bar{\rho}_s(\rho_{\varsigma'}(u))$. By the base case, there exists some $s'\in\{s,1\}$ such that $\bar{\rho}_s(\rho_{\varsigma'}(u))=\rho_{s'}(\rho_{\varsigma'}(u))=\rho_{s'\varsigma'}(u)$. Taking $\varsigma=s'\varsigma'$ completes the induction.

\eqref{061}~Given $v\in\bar{G}u$, there exists $\sigma\in\bar{G}$ such that $v=\bar{\rho}_\sigma(u)$. By \eqref{060}, there exists $\varsigma\in G$ such that $v=\rho_\varsigma(u)$, which implies $v\in Gu$, and so $\bar{G}u\subseteq Gu$. Since both group actions partition $M$ into disjoint orbits, this inclusion precisely means that the partition into $\bar{G}$-orbits is a refinement of the partition into $G$-orbits.

\eqref{062}~Let $w\in M$ and $s\in S$, and consider $u,v\in M$ such that $uv=w$ and $\tau(u,v)\neq0$. If $w\in M_s$, since $M\setminus M_s$ is an ideal of $M$ with respect to $\tau$, we have $u,v\in M_s$, and so $\bar{\rho}_s(uv)=\rho_s(uv)=\rho_s(u)\rho_s(v)=\bar{\rho}_s(u)\bar{\rho}_s(v)$. If $w\notin M_s$, then $\bar{\rho}_s(uv)=uv=\bar{\rho}_1(u)\bar{\rho}_1(v)$. Thus, for each generator $s\in S$ and each $w\in M$, the pair $(\alpha_s,\varsigma_s)\in\{(s,s),(1,1)\}$ depends only on $w$, and by~\eqref{028} on $M_s$ it satisfies $\tau(\bar{\rho}_{\alpha_s}(u),\bar{\rho}_{\varsigma_s}(v))\neq0$ whenever $\tau(u,v)\neq0$. Consequently, these uniform choices compose along any word in the generators $S$, showing that $\bar{\rho}$ is quasi-multiplicative with respect to $\tau$. Moreover, as $\bar{\rho}_s$ preserves the bigrading for each generator, $\bar{\rho}$ is a bigraded quasi-multiplicative action.
\end{proof}

For each $s\in S$, let $\bar{\omega}_s:M\to\field^\times$ be defined by $\bar{\omega}_s(u)=\omega(s,u)$ if $u\in M_s$, and $\bar{\omega}_s(u)=1$ otherwise. Let $\bar{\omega}^*:F(S)\times M\to\field^\times$ be the unique map extending the assignments $s\mapsto\bar{\omega}_s$ and determined by the $1$-cocycle condition $\bar{\omega}^*(gh,u)=\bar{\omega}^*(g,\bar{\rho}_h^*(u))\,\bar{\omega}^*(h,u)$ for all $g,h\in F(S)$. If $\bar{\omega}^*(r,u)=\bar{\omega}^*(t,u)$ for all $(r,t)\in R$ and all $u\in M$, then there exists a unique well-defined weight function $\bar{\omega}:\bar{G}\times M\to\field^\times$ such that $\bar{\omega}(s,u)=\bar{\omega}_s(u)$ for all $s\in S$. We call $\bar{\omega}$ the \emph{restricted weight function} for $\bar{\rho}$.

\begin{pro}\label{055}
\begin{enumerate}
\item\label{063} For any $\sigma\in\bar{G}$ and $u\in M$, we have $\bar{\omega}(\sigma,u)=\omega(\varsigma,u)$, where $\varsigma$ is the element in $G$ found in Proposition~\ref{054}\eqref{060}.
\item\label{064} If $u\in M$ is $\omega$-admissible and $v\in Gu$, then $v$ is $\bar{\omega}$-admissible.
\item\label{065} The pair $(\tau,\bar{\omega})$ is multiplicative for the restricted action $\bar{\rho}$.
\end{enumerate}
\end{pro}
\begin{proof}
\eqref{063}~We follow the exact same induction as in the proof of Proposition~\ref{054}\eqref{060}. If $\sigma=1$, then $\varsigma=1$ and $\bar{\omega}(1,u)=1=\omega(1,u)$. If $\sigma=s$ for some generator $s\in S$, the element $\varsigma\in\{s,1\}$ chosen in Proposition~\ref{054}\eqref{060} perfectly matches the piecewise definition of the weight, yielding $\bar{\omega}(s,u)=\omega(\varsigma,u)$. Suppose the claim holds for $\sigma'\in\bar{G}$, that is, $\bar{\omega}(\sigma',u)=\omega(\varsigma',u)$. If $\sigma=s\sigma'$, as in the proof of Proposition~\ref{054}\eqref{060}, we take $\varsigma=s'\varsigma'$ with $s'\in\{s,1\}$, where $\bar{\rho}_{\sigma'}(u)=\rho_{\varsigma'}(u)$ and $\bar{\rho}_s(\rho_{\varsigma'}(u))=\rho_{s'}(\rho_{\varsigma'}(u))$. By the 1-cocycle condition, $\bar{\omega}(\sigma,u)=\bar{\omega}(s,\bar{\rho}_{\sigma'}(u))\,\bar{\omega}(\sigma',u)=\bar{\omega}(s,\rho_{\varsigma'}(u))\,\omega(\varsigma',u)=\omega(s',\rho_{\varsigma'}(u))\,\omega(\varsigma',u)=\omega(\varsigma,u)$.

\eqref{064}~Since $u$ is $\omega$-admissible and $v\in Gu$, the element $v$ is also $\omega$-admissible. If $\sigma\in\bar{G}_{(v)}$, by \eqref{063} and Proposition~\ref{054}\eqref{060}, there exists $\varsigma\in G$ such that $\rho_\varsigma(v)=\bar{\rho}_\sigma(v)=v$ and $\bar{\omega}(\sigma,v)=\omega(\varsigma,v)$. Since $v$ is $\omega$-admissible and $\varsigma\in G_{(v)}$, we must have $\bar{\omega}(\sigma,v)=\omega(\varsigma,v)=1$, which proves that $v$ is $\bar{\omega}$-admissible.

\eqref{065}~Let $w\in M$ and $s\in S$, and consider $u,v\in M$ such that $uv=w$ and $\tau(u,v)\neq0$. If $w\in M_s$, since $M\setminus M_s$ is an ideal with respect to $\tau$, we have $u,v\in M_s$. Because $(\tau,\omega)$ is a multiplicative pair and $\rho_s(uv)=\rho_s(u)\rho_s(v)$, we have $\tau(u,v)\,\bar{\omega}(s,uv)=\tau(u,v)\,\omega(s,uv)=\tau(\rho_s(u),\rho_s(v))\,\omega(s,u)\omega(s,v)=\tau(\bar{\rho}_s(u),\bar{\rho}_s(v))\,\bar{\omega}(s,u)\,\bar{\omega}(s,v)$. If $w\notin M_s$, then $\bar{\rho}_s(uv)=uv$ and $\bar{\omega}(s,uv)=1$, hence $\tau(u,v)\,\bar{\omega}(s,uv)=\tau(u,v)=\tau(\bar{\rho}_1(u),\bar{\rho}_1(v))\,\bar{\omega}(1,u)\,\bar{\omega}(1,v)$. Since the uniform witnesses $(\alpha_s,\varsigma_s)\in\{(s,s),(1,1)\}$ preserve the non-vanishing of $\tau$ and compose along words in the generators $S$, induction on the word length guarantees that $(\tau,\bar{\omega})$ is a multiplicative pair for $\bar{\rho}$.
\end{proof}

By Proposition~\ref{054}\eqref{062} and Proposition~\ref{055}\eqref{065}, the restricted action $\bar{\rho}:\bar{G}\to\sym(M)$ together with the weight function $\bar{\omega}:\bar{G}\times M\to\field^\times$ satisfy condition \eqref{018}. This ensures the existence of the twisted algebra of $\bar{G}$-invariant functions, denoted by $\funcs{\field}{M}^{\bar{G}}$. By Proposition~\ref{054}\eqref{061}, the condition of $\bar{G}$-invariance is less restrictive than $G$-invariance. Consequently, we obtain the natural inclusion $\funcs{\field}{M}^G \subseteq \funcs{\field}{M}^{\bar{G}}$.

\begin{thm}\label{057}
Assume that every fibre $M^d$ intersects only finitely many $\bar{\omega}$-admissible $\bar{G}$-orbits. Then:
\begin{enumerate}
\item\label{058}The algebra $\funcs{\field}{M}^{\bar{G}}$ admits a monomial basis formed by the elements\[\bar{m}_t=\sum_{u\in\bar{G}t}\bar{\omega}(\sigma, t)u,\]where $t$ runs over a system of representatives $\bar{T}$ of the $\bar{\omega}$-admissible $\bar{G}$-orbits, and $\sigma\in\bar{G}$ is any element such that $\bar{\rho}_\sigma(t)=u$.
\item\label{059}For any $t\in T$, the corresponding monomial $G$-invariant function $m_t$ decomposes as\[m_t=\sum_{u\in T_t}\omega(\sigma,t)\,\bar{m}_u,\]where $T_t$ is a set of representatives for the $\bar{\omega}$-admissible $\bar{G}$-orbits contained in $Gt$, and $\sigma\in G$ is any element such that $\rho_\sigma(t)=u$. 
\end{enumerate}
\end{thm}
\begin{proof}
\eqref{058} The hypothesis that every fibre $M^d$ intersects only finitely many $\bar{\omega}$-admissible $\bar{G}$-orbits guarantees that condition \eqref{051} is satisfied. Consequently, by Theorem~\ref{050}, the set of $\bar{G}$-invariant functions $\bar{m}_t$ with $t\in\bar{T}$ constitute the canonical monomial basis for the algebra $\funcs{\field}{M}^{\bar{G}}$.

\eqref{059} The natural inclusion $\funcs{\field}{M}^G\subseteq\funcs{\field}{M}^{\bar{G}}$ implies that the monomial $G$-invariant function $m_t$ can be expressed as a unique linear combination in the monomial basis $\{\bar{m}_u\mid u\in\bar{T}\}$. By Proposition~\ref{055}\eqref{064}, every element in $Gt$ is $\bar{\omega}$-admissible, meaning that $Gt$ decomposes entirely into a disjoint union of $\bar{\omega}$-admissible $\bar{G}$-orbits $\bar{G}u$ with $u\in T_t$. For each $v\in Gt$, choose $\sigma_v\in G$ such that $\rho_{\sigma_v}(t)=v$. The function $m_t$ can be written as:
\[m_t=\sum_{v\in Gt}\omega(\sigma_v,t)v=\sum_{u\in T_t}\hat{m}_u,\quad\text{where}\quad\hat{m}_u=\sum_{v\in\bar{G}u}\omega(\sigma_v,t)v.\]
Note that $\hat{m}_u$ is the restricted function of $m_t$ with respect to the set $\bar{G}u$, defined in Remark~\ref{049}. Because $m_t$ is $\bar{G}$-invariant and $\bar{G}u$ is $\bar{G}$-stable, Remark~\ref{049} implies that each $\hat{m}_u$ is a $\bar{G}$-invariant function. For any $v\in\bar{G}u$, choose $\bar{\sigma}_v\in\bar{G}$ such that $\bar{\rho}_{\bar{\sigma}_v}(u)=v$. By Proposition~\ref{054}\eqref{060} and Proposition~\ref{055}\eqref{063}, there exists $\varsigma_v\in G$ such that $\rho_{\varsigma_v}(u)=v$ and $\bar{\omega}(\bar{\sigma}_v,u)=\omega(\varsigma_v,u)$. Observe that $\rho_{\varsigma_v\sigma_u}(t)=\rho_{\varsigma_v}(\rho_{\sigma_u}(t))=\rho_{\varsigma_v}(u)=v$. Since $t$ is $\omega$-admissible, Corollary~\ref{066} and the 1-cocycle condition imply that $\omega(\sigma_v,t)=\omega(\varsigma_v\sigma_u,t)=\omega(\varsigma_v,\rho_{\sigma_u}(t))\,\omega(\sigma_u,t)=\omega(\varsigma_v,u)\,\omega(\sigma_u,t)=\bar{\omega}(\bar{\sigma}_v,u)\,\omega(\sigma_u,t)$. This shows that $\hat{m}_u=\omega(\sigma_u,t)\bar{m}_u$, which yields the desired unique linear combination.
\end{proof}

\subsection{Products}\label{120}

The primary goal of this subsection is to establish that the twisted algebra of invariant functions on a direct product of monoids is isomorphic to the tensor product of their respective twisted algebras of invariants.\medskip

Let $M_1,\ldots,M_k$ be monoids, each equipped with a $\G$-grading $\{\ubar{M}_{i,x}\mid x\in\G\}$ with respect to the trivial twisting. Their direct product $M_1\times\cdots\times M_k$ naturally inherits a $\G$-grading with degree map defined by $|(u_1,\ldots,u_k)|=|u_1|\cdots|u_k|$. Accordingly, the homogeneous component of degree $x\in\G$ is given by the following disjoint union\begin{equation}\label{090}\bigsqcup_{x_1\cdots x_k=x}\ubar{M}_{1,x_1}\times\cdots\times\ubar{M}_{k,x_k}.\end{equation}

\begin{pro}\label{083}
The map $\tau^\phi:(M_1\times\cdots\times M_k)\times(M_1\times\cdots\times M_k)\to\field$ defined by\[\tau^\phi(u,v)=\prod_{i\,<\,j}\phi(|u_j|,|v_i|),\quad\text{where}\quad u=(u_1,\ldots,u_k)\quad\text{and}\quad v=(v_1,\ldots,v_k),\]defines a twisting on the direct product $M_1\times\cdots\times M_k$.
\end{pro}
\begin{proof}
Let $u=(u_1,\ldots,u_k)$, $v=(v_1,\ldots,v_k)$, and $w=(w_1,\ldots,w_k)$ be elements in $M_1\times\cdots\times M_k$. First, observe that\[\tau^\phi(uv,w)=\prod_{i\,<\,j}\phi(|u_jv_j|,|w_i|)\stackrel{\eqref{082}}{=}\prod_{i\,<\,j}\phi(|u_j|,|w_i|)\,\phi(|v_j|,|w_i|)=\tau^\phi(u,w)\,\tau^\phi(v,w).\]Similarly, we obtain that $\tau^\phi(u,vw)=\tau^\phi(u,v)\,\tau^\phi(u,w)$. Consequently, we have\[\tau^\phi(u,v)\,\tau^\phi(uv,w)=\tau^\phi(u,v)\,\tau^\phi(u,w)\,\tau^\phi(v,w)=\tau^\phi(u,vw)\,\tau^\phi(v,w).\]On the other hand, since $1=(1_{_{M_1}},\ldots,1_{_{M_k}})$, we have\[\tau^\phi(1,u)=\prod_{i\,<\,j}\phi(|1_{_{M_j}}|,|u_i|)=\prod_{i\,<\,j}\phi(1,|u_i|)\,\stackrel{\eqref{082}}{=}\,1\,\stackrel{\eqref{082}}{=}\,\prod_{i\,<\,j}\phi(|u_j|,1)=\prod_{i\,<\,j}\phi(|u_j|,|1_{_{M_i}}|)=\tau^\phi(u,1).\]Therefore, $\tau^\phi$ is a twisting on the direct product.
\end{proof}

Suppose that each $M_i$ is equipped with a twisting $\tau^{\,i}:M_i\times M_i\to\field$ and an action $\rho^i:G_i\to\sym(M_i)$ by a group $G_i$. Observe that, together with the group $\G$, each monoid $M_i$ satisfies Condition~\eqref{015}.\medskip

The following proposition follows immediately from the conditions in \eqref{001}.
\begin{pro}\label{075}
\begin{enumerate}
\item If $\tau$ and $\tau'$ are two twistings on a monoid, then $(\tau\tau')(u,v)=\tau(u,v)\,\tau'(u,v)$ defines another twisting on the same monoid. Consequently, the set of all twistings on a monoid forms a commutative monoid under pointwise multiplication.
\item\label{009}The map $\tau^\times:(M_1\times\cdots\times M_k)\times(M_1\times\cdots\times M_k)\to\field$ defined by $\tau^\times(u,v)=\tau^{\,1}(u_1,v_1)\cdots\tau^{\,k}(u_k,v_k)$, where $u=(u_1,\ldots,u_k)$ and $v=(v_1,\ldots,v_k)$, constitutes a twisting on the direct product $M_1\times\cdots\times M_k$.
\end{enumerate}
\end{pro}
We denote by $\tau_\phi$ the \emph{combined twisting} on the direct product $M_1\times\cdots\times M_k$ defined by $\tau_\phi=\tau^\times\tau^\phi$, where $\tau^\times$ and $\tau^\phi$ are given by Proposition~\ref{075}\eqref{009} and Proposition~\ref{083}, respectively. Consequently, the direct product $M_1\times\cdots\times M_k$, together with the group $\G$ and equipped with the twisting $\tau_\phi$, satisfies Condition~\eqref{015} as well.

\begin{pro}\label{084}
\begin{enumerate}
\item\label{085}There is a group action $\rho^\otimes:G_1\times\cdots\times G_k\to\sym(M_1\times\cdots\times M_k)$, called the \emph{product action}, defined by $\rho^\otimes_\sigma(u_1,\ldots,u_k)=(\rho_{\sigma_1}^{\,1}(u_1),\ldots,\rho_{\sigma_k}^{\,k}(u_k))$ for all $\sigma=(\sigma_1,\ldots,\sigma_k)$. Consequently, for any $u=(u_1,\ldots,u_k)\in M_1\times\cdots\times M_k$, its orbit and stabilizer decompose as\begin{equation}\label{158}(G_1\times\cdots\times G_k)u=G_1u_1\times\cdots\times G_ku_k\quad\text{and}\quad(G_1\times\cdots\times G_k)_{(u)}=(G_1)_{(u_1)}\times\cdots\times(G_k)_{(u_k)}.\,\,\end{equation}
\item\label{086}If each action $\rho^i$ is equipped with a weight function $\omega^{\,i}:G_i\times M_i\to\field^\times$, there is a weight function $\omega^\otimes:(G_1\times\cdots\times G_k)\times(M_1\times\cdots\times M_k)\to\field^\times$, called the \emph{product weight function}, defined by $\omega^\otimes(\sigma,u)=\omega^{\,1}(\sigma_1,u_1)\cdots\omega^{\,k}(\sigma_k,u_k)$ for all $\sigma=(\sigma_1,\ldots,\sigma_k)$ and $u=(u_1,\ldots,u_k)$.
\item\label{087}Assume that each action $\rho^i$ is quasi-multiplicative and $\G$-graded, that is, $|\rho_{\sigma_i}^{\,i}(u_i)|=|u_i|$ for all $\sigma_i\in G_i$ and $u_i\in M_i$. If each $(\tau^{\,i},\omega^{\,i})$ is a multiplicative pair for $\rho^i$, then the product action $\rho^\otimes$ is quasi-multiplicative and $\G$-graded, and $(\tau_\phi,\omega^\otimes)$ is a multiplicative pair.
\end{enumerate}
\end{pro}
\begin{proof}
\eqref{085} The verification that $\rho^\otimes$ defines a valid group action on $M_1\times\cdots\times M_k$ is direct. By definition, we have $(v_1,\dots,v_k)=\rho_\sigma^\otimes(u)$ for some $\sigma=(\sigma_1,\dots,\sigma_k)$ if and only if $v_i=\rho_{\sigma_i}^i(u_i)$ for all $i\in[k]$. Similarly, $\rho_\sigma^\otimes(u)=u$ if and only if each $\rho_{\sigma_i}^i(u_i)=u_i$. This yields the decompositions of the orbit and the stabilizer.

\eqref{086} Let $\sigma=(\sigma_1,\ldots,\sigma_k)$ and $\varsigma=(\varsigma_1,\ldots,\varsigma_k)$ be elements in $G_1\times\cdots\times G_k$, and let $u=(u_1,\ldots,u_k)\in M_1\times\cdots\times M_k$. Since each $\omega^{\,i}$ is a weight function, we have\[\begin{array}{rclcl}
\omega^\otimes(\sigma\varsigma,u)
&=&\omega^{\,1}(\sigma_1\varsigma_1,u_1)\cdots\omega^{\,k}(\sigma_k\varsigma_k,u_k)\\
&\stackrel{\eqref{003}}{=}&\omega^{\,1}(\sigma_1,\rho_{\varsigma_1}^{\,1}(u_1))\,\omega^{\,1}(\varsigma_1,u_1)\cdots\omega^{\,k}(\sigma_k,\rho_{\varsigma_k}^{\,k}(u_k))\,\omega^{\,k}(\varsigma_k,u_k)&=&\omega^\otimes(\sigma,\rho^\otimes_\varsigma(u))\,\omega^\otimes(\varsigma,u).\end{array}\]Furthermore, $\omega^\otimes(\sigma,1)=\omega^{\,1}(\sigma_1,1_{_{M_1}})\cdots\omega^{\,k}(\sigma_k,1_{_{M_k}})=1$. Thus, $\omega^\otimes$ is a well-defined weight function.

\eqref{087}~Let $\sigma=(\sigma_1,\dots,\sigma_k)\in G_1\times\cdots\times G_k$ and $w=(w_1,\dots,w_k)\in M_1\times\cdots\times M_k$. Observe that since the degree map on $M_1\times\cdots\times M_k$ is defined component-wise, for any $u=(u_1,\dots,u_k)$ we have\[|\rho_\sigma^\otimes(u)|=|\rho_{\sigma_1}^1(u_1)|\cdots|\rho_{\sigma_k}^k(u_k)|=|u_1|\cdots|u_k|=|u|,\]which proves that the product action preserves the $\G$-grading. Now consider $u=(u_1,\dots,u_k)$ and $v=(v_1,\dots,v_k)$ such that $uv=w$ and $\tau_\phi(u,v)\neq0$. Since $\tau_\phi(u,v)=\tau^1(u_1,v_1)\cdots\tau^k(u_k,v_k)\,\tau^\phi(u,v)$, we have $u_iv_i=w_i$ and $\tau^i(u_i,v_i)\neq0$ for all $i\in[k]$. By the quasi-multiplicativity of each $\rho^i$ on $w_i$, there exist $\alpha_i,\varsigma_i\in G_i$ depending only on $\sigma_i$ and $w_i$ such that $\rho_{\sigma_i}^i(u_iv_i)=\rho_{\alpha_i}^i(u_i)\,\rho_{\varsigma_i}^i(v_i)$. Setting $\alpha=(\alpha_1,\dots,\alpha_k)$ and $\varsigma=(\varsigma_1,\dots,\varsigma_k)$, which depend only on $\sigma$ and $w$, we obtain $\rho_\sigma^\otimes(uv)=\rho_\alpha^\otimes(u)\,\rho_\varsigma^\otimes(v)$. Since each $\rho^i$ preserves the $\G$-grading, we have $|\rho_{\alpha_j}^j(u_j)|=|u_j|$ and $|\rho_{\varsigma_i}^i(v_i)|=|v_i|$, which implies that $\tau^\phi(\rho_\alpha^\otimes(u),\rho_\varsigma^\otimes(v))=\tau^\phi(u,v)$. Using this invariance and since each pair $(\tau^i,\omega^i)$ is multiplicative, we have\[\begin{array}{rcl}\tau_\phi(u,v)\,\omega^\otimes(\sigma,uv)
&=&\tau^\times(u,v)\,\tau^\phi(u,v)\,\omega^\otimes(\sigma,uv)\\
&=&\tau^\phi(u,v)\,\tau^{\,1}(u_1,v_1)\,\omega^{\,1}(\sigma_1,u_1v_1)\cdots\tau^{\,k}(u_k,v_k)\,\omega^{\,k}(\sigma_k,u_kv_k)\\
&\stackrel{\eqref{028}}{=}&\tau^\phi(u,v)\,\tau^{\,1}(\rho_{\alpha_1}^{\,1}(u_1),\rho_{\varsigma_1}^{\,1}(v_1))\,\omega^{\,1}(\alpha_1,u_1)\,\omega^{\,1}(\varsigma_1,v_1)\cdots\tau^{\,k}(\rho_{\alpha_k}^{\,k}(u_k),\rho_{\varsigma_k}^{\,k}(v_k))\\&&\omega^{\,k}(\alpha_k,u_k)\,\omega^{\,k}(\varsigma_k,v_k)\\
&=&\tau^\phi(\rho_\alpha^\otimes(u),\rho_\varsigma^\otimes(v))\,\tau^\times(\rho_\alpha^\otimes(u),\rho_\varsigma^\otimes(v))\,\omega^\otimes(\alpha,u)\,\omega^\otimes(\varsigma,v)\\
&=&\tau_\phi(\rho_\alpha^\otimes(u),\rho_\varsigma^\otimes(v))\,\omega^\otimes(\alpha,u)\,\omega^\otimes(\varsigma,v).
\end{array}\]Therefore, $(\tau_\phi,\omega^\otimes)$ is a multiplicative pair for the product action.
\end{proof}

Assume now that each monoid $M_i$ is locally finite and graded, both with respect to the twisting $\tau^{\,i}$, with degree map $\deg:M_i\to\N$ and grading $\{M_i^d\mid d\in\N\}$, where each $M_i^d$ intersects only finitely many fibres $\ubar{M}_{i,x}$. This implies that each $M_i$ satisfies Conditions~\eqref{016} and~\eqref{017}. Similarly to the $\G$-grading, their direct product $M_1\times\cdots\times M_k$ naturally inherits an $\N$-grading with degree map defined by $\deg(u_1,\ldots,u_k)=\deg(u_1)+\cdots+\deg(u_k)$. Accordingly, the homogeneous component of degree $d\in\N$ is given by the disjoint union\begin{equation}\label{091}\bigsqcup_{d_1+\cdots+d_k=d}M_1^{\,d_1}\times\cdots\times M_k^{\,d_k}.\end{equation}Furthermore, we assume that each action $\rho^i$ is quasi-multiplicative, bigraded, and equipped with a weight function $\omega^{\,i}$, such that $(\tau^{\,i},\omega^{\,i})$ is a multiplicative pair. This implies that each group $G_i$ satisfies Condition~\eqref{018}. Additionally, we assume that for each $i\in[k]$ there exists a system of representatives $T_i\subseteq M_i$ for the $\omega^i$-admissible $G_i$-orbits such that $T_i\cap M_i^d$ is finite for all $d\in\N$, so that Condition~\eqref{051} holds for all $i$, and so, by Theorem~\ref{050}, each algebra $\funcs{\field}{M_i}^{G_i}$ is endowed with the monomial basis $\{m_t^{\,i}\mid t\in T_i\}$.

\begin{pro}\label{092}
Consider the monoid $M_1\times\cdots\times M_k$ equipped with the combined twisting $\tau_\phi$, and the group $G_1\times\cdots\times G_k$ endowed with the product action $\rho^\otimes$ and the product weight function $\omega^\otimes$. Then
\begin{enumerate}
\item\label{093}The monoid $M_1\times\cdots\times M_k$ is locally finite and each homogeneous component of its $\N$-grading intersects only finitely many fibres of its $\G$-grading. Consequently, $M_1\times\cdots\times M_k$ satisfies Conditions~\eqref{016}~and~\eqref{017}.
\item\label{094}The product action $\rho^\otimes$ is bigraded. Consequently, $G_1\times\cdots\times G_k$ satisfies Condition~\eqref{018}.
\item\label{095}The direct product $T_1\times\cdots\times T_k$ is a system of representatives for the $\omega^\otimes$-admissible $(G_1\times\cdots\times G_k)$-orbits, and each fibre of the $\N$-grading of $M_1\times\cdots\times M_k$ intersects only finitely many of these orbits. Consequently, $T_1\times\cdots\times T_k$ satisfies Condition~\eqref{051}.
\end{enumerate}
\end{pro}
\begin{proof}
\eqref{093} A factorization $vw$ of $u=(u_1,\ldots,u_k)$ with $\tau_\phi(v,w)\neq0$ in $M_1\times\cdots\times M_k$ corresponds to component-wise factorizations of each $u_i=v_iw_i$ with $\tau^{\,i}(v_i,w_i)\neq0$ in $M_i$. Since each $M_i$ is locally finite with respect to $\tau^{\,i}$, each $u_i$ admits only finitely many factorizations, implying that $u$ admits only finitely many factorizations. Thus, $M_1\times\cdots\times M_k$ is locally finite with respect to $\tau_\phi$. Since $M_1\times\cdots\times M_k$ is $\G$-graded as in~\eqref{090}, it satisfies Condition~\eqref{016}. Furthermore, for each fixed $d\in\N$, the disjoint union in~\eqref{091} runs over finitely many tuples $(d_1,\dots,d_k)$ with $d_1+\cdots+d_k=d$, and for each tuple, every factor $M_i^{d_i}$ intersects only finitely many $\G$-fibres $\ubar{M}_{i,x_i}$. Hence, only finitely many product $\G$-degrees $x=x_1\cdots x_k$ occur in $(M_1\times\cdots\times M_k)^d$, so Condition~\eqref{017} holds.

\eqref{094} Let $\sigma=(\sigma_1,\ldots,\sigma_k)\in G_1\times\cdots\times G_k$ and $u=(u_1,\ldots,u_k)\in M_1\times\cdots\times M_k$. By definition of the degree map on the direct product and the fact that each $\rho^i$ is graded, we have\[\deg(\rho^\otimes_\sigma(u))=\deg(\rho_{\sigma_1}^{\,1}(u_1))+\cdots+\deg(\rho_{\sigma_k}^{\,k}(u_k))=\deg(u_1)+\cdots+\deg(u_k)=\deg(u).\]Therefore, $\rho^\otimes$ preserves the $\N$-grading. This result, combined with Proposition~\ref{084}, implies that $\rho^\otimes$ is bigraded and $(\tau_\phi,\omega^\otimes)$ is a multiplicative pair, ensuring that Condition~\eqref{018} is satisfied.

\eqref{095}~By Proposition~\ref{084}\eqref{085}, the orbit of any $u=(u_1,\ldots,u_k)\in M_1\times\cdots\times M_k$ is $G_1u_1\times\cdots\times G_ku_k$, and its stabilizer is $(G_1)_{(u_1)}\times\cdots\times(G_k)_{(u_k)}$. By definition, $\omega^\otimes(\sigma,u)=1$ for all $\sigma$ in the stabilizer if and only if $\omega^{\,i}(\sigma_i,u_i)=1$ for all $\sigma_i\in (G_i)_{(u_i)}$. This proves that $u$ is $\omega^\otimes$-admissible if and only if each $u_i$ is $\omega^{\,i}$-admissible. Since each $T_i$ contains exactly one representative for each admissible $G_i$-orbit, $T_1\times\cdots\times T_k$ contains exactly one representative for each admissible orbit of the direct product. By definition, each homogeneous component $(M_1\times\cdots\times M_k)^d$ is the finite disjoint union of products $M_1^{d_1}\times\cdots\times M_k^{d_k}$ with $d_1+\cdots+d_k=d$. By~\eqref{051} on each $M_i$, each intersection $T_i\cap M_i^{d_i}$ is finite, and therefore $(T_1\times\cdots\times T_k)\cap(M_1\times\cdots\times M_k)^d$ is finite, ensuring that Condition~\eqref{051} holds. %By~\eqref{093}, any intersection of a $\G$-fibre with an $\N$-homogeneous component is a finite union of products of local intersections $\ubar{M}_{i,x_i}\cap M_i^{\,d_i}$. By~\eqref{051} on each $M_i$, these local intersections contain only finitely many elements of $T_i$. Therefore, their finite union contains only finitely many elements of $T_1\times\cdots\times T_k$, ensuring that Condition~\eqref{051} holds.
\end{proof}

Proposition~\ref{092} implies that the monoid $M_1\times\cdots\times M_k$ equipped with the combined twisting $\tau_\phi$, and the group $G_1\times\cdots\times G_k$ endowed with the product action $\rho^\otimes$ and weight function $\omega^\otimes$, satisfy conditions \eqref{015} to \eqref{051} with respect to the system of representatives $T_1\times\cdots\times T_k$. Consequently, by Theorem~\ref{050}, the twisted algebra of invariant functions $\funcs{\field}{M_1\times\cdots\times M_k}^{G_1\times\cdots\times G_k}$ admits the monomial basis $\{m_t^\otimes\mid t\in T_1\times\cdots\times T_k\}$.

\begin{pro}\label{159}
Consider the monoid $M_1\times\cdots\times M_k$ equipped with the combined twisting $\tau_\phi$, and the group $G_1\times\cdots\times G_k$ endowed with the product action $\rho^\otimes$ and the product weight function $\omega^\otimes$. Let $t=(t_1,\dots,t_k)\in T_1\times\dots\times T_k$ and $u=(u_1,\dots,u_k)\in M_1\times\dots\times M_k$. Then\[m_t^\otimes(u)=m_{t_1}^1(u_1)\cdots m_{t_k}^k(u_k),\]where each $m_{t_i}^i$ is the monomial invariant function associated with $t_i\in T_i$.
\end{pro}
\begin{proof}
Assume that $u\in(G_1\times\cdots\times G_k)t$. By Proposition~\ref{084}\eqref{085}, there exist elements $\sigma_i\in G_i$ such that $u_i=\rho_{\sigma_i}^i(t_i)$, which implies that $u=\rho_\sigma^\otimes(t)$ for $\sigma=(\sigma_1,\cdots,\sigma_k)\in G_1\times\cdots\times G_k$. By definition of the monomial invariant function, we have $m_t^\otimes(u)=\omega^\otimes(\sigma,t)$. Then\[m_t^\otimes(u)=\omega^\otimes(\sigma,t)=\omega^1(\sigma_1,t_1)\cdots\omega^k(\sigma_k,t_k)=m_{t_1}^1(u_1)\cdots m_{t_k}^k(u_k).\]Assume now that $u\notin(G_1\times\cdots\times G_k)t$. By definition, we have $m_t^\otimes(u)=0$. On the other hand, Proposition~\ref{084}\eqref{085} implies that there exists $j\in[k]$ such that $u_j\notin G_j t_j$. Consequently, $m_{t_j}^j(u_j)=0$. Therefore, $m_{t_1}^1(u_1)\cdots m_{t_k}^k(u_k)=0=m_t^\otimes(u)$, which concludes the proof.
\end{proof}

\begin{thm}\label{088}
The twisted algebra of invariant functions $\funcs{\field}{M_1\times\cdots\times M_k}^{G_1\times\cdots\times G_k}$ is isomorphic to the tensor product of algebras $\funcs{\field}{M_1}^{G_1}\otimes_\phi\cdots\otimes_\phi\funcs{\field}{M_k}^{G_k}$.
\end{thm}
\begin{proof}
We define the following linear map on basis elements\[\Phi:\funcs{\field}{M_1\times\cdots\times M_k}^{G_1\times\cdots\times G_k}\to\funcs{\field}{M_1}^{G_1}\otimes_\phi\cdots\otimes_\phi\funcs{\field}{M_k}^{G_k}\quad\text{such that}\quad\Phi(m_{(u_1,\ldots,u_k)}^\otimes)=m_{u_1}^{\,1}\otimes\cdots\otimes m_{u_k}^{\,k}.\]Since $\Phi$ maps a basis bijectively to a basis, it is a vector space isomorphism. To prove it is an algebra homomorphism, let $u,v\in T$. By Proposition~\ref{052}, the structure constant of the product $m_u^\otimes\cdot m_v^\otimes$ is given by\[c_{u,v}^t=\sum_{\substack{(x,y)\in Gu\times Gv\\xy=t}}\tau_\phi(x,y)\,m_u^\otimes(x)\,m_v^\otimes(y)\,\,\,\,\refeq{158}\sum_{\substack{(x_i,y_i)\in G_iu_i\times G_iv_i\\x_iy_i=t_i}}\!\!\!\tau^\phi(x,y)\,\tau^\times(x,y)\prod_{i=1}^km_{u_i}^{\,i}(x_i)\,m_{v_i}^{\,i}(y_i),\]where $u=(u_1,\ldots,u_k)$, $v=(v_1,\ldots,v_k)$, $t=(t_1,\ldots,t_k)$, $x=(x_1,\ldots,x_k)$ and $y=(y_1,\ldots,y_k)$. Since each $\rho^i$ is $\G$-graded, the elements $x_j\in G_ju_j$ and $y_i\in G_iv_i$ satisfy $|x_j|=|u_j|$ and $|y_i|=|v_i|$. Consequently, we have\[\tau^\phi(x,y)=\prod_{i\,<\,j}\phi(|x_j|,|y_i|)=\prod_{i\,<\,j}\phi(|u_j|,|v_i|)=\tau^\phi(u,v).\]Factoring out this constant, the sum splits by coordinates\[c_{u,v}^t=\,\tau^\phi(u,v)\prod_{i=1}^k\Bigg(\,\sum_{\substack{(x_i,y_i)\in G_iu_i\times G_iv_i\\x_iy_i=t_i}}\!\!\!\!\tau^{\,i}(x_i,y_i)\,m_{u_i}^{\,i}(x_i)\,m_{v_i}^{\,i}(y_i)\Bigg)=\tau^\phi(u,v)\prod_{i=1}^k c_{u_i,v_i}^{t_i}.\]Applying $\Phi$ yields\[\Phi(m_u^\otimes\cdot m_v^\otimes)=\tau^\phi(u,v)\sum_{t\in T}\Bigg(\prod_{i=1}^k c_{u_i,v_i}^{t_i}\Bigg)m_{t_1}^{\,1}\otimes\cdots\otimes m_{t_k}^{\,k}.\]
Now, multiplying the images in the tensor product, and recalling that $|m_{u_i}^{\,i}|=|u_i|$, we obtain\[\begin{array}{rcl}
\Phi(m_u^\otimes)\,\Phi(m_v^\otimes)
&=&(m_{u_1}^{\,1}\otimes\cdots\otimes m_{u_k}^{\,k})\,(m_{v_1}^{\,1}\otimes\cdots\otimes m_{v_k}^{\,k})\\[0.1cm]
&=&\displaystyle\tau^\phi(u,v)(m_{u_1}^{\,1}m_{v_1}^{\,1}\otimes\cdots\otimes m_{u_k}^{\,k}m_{v_k}^{\,k})\\[0.1cm]
&=&\displaystyle\tau^\phi(u,v)\Bigg(\sum_{t_1\in T_1}c_{u_1,v_1}^{t_1}m_{t_1}^{\,1}\Bigg) \otimes\cdots\otimes\Bigg(\sum_{t_k\in T_k}c_{u_k,v_k}^{t_k} m_{t_k}^{\,k}\Bigg)\\
&=&\displaystyle\tau^\phi(u,v)\sum_{t\in T}\Bigg(\prod_{i=1}^kc_{u_i,v_i}^{t_i}\Bigg)\,m_{t_1}^{\,1}\otimes\cdots\otimes m_{t_k}^{\,k}\,\,\,=\,\,\,\Phi(m_u^\otimes\cdot m_v^\otimes).
\end{array}\]Given $t=(t_1,\ldots,t_k)\in T$, the basis element $m_t^\otimes$ is homogeneous of bidegree $(|t|,\deg(t))=(|t_1|\cdots|t_k|,\allowbreak\deg(t_1)+\cdots+\deg(t_k))$. By Remark~\ref{081} and~\eqref{096}, in the tensor product, the $\G$-degree of a pure tensor is the product of the local $\G$-degrees, and the $\N$-degree is the sum of the local $\N$-degrees. Therefore, the bidegree of $\Phi(m_t^\otimes)=m_{t_1}^{\,1}\otimes\cdots\otimes m_{t_k}^{\,k}$ evaluates precisely to $(|t_1|\cdots|t_k|,\deg(t_1)+\cdots+\deg(t_k))$. Since $\Phi$ preserves the bidegree of the basis elements, it is an isomorphism of bigraded algebras.
\end{proof}

Specializing to the case where all factors coincide immediately yields the following result.
\begin{crl}\label{097}
For every positive integer $k$, there exists an isomorphism of bigraded algebras\[\Phi_k:\funcs{\field}{M\times\cdots\times M}^{G\times\cdots\times G}\to\underbrace{\funcs{\field}{M}^G\otimes_\phi\cdots\otimes_\phi\funcs{\field}{M}^G}_{k\text{ times}},\]defined on the monomial basis elements by $\Phi_k(m_{(t_1,\ldots,t_k)}^\otimes)=m_{t_1}\otimes\cdots\otimes m_{t_k}$ for all $t_1,\ldots,t_k\in T$.
\end{crl}

\subsection{Diagonal invariants}\label{183}

In the previous subsection, we considered the product action of multiple groups acting independently on a direct product of monoids. We now turn our attention to the scenario where a single group acts simultaneously on direct and free products. In this subsection, we formalize the lifting of local twistings and the construction of diagonal actions on multiple products, establishing the framework required to obtain well-defined twisted algebras.

As in Subsection~\ref{120}, let $M_1,\ldots,M_k$ be monoids, each equipped with a $\G$-grading $\{\ubar{M}_{i,x}\mid x\in\G\}$. Furthermore, assume that each $M_i$ is endowed with a twisting $\tau^{\,i}:M_i\times M_i\to\field$ with respect to which it is locally finite and graded with $\N$-grading $\{M_i^d\mid d\in\N\}$, such that each homogeneous component $M_i^d$ intersects only finitely many fibres $\ubar{M}_{i,x}$. In addition, assume that each $M_i$ is endowed with a multiplicative bigraded action $\rho^i:G\to\sym(M_i)$ and a weight function $\omega^{\,i}:G\times M_i\to\field^\times$ such that $(\tau^{\,i},\omega^{\,i})$ forms a multiplicative pair. Finally, let $T_i\subseteq M_i$ be a system of representatives for the $\omega^i$-admissible $G$-orbits such that $T_i\cap M_i^d$ is finite for all $d\in\N$. Consequently, Conditions~\eqref{015} to~\eqref{051} hold for each $i$. Building on this, by Propositions~\ref{083},~\ref{075}~and~\ref{092}\eqref{093}, the direct product $M_1\times\cdots\times M_k$ can be equipped with the twisting $\tau_\phi=\tau^\times\tau^\phi$ and the gradings in \eqref{090} and~\eqref{091}; that is, $|(u_1,\ldots,u_k)|=|u_1|\cdots|u_k|$ and $\deg(u_1,\ldots,u_k)=\deg(u_1)+\cdots+\deg(u_k)$, such that Conditions~\eqref{015}~to~\eqref{017} hold for $M_1\times\cdots\times M_k$.

\begin{pro}\label{076}
\begin{enumerate}
\item\label{077}There is a multiplicative bigraded group action $\rho^\times:G\to\sym(M_1\times\cdots\times M_k)$, called the \emph{diagonal action}, defined by $\rho^\times_\sigma(u_1,\ldots,u_k)=(\rho_\sigma^{\,1}(u_1),\ldots,\rho_\sigma^{\,k}(u_k))$ for all $\sigma\in G$.
\item\label{078}There is a weight function $\omega^\times:G\times(M_1\times\cdots\times M_k)\to\field^\times$, called the \emph{diagonal weight function}, defined by $\omega^\times(\sigma,(u_1,\ldots,u_k))=\omega^{\,1}(\sigma,u_1)\cdots\omega^{\,k}(\sigma,u_k)$ for all $\sigma\in G$, such that $(\tau_\phi,\omega^\times)$ is a multiplicative pair.
\end{enumerate}\medskip
\noindent Consequently, Condition~\eqref{018} holds for the group $G$ and the monoid $M_1\times\cdots\times M_k$.
\end{pro}
\begin{proof}
The fact that $\rho^\times$ is a multiplicative bigraded action and $\omega^\times$ is a well-defined weight function follows directly from their component-wise definitions and the hypotheses on $M_i$. To verify that $(\tau_\phi,\omega^\times)$ is a multiplicative pair, let $\sigma\in G$ and consider $u=(u_1,\dots,u_k)$ and $v=(v_1,\dots,v_k)$ in $M_1\times\cdots\times M_k$ with $\tau_\phi(u,v)\neq0$. Since $\rho^\times$ is $\G$-graded, by definition of $\tau^\phi$, it follows that $\tau^\phi(\rho^\times_\sigma(u),\rho^\times_\sigma(v))=\tau^\phi(u,v)$. Then, because each $(\tau^{\,i},\omega^{\,i})$ is a multiplicative pair, we have\[\begin{array}{rcl}
\tau_\phi(u,v)\,\omega^\times(\sigma,uv)
&=&\tau^\phi(u,v)\,\tau^{\,1}(u_1,v_1)\,\omega^{\,1}(\sigma,u_1v_1)\cdots\tau^{\,k}(u_k,v_k)\,\omega^{\,k}(\sigma,u_kv_k)\\[0.1cm]
&\stackrel{\eqref{028}}{=}&\tau^\phi(\rho^\times_\sigma(u),\rho^\times_\sigma(v))\,\tau^{\,1}(\rho_\sigma^{\,1}(u_1),\rho_\sigma^{\,1}(v_1))\,\omega^{\,1}(\sigma,u_1)\,\omega^{\,1}(\sigma,v_1)\cdots\\[0.1cm]
&&\cdots\tau^{\,k}(\rho_\sigma^{\,k}(u_k),\rho_\sigma^{\,k}(v_k))\,\omega^{\,k}(\sigma,u_k)\,\omega^{\,k}(\sigma,v_k)\\[0.15cm]
&=&\tau_\phi(\rho^\times_\sigma(u),\rho^\times_\sigma(v))\,\omega^\times(\sigma,u)\,\omega^\times(\sigma,v).
\end{array}\]Therefore, $(\tau_\phi,\omega^\times)$ is a multiplicative pair for the diagonal action $\rho^\times$.
\end{proof}
Proposition~\ref{076}, together with Proposition~\ref{092}\eqref{093}, ensures the existence of the \emph{twisted algebra of diagonal invariants} $\funcs{\field}{M_1\times\cdots\times M_k}^G$. Furthermore, if $T^\times \subseteq M_1\times\cdots\times M_k$ is a system of representatives for the $\omega^\times$-admissible diagonal $G$-orbits satisfying Condition~\eqref{051}, Theorem~\ref{050} guarantees that this algebra admits the monomial basis $\{m_t^\times\mid t\in T^\times\}$.\medskip

Having established that the direct product of monoids satisfies the necessary conditions to yield a well-defined twisted algebra of invariants, we now turn our attention to free products.

In general, if $M_1,\ldots,M_k$ are $\P$-graded monoids with respect to the trivial twisting, with respective degree maps $\nu_1,\ldots,\nu_k$, then these extend to a degree map $\nu^*=\nu_1*\cdots*\nu_k$, which endows the free product $M_1*\cdots*M_k$ with a $\P$-grading \cite[Proposition~8.2.1]{How95}. Explicitly, the $\P$-degree of an alternating product $u=u_1\cdots u_q$ with $u_j\in M_{i_j}\setminus\{1\}$ is given by $\nu^*(u)=\nu_{i_1}(u_1)\cdots\nu_{i_q}(u_q)$. In particular, the free product $M_1*\cdots*M_k$ naturally inherits a $\G$-grading given by $|u|=|u_1|\cdots|u_q|$.

\begin{pro}\label{067}
There exists a unique twisting $\tau^*:(M_1*\cdots*M_k)\times(M_1*\cdots*M_k)\to\field$, called the \emph{free twisting}, such that its restriction to each $M_i\times M_i$ coincides with $\tau^{\,i}$, and for any $u,v\in (M_1*\cdots*M_k)\setminus\{1\}$, if the central factors of $u$ and $v$ belong to distinct monoids, then $\tau^*(u,v)=1$. Consequently, Condition~\eqref{015} holds for the free product $M_1*\cdots*M_k$.
\end{pro}
\begin{proof}
We first establish uniqueness. Suppose $\gamma$ is a twisting on $M_1*\cdots*M_k$ satisfying the conditions of $\tau^*$. If $u=1$ or $v=1$, Condition~\eqref{001} forces $\gamma(u,v)=1$. If the central factors of $u$ and $v$ belong to distinct monoids, then $\gamma(u,v)=1$. Otherwise, we can decompose $u=u'x$ and $v=yv'$ such that $x,y\in M_i\setminus\{1\}$ for some $i\in[k]$. Because $\gamma(u',x)=1=\gamma(y,v')$, applying~\eqref{001} repeatedly yields\[\gamma(u,v)=\gamma(u',x)\,\gamma(u'x,yv')\,\gamma(y,v')\refeq{001}\gamma(u',xyv')\,\gamma(x,yv')\,\gamma(y,v')\refeq{001}\gamma(u',xyv')\,\gamma(x,y)\,\gamma(xy,v'),\]where $\gamma(x,y)=\tau^i(x,y)$. If $\tau^i(x,y)=0$, this product vanishes, so $\gamma(u,v)=0=\tau^i(x,y)$. If $\tau^i(x,y)\neq0$, then $\deg(xy)=\deg(x)+\deg(y)>0$ because $M_i^0=\{1\}$, so $xy\in M_i\setminus\{1\}$; thus $\gamma(u',xyv')=1=\gamma(xy,v')$ since the central factors belong to distinct monoids, yielding again $\gamma(u,v)=\tau^i(x,y)$. Thus, the values of $\gamma$ are uniquely determined by the local twistings $\tau^i$. To establish existence, for each $u,v\in M_1*\cdots*M_k$, we define\begin{equation}\label{072}
\tau^*(u,v)=\begin{cases}\tau^i(x,y)&\text{if the central factors $x,y$ belong to }M_i\setminus\{1\},\\1&\text{otherwise.}\end{cases}\end{equation}By construction, this map satisfies the stated properties. To complete the proof, we verify that $\tau^*$ satisfies the $2$-cocycle condition $\tau^*(u,v)\,\tau^*(uv,w)=\tau^*(u,vw)\,\tau^*(v,w)$. If any element is $1$, the equality holds trivially. Consider $u,v,w\in(M_1*\cdots*M_k)\setminus\{1\}$,
and let $q$ be the length of the alternating product of $v$. If $q\geq2$, we can write $v=yv'z$, where $y\in M_i\setminus\{1\}$ and $z\in M_j\setminus\{1\}$ are the first and last factors in the alternating product of $v$. Let $x$ be the last factor of $u$ and $t$ be the first factor of $w$. If $x\in M_i$ and $\tau^i(x,y)=0$, then $\tau^*(u,v)=0$, while on the right-hand side either $t\in M_j$ and $\tau^j(z,t)=0$, so $\tau^*(v,w)=0$, or else $vw$ still begins with $y$, since $zt\neq1$ when $t\in M_j$ and $\tau^j(z,t)\neq0$, which forces $\tau^*(u,vw)=\tau^i(x,y)=0$; hence both sides vanish. Symmetrically, both sides vanish if $t\in M_j$ and $\tau^j(z,t)=0$. Otherwise, the reduced word of $uv$ still ends in $z$ and the reduced word of $vw$ still begins with $y$, so\[\tau^*(u,v)\,\tau^*(uv,w)=\tau^*(u,v)\,\tau^*(z,w)=\tau^*(u,v)\,\tau^*(v,w)=\tau^*(u,y)\,\tau^*(v,w)=\tau^*(u,vw)\,\tau^*(v,w).\]If $q=1$, then $v\in M_j\setminus\{1\}$ for some $j\in[k]$. Let $x$ and $t$ denote the central factors of $u$ and $w$ respectively. If $x,t\notin M_j$, then $\tau^*(u,v)\,\tau^*(uv,w)=1=\tau^*(u,vw)\,\tau^*(v,w)$. If $x\in M_j$ and $t\notin M_j$, then $vw$ begins with $v\in M_j$, so $\tau^*(u,vw)\,\tau^*(v,w)=\tau^j(x,v)\cdot1=\tau^j(x,v)$, while $\tau^*(u,v)\,\tau^*(uv,w)=\tau^j(x,v)\,\tau^*(uv,w)$ also equals $\tau^j(x,v)$ as $\tau^j(x,v)=0$ when $xv=1$, and $xv\in M_j\setminus\{1\}$ gives $\tau^*(uv,w)=1$ when $\tau^j(x,v)\neq0$. The case $x\notin M_j$ and $t\in M_j$ is symmetric. Finally, if $x,v,t\in M_j$, the $2$-cocycle identity for $\tau^j$ gives\[\tau^*(u,v)\,\tau^*(uv,w)=\tau^j(x,v)\,\tau^j(xv,t)\refeq{001}\tau^j(x,vt)\,\tau^j(v,t)=\tau^*(u,vw)\,\tau^*(v,w),\]since $\tau^*(u',1)=1=\tau^*(1,w')$ coincides with $\tau^j(1,t)=1=\tau^j(x,1)$ even if $xv=1$ or $vt=1$. This completes the proof.
\end{proof}

\begin{pro}\label{184}
The free product $M_1*\cdots*M_k$ is a graded monoid with respect to $\tau^*$, with degree map $\deg:M_1*\cdots*M_k\to\N$ given on alternating products $u=u_1\cdots u_q$, with $u_j\in M_{i_j}\setminus\{1\}$ and $i_j\neq i_{j+1}$, by\[\deg(u)=\deg(u_1)+\cdots+\deg(u_q),\]and $\deg(1)=0$. Moreover, $(M_1*\cdots*M_k)^0=\{1\}$.
\end{pro}
\begin{proof}
By definition, $\deg(1)=0$. Since each $M_i^0=\{1\}$, for any $u=u_1\cdots u_q\neq1$ we have $\deg(u_j)\geq1$ for all $j$, so $\deg(u)\geq q\geq1$. Hence $(M_1*\cdots*M_k)^0=\{1\}$. Let $u,v\in M_1*\cdots*M_k$ such that $\tau^*(u,v)\neq 0$. If $u=1$ or $v=1$, the equality $\deg(uv)=\deg(u)+\deg(v)$ is immediate. Assume $u,v\neq 1$, and write them as $u=u'x$ and $v=yv'$, where $x\in M_i\setminus\{1\}$ and $y\in M_j\setminus\{1\}$ are their central factors. If $i\neq j$, then $uv=u'xyv'$ is already an alternating product. Thus, $\deg(uv)=\deg(u')+\deg(x)+\deg(y)+\deg(v')=\deg(u)+\deg(v)$. If $i=j$, then by definition $\tau^*(u,v)=\tau^{\,i}(x,y)\neq 0$. Since $M_i$ is $\N$-graded with respect to $\tau^{\,i}$, it follows that $\deg(xy)=\deg(x)+\deg(y)$. Moreover, since $M_i^0=\{1\}$, we have $xy\neq 1$. Thus, the alternating product of $uv$ is $u'(xy)v'$, and we obtain $\deg(uv)=\deg(u')+\deg(xy)+\deg(v')=\deg(u')+\deg(x)+\deg(y)+\deg(v')=\deg(u)+\deg(v)$. Therefore, $M_1*\cdots*M_k$ is graded with respect to $\tau^*$.
\end{proof}

Given $t\in M_1*\cdots*M_k$, any factorization $uv=t$ with $\tau^*(u,v)\neq 0$ is either a concatenation of alternating subproducts of $t$, or involves a central reduction $xy=w$ for some factor $w$ of $t$ with $\tau^{\,i}(x,y)\neq 0$. Since $t$ has finite length and each $M_i$ is locally finite with respect to $\tau^{\,i}$, there are only finitely many such pairs $(u,v)$. Thus, $M_1*\cdots*M_k$ is locally finite with respect to $\tau^*$. Moreover, since $M_i^0=\{1\}$, words of $\N$-degree $d$ have length at most $d$. As their letters come from components that intersect only finitely many $\G$-fibres, their product in $\G$ takes only finitely many values. Consequently, Conditions~\eqref{016} and~\eqref{017} hold for the free product.

\begin{pro}\label{068}
\begin{enumerate}
\item\label{069} Assume that each $\rho^i:G\to\Aut(M_i)$ acts by monoid automorphisms. Then there is a unique bigraded group action $\rho^*:G\to\Aut(M_1*\cdots*M_k)$ by monoid automorphisms, called the \emph{free diagonal action}, such that $\rho_\sigma^*(u)=\rho_\sigma^i(u)$ for all $i\in[k]$ and $u\in M_i$.
\item\label{070}There is a unique weight function $\omega^*:G\times(M_1*\cdots*M_k)\to\field^\times$ for $\rho^*$, called the \emph{free diagonal weight function}, such that its restriction to each $G\times M_i$ coincides with $\omega^{\,i}$, and $\omega^*(\sigma,uv)=\omega^*(\sigma,u)\,\omega^*(\sigma,v)$ whenever $u,v\in(M_1*\cdots*M_k)\setminus\{1\}$ and their central factors belong to distinct monoids.
\item\label{071}$(\tau^*,\omega^*)$ is a multiplicative pair for the free diagonal action $\rho^*$.
\end{enumerate}\medskip
\noindent Consequently, Condition~\eqref{018} holds for the group $G$ and the free product $M_1*\cdots*M_k$.
\end{pro}
\begin{proof}
\eqref{069}~Since each $M_i$ embeds as a submonoid of $M_1*\cdots*M_k$, we may, for every $\sigma\in G$, regard $\rho^i_\sigma$ as a monoid monomorphism $M_i\hookrightarrow M_1*\cdots*M_k$. By \cite[Proposition~8.2.1]{How95}, there exists a unique endomorphism $\rho^*_\sigma:M_1*\cdots*M_k\to M_1*\cdots*M_k$ such that $\rho^*_\sigma(u)=\rho^i_\sigma(u)$ for all $i\in[k]$ and $u\in M_i$. Note that for $\sigma,\varsigma\in G$ and $u\in M_i$, since $\rho^i_\varsigma(u)\in M_i$, we have $(\rho^*_\sigma\circ\rho^*_\varsigma)(u)=\rho^*_\sigma(\rho^i_\varsigma(u))=\rho^i_\sigma(\rho^i_\varsigma(u))=\rho^i_{\sigma\varsigma}(u)=\rho^*_{\sigma\varsigma}(u)$. Because $M_1*\cdots*M_k$ is generated by the submonoids $M_i$, the uniqueness of the extension forces $\rho^*_\sigma\circ\rho^*_\varsigma=\rho^*_{\sigma\varsigma}$. Similarly, $\rho^*_1(u)=\rho^i_1(u)=u$ forces $\rho^*_1=\id$. This implies that each $\rho^*_\sigma$ is invertible with inverse $\rho^*_{\sigma^{-1}}$, and so the map $\rho^*:G\to\Aut(M_1*\cdots*M_k)$ defined by $\sigma\mapsto\rho^*_\sigma$ is a well-defined multiplicative group action. Furthermore, because each $\rho^i$ preserves the bigrading on $M_i$, computing degrees factor-by-factor over alternating products ensures that $\rho^*$ preserves the induced bigrading on $M_1*\cdots*M_k$.

\eqref{070}~For every $\sigma\in G$, we set $\omega^*(\sigma,1)=1$. Let $u\in(M_1*\cdots*M_k)\setminus\{1\}$ with $u=u_1\cdots u_q$ its representation as an alternating product, where $u_r\in M_{i_r}\setminus\{1\}$. We define $\omega^*(\sigma,u)=\omega^{\,i_1}(\sigma,u_1)\cdots\omega^{\,i_q}(\sigma,u_q)$, which ensures uniqueness. To verify existence, we must show that this $\omega^*$ satisfies the 1-cocycle condition in \eqref{003}. Since $\rho^*_\varsigma$ acts component-wise and preserves each submonoid $M_i$, the element $\rho^*_\varsigma(u)=\rho^*_\varsigma(u_1)\cdots\rho^*_\varsigma(u_q)$ is precisely the alternating product of $\rho^*_\varsigma(u)$. Therefore,\[\begin{array}{rclll}
\omega^*(\sigma\varsigma,u)&=&\omega^{\,i_1}(\sigma\varsigma,u_1)\cdots\omega^{\,i_q}(\sigma\varsigma,u_q)\\
&\stackrel{\eqref{003}}{=}&\omega^{\,i_1}(\sigma,\rho^*_\varsigma(u_1))\,\omega^{\,i_1}(\varsigma,u_1)\cdots\omega^{\,i_q}(\sigma,\rho^*_\varsigma(u_q))\,\omega^{\,i_q}(\varsigma,u_q)&=&\omega^*(\sigma,\rho^*_\varsigma(u))\,\omega^*(\varsigma,u).
\end{array}\]

\eqref{071}~We must verify that $\tau^*(u,v)\,\omega^*(\sigma,uv)=\tau^*(\rho_\sigma^*(u),\rho_\sigma^*(v))\,\omega^*(\sigma,u)\,\omega^*(\sigma,v)$ for all $u,v\in M_1*\cdots*M_k$ with $\tau^*(u,v)\neq0$. If $u=1$ or $v=1$, or if their central factors belong to distinct monoids, then $\tau^*(u,v)=1$ and $\tau^*(\rho^*_\sigma(u),\rho^*_\sigma(v))=1$. In this case, $\omega^*(\sigma,uv)=\omega^*(\sigma,u)\,\omega^*(\sigma,v)$ by definition, so the equality holds trivially. Otherwise, decompose $u=u'x$ and $v=yv'$ such that $x,y\in M_i\setminus\{1\}$ are the central factors. Since $\tau^*(u,v)=\tau^i(x,y)\neq0$ and $M_i^0=\{1\}$, we have $\deg(xy)=\deg(x)+\deg(y)>0$ and thus $xy\neq1$, which gives $\tau^*(u,v)=\tau^i(x,y)$ and $\tau^*(\rho^*_\sigma(u),\rho^*_\sigma(v))=\tau^{\,i}(\rho^*_\sigma(x),\rho^*_\sigma(y))$. Moreover, the alternating reduced word for $uv$ is $u'(xy)v'$, so $\omega^*(\sigma,uv)=\omega^*(\sigma,u')\,\omega^{\,i}(\sigma,xy)\,\omega^*(\sigma,v')$. Using that $(\tau^{\,i},\omega^{\,i})$ is a multiplicative pair on $M_i$, we obtain\[\begin{array}{rcl}
\tau^*(u,v)\,\omega^*(\sigma,uv)&=&\omega^*(\sigma,u')\,\tau^{\,i}(x,y)\,\omega^{\,i}(\sigma,xy)\,\omega^*(\sigma,v')\\
&\stackrel{\eqref{028}}{=}&\omega^*(\sigma,u')\,\tau^{\,i}(\rho^*_\sigma(x),\rho^*_\sigma(y))\,\omega^{\,i}(\sigma,x)\,\omega^{\,i}(\sigma,y)\,\omega^*(\sigma,v')\\[0.08cm]
&=&\tau^*(\rho^*_\sigma(u),\rho^*_\sigma(v))\,\omega^*(\sigma,u)\,\omega^*(\sigma,v).
\end{array}\]This completes the proof.
\end{proof}

\begin{rem}
Note that the construction of the free diagonal action in Proposition~\ref{068} relies on the base actions $\rho^i$ acting by genuine monoid automorphisms on each $M_i$. Indeed, the universal property of the free product $M_1*\cdots*M_k$~\cite[Proposition~8.2.1]{How95} requires each map $\rho_\sigma^i:M_i\hookrightarrow M_1*\cdots*M_k$ to be a monoid homomorphism on all pairs $x,y\in M_i$, so neither actions that are only multiplicative with respect to $\tau^i$ nor general quasi-multiplicative actions extend canonically to free products.
\end{rem}

Proposition~\ref{068} ensures the existence of the \emph{twisted algebra of free diagonal invariants} $\funcs{\field}{M_1*\cdots*M_k}^G$. Furthermore, if $T^* \subseteq M_1*\cdots*M_k$ is a system of representatives for the $\omega^*$-admissible $G$-orbits satisfying Condition~\eqref{051}, Theorem~\ref{050} guarantees that this algebra admits the monomial basis $\{m_t^*\mid t\in T^*\}$.

\subsection{Hopf structure}\label{191}

In this subsection, we endow the twisted algebra of invariants $\funcs{\field}{M}^G$ with a Hopf structure. To this end, we first introduce an algebra homomorphism $\nabla$ that induces the required coproduct.

\begin{pro}\label{126}
Let $\psf:M\to M\times M$ be a twisting-preserving bigraded monoid homomorphism with finite fibres, and let $\hsf:G\times G\to G$ be an injective group homomorphism, such that\begin{equation}\label{125}\psf(\rho_{\hsf(\sigma,\varsigma)}(u))=\rho^\otimes_{(\sigma,\varsigma)}(\psf(u))\quad\text{and}\quad\omega(\hsf(\sigma,\varsigma),u)=\omega^\otimes((\sigma,\varsigma),\psf(u))\quad\text{for all}\quad u\in M\,\,\,\text{ and }\,\,\,\sigma,\varsigma\in G.\end{equation}Then, there exists a bigraded algebra homomorphism $\nabla:\funcs{\field}{M}^G\to\funcs{\field}{M\times M}^{G\times G}$ defined as follows:\begin{equation}\label{128}\nabla(g)(w)\,\,=\sum_{u\in\psf^{-1}(\{w\})}g(u)\,\,\quad\text{for all}\quad g\in\funcs{\field}{M}^G\quad\text{and}\quad w\in M\times M.\end{equation}
\end{pro}
\begin{proof}
Let $\sigma,\varsigma\in G$ and $w\in M\times M$. By the first condition in~\eqref{125}, the action of $\hsf(\sigma,\varsigma)\in G$ restricts to a bijection between the fibres $\psf^{-1}(\{w\})$ and $\psf^{-1}(\{\rho^\otimes_{(\sigma,\varsigma)}(w)\})$, indeed\[\rho_{\hsf(\sigma,\varsigma)}(\psf^{-1}(\{w\}))=\{\rho_{\hsf(\sigma,\varsigma)}(u)\mid\psf(u)=w\}\stackrel{\eqref{125}}{=}\{\rho_{\hsf(\sigma,\varsigma)}(u)\mid\psf(\rho_{\hsf(\sigma,\varsigma)}(u))=\rho^\otimes_{(\sigma,\varsigma)}(w)\}=\psf^{-1}(\{\rho^\otimes_{(\sigma,\varsigma)}(w)\}).\]Then, for every $g\in\funcs{\field}{M}^G$, we have\[\nabla(g)(\rho^\otimes_{(\sigma,\varsigma)}(w))\,\,\,=\sum_{v\in\psf^{-1}(\{\rho^\otimes_{(\sigma,\varsigma)}(w)\})}\!\!\!g(v)\,\,\,=\sum_{v\in\rho_{\hsf(\sigma,\varsigma)}(\psf^{-1}(\{w\}))}\!\!\!g(v)\,\,\,=\sum_{u\in\psf^{-1}(\{w\})}\!\!\!g(\rho_{\hsf(\sigma,\varsigma)}(u)).\]Since $g$ is $G$-invariant, by Remark~\ref{044}, the second condition in~\eqref{125} and the fact that $\psf(u)=w$, it follows that\[\nabla(g)(\rho^\otimes_{(\sigma,\varsigma)}(w))\,\,\stackrel{\eqref{043}}{=}\!\!\!\sum_{u\in\psf^{-1}(\{w\})}\!\!\omega(\hsf(\sigma,\varsigma),u)\,g(u)\,\,\stackrel{\eqref{125}}{=}\!\!\!\sum_{u\in\psf^{-1}(\{w\})}\!\!\omega^\otimes((\sigma,\varsigma),\psf(u))\,g(u)\,\,=\,\,\omega^\otimes((\sigma,\varsigma),w)\,\nabla(g)(w).\]This proves that $\nabla(g)$ is $(G\times G)$-invariant. Now, for $g,h\in\funcs{\field}{M}^G$, we have\[(\nabla(g)\cdot\nabla(h))(w)=\sum_{uv=w}\tau_\phi(u,v)\,\nabla(g)(u)\,\nabla(h)(v)=\sum_{uv=w}\tau_\phi(u,v)\,\Bigg(\sum_{a\in\psf^{-1}(\{u\})}g(a)\Bigg)\Bigg(\sum_{b\in\psf^{-1}(\{v\})}h(b)\Bigg).\]Because $\psf$ is a monoid homomorphism, for any pair $a,b$ such that $\psf(a)=u$ and $\psf(b)=v$, setting $c=ab$ yields $\psf(c)=\psf(a)\psf(b)=uv=w$. Then, as $\psf$ is twisting-preserving, we obtain\[(\nabla(g)\cdot\nabla(h))(w)\,\,=\sum_{c\in\psf^{-1}(\{w\})}\sum_{ab=c}\tau_\phi(\psf(a),\psf(b))\,g(a)\,h(b)\,\,\stackrel{\eqref{125}}{=}\sum_{c\in\psf^{-1}(\{w\})}(g\cdot h)(c)=\nabla(g\cdot h)(w).\]Furthermore, by definition and the fact that $\psf(1_{_M})=1_{_{M\times M}}$, we have\[\nabla(1)(w)=\sum_{u\in\psf^{-1}(\{w\})}\delta_{1_{_M}}(u)=\delta_{1_{_{M\times M}}}(w).\]Assume that $g\in\funcs{\field}{M_x^d}^G$, that is, $g(u)=0$ for all $u\in M$ satisfying $\bdg(u)\neq(x,d)$. Since $\psf$ is bigraded, for all $w\in M\times M$ and $u\in\psf^{-1}(\{w\})$, we have $|u|=|\psf(u)|=|w|$ and $\deg(u)=\deg(\psf(u))=\deg(w)$. Hence, if $\bdg(w)\neq(x,d)$ then $\bdg(u)\neq(x,d)$ as well, which implies $g(u)=0$ for all $u\in\psf^{-1}(\{w\})$. Thus,\[\nabla(g)(w)=\sum_{u\in\psf^{-1}(\{w\})}g(u)=0.\]This shows that $\nabla(g)$ is homogeneous of bidegree $(x,d)$. Therefore, $\nabla$ is a bigraded algebra homomorphism, concluding the proof.
\end{proof}

Before establishing the Hopf structure, we explicitly describe the action of $\nabla$ on the monomial basis.
\begin{lem}\label{129}
For each $u\in T$, we have\[\nabla(m_u)=\sum_{(v,w)\in T\times T}d_{v,w}^{\,u}\,m^\otimes_{(v,w)}\,,\quad\text{where}\qquad d_{v,w}^{\,u}=\sum_{a\in\psf^{-1}(\{(v,w)\})\cap Gu}m_u(a).\]More generally, for any elements $b=\rho_\sigma(v)$ and $c=\rho_\varsigma(w)$ in the orbits $Gv$ and $Gw$, respectively, we have\[d_{v,w}^{\,u}=\omega(\sigma,v)^{-1}\omega(\varsigma,w)^{-1}\!\!\sum_{a\in\psf^{-1}(\{(b,c)\})\cap Gu}m_u(a).\]
\end{lem}
\begin{proof}
By Proposition~\ref{126}, $\nabla(m_u)$ is a $(G\times G)$-invariant function, so it can be uniquely expanded as a linear combination of the monomial basis:\[\nabla(m_u)=\sum_{(v,w)\in T\times T}d_{v,w}^{\,u}\,m^\otimes_{(v,w)}\quad\text{for some}\quad d_{v,w}^{\,u}\in\field.\]By construction, the monomial function $m^\otimes_{(v,w)}$ is supported on the orbit of $(v,w)$ and $m^\otimes_{(v,w)}(v,w)=1$. Since distinct elements in $T\times T$ belong to distinct orbits, evaluating both sides at $(v,w)$ yields\begin{equation}\label{131}d_{v,w}^{\,u}\,=\,\nabla(m_u)(v,w)\,\stackrel{\eqref{128}}{=}\sum_{a\in\psf^{-1}(\{(v,w)\})}m_u(a)\,\,=\sum_{a\in\psf^{-1}(\{(v,w)\})\cap Gu}m_u(a),\end{equation}where the last equality is obtained because $\supp(m_u)=Gu$. 

Now, consider the elements $b=\rho_\sigma(v)$ and $c=\rho_\varsigma(w)$ for some $\sigma,\varsigma\in G$. Since $\nabla(m_u)$ is $(G\times G)$-invariant, evaluating at the pair $(b,c)=\rho^\otimes_{(\sigma,\varsigma)}(v,w)$ gives\[\nabla(m_u)(b,c)=\nabla(m_u)(\rho^\otimes_{(\sigma,\varsigma)}(v,w))=\omega^\otimes((\sigma,\varsigma),(v,w))\nabla(m_u)(v,w)\stackrel{\eqref{131}}{=}\omega(\sigma,v)\,\omega(\varsigma,w)\,d_{v,w}^{\,u}.\]Hence,\[d_{v,w}^{\,u}=\omega(\sigma,v)^{-1}\omega(\varsigma,w)^{-1}\nabla(m_u)(b,c)\stackrel{\eqref{128}}{=}\omega(\sigma,v)^{-1}\omega(\varsigma,w)^{-1}\sum_{a\in\psf^{-1}(\{(b,c)\})\cap Gu}m_u(a),\]concluding the proof.
\end{proof}

A monoid homomorphism $\psf:M\to M\times M$ is said to be \emph{coassociative} if for all $u\in M\times M\times M$ and every representative $t\in T$, it satisfies\begin{equation}\label{041}\sum_{v\in((\psf\times\id)\circ\psf)^{-1}(\{u\})\cap Gt}\!\!\omega(\sigma_v,t)\,\,\,=\!\!\sum_{v\in((\id\times\psf)\circ\psf)^{-1}(\{u\})\cap Gt}\!\!\omega(\sigma_v,t),\quad\text{where}\quad \rho_{\sigma_v}(t)=v.\end{equation}Furthermore, we say that $\psf$ is \emph{counital} if for all $t,t'\in T$, it satisfies\begin{equation}\label{089}\sum_{v\in\psf^{-1}(\{(1,t')\})\cap Gt}\!\!\omega(\sigma_v,t)\,\,\,=\,\,\,\delta_{t,t'}\,\,\,=\!\!\sum_{v\in\psf^{-1}(\{(t',1)\})\cap Gt}\!\!\omega(\sigma_v,t),\quad\text{where}\quad\rho_{\sigma_v}(t)=v.\end{equation}

A \emph{coproduct structure} on $M$ is a pair $(\psf,\hsf)$ satisfying the conditions of Proposition~\ref{126}, such that $\psf$ is coassociative and counital.

\begin{lem}\label{104}
Let $\psf:M\to M\times M$ be a monoid homomorphism and let $(\psf,\hsf)$ be a pair satisfying the conditions of Proposition~\ref{126}. Then:
\begin{enumerate}
\item\label{116}$\psf$ is coassociative if and only if for every $G$-invariant function $g$ and all $u\in M\times M\times M$, it satisfies\begin{equation}\label{108}\sum_{a\in((\psf\times\id)\circ\psf)^{-1}(\{u\})}g(a)=\sum_{a\in((\id\times\psf)\circ\psf)^{-1}(\{u\})}g(a).\end{equation}
\item\label{117}$\psf$ is counital if and only if for every $G$-invariant function $g$ and all $u\in M$, it satisfies\begin{equation}\label{122}\sum_{v\in\psf^{-1}(\{(1,u)\})}\!\!g(v)\,\,\,=\,\,\,g(u)\,\,\,=\!\!\sum_{v\in\psf^{-1}(\{(u,1)\})}\!\!g(v).\end{equation}
\end{enumerate}
\end{lem}
\begin{proof}
By Theorem~\ref{050}, the identities hold for all $G$-invariant $g$ if and only if they hold for every basis element $m_t$, where $m_t(v)$ is non-zero if and only if $v\in Gt$, in which case $m_t(v)=\omega(\sigma_v,t)$ with $\rho_{\sigma_v}(t)=v$. Evaluating~\eqref{108} at $g=m_t$ exactly yields~\eqref{041}. Furthermore, evaluating~\eqref{122} at $g=m_t$ and restricting to $u=t'\in T$ yields~\eqref{089}. Conversely, since $\nabla(m_t)$ is $(G\times G)$-invariant and $\omega(\sigma,1)=1$, the maps $u\mapsto\nabla(m_t)(1,u)$ and $u\mapsto\nabla(m_t)(u,1)$ are $G$-invariant functions on $M$; hence, if they agree with $m_t$ on $T$, they agree with $m_t$ on all $\omega$-admissible orbits by Proposition~\ref{048} and vanish on all non-admissible orbits by Proposition~\ref{045}, establishing~\eqref{122} for all $u\in M$.
\end{proof}

In what follows, we require the following additional condition that, in particular, ensures the existence of $\nabla$.
\stepcounter{condition}
\begin{enumerate}[start=\thecondition,label={(\condition)},ref={\condition}]
\item\label{127}\conditions{5}
\end{enumerate}

\begin{thm}\label{130}
The algebra homomorphism $\Delta=\Phi_2\circ\nabla$ defines a coproduct on $\funcs{\field}{M}^G$, which, together with the counit $\varepsilon:\funcs{\field}{M}^G\to\field$ given by $\varepsilon(g)=g(1)$, endows $\funcs{\field}{M}^G$ with the structure of a Hopf algebra.
\end{thm}
\begin{proof}
By Lemma~\ref{129}, for every $u\in T$, we have\begin{equation}\label{135}\Delta(m_u)=\Phi_2\Bigg(\sum_{(v,w)\in T\times T}d_{v,w}^{\,u}m^\otimes_{(v,w)}\Bigg)=\sum_{(v,w)\in T\times T}d_{v,w}^{\,u}\,m_v\otimes m_w.\end{equation}Then, applying $\Delta\otimes\id$ yields\[\begin{array}{rclcl}
(\Delta\otimes\id)(\Delta(m_u))
&=&\displaystyle\sum_{(v,w)\in T\times T}d_{v,w}^{\,u}\,\Delta(m_v)\otimes m_w
&=&\displaystyle\sum_{(v,w)\in T\times T}d_{v,w}^{\,u}\Bigg(\sum_{(a,b)\in T\times T}d_{a,b}^{\,v}\,m_a\otimes m_b\Bigg)\otimes m_w\\
%&&&=&\displaystyle\sum_{(v,w)\in T\times T}\sum_{(a,b)\in T\times T}d_{v,w}^{\,u}\,d_{a,b}^{\,v}\,m_a\otimes m_b\otimes m_w\\
%&&&=&\displaystyle\sum_{(v,w,a,b)\in T\times T\times T\times T}d_{v,w}^{\,u}\,d_{a,b}^{\,v}\,m_a\otimes m_b\otimes m_w\\
%&&&=&\displaystyle\sum_{(v,c,a,b)\in T\times T\times T\times T}d_{v,c}^{\,u}\,d_{a,b}^{\,v}\,m_a\otimes m_b\otimes m_c\\[1.5cm]
&&&=&\displaystyle\sum_{a,b,c\in T}\Bigg(\sum_{v\in T}d_{v,c}^{\,u}\,d_{a,b}^{\,v}\Bigg)m_a\otimes m_b\otimes m_c.
\end{array}\]Similarly, applying $\id\otimes\Delta$ gives
\[(\id\otimes\Delta)(\Delta(m_u))=\sum_{(v,w)\in T\times T}d_{v,w}^{\,u}\,m_v\otimes\Delta(m_w)\,\,=\sum_{a,b,c\in T}\Bigg(\sum_{w\in T}d_{a,w}^{\,u}d_{b,c}^{\,w}\Bigg)m_a\otimes m_b\otimes m_c.\]
Since the simple tensors $m_a\otimes m_b\otimes m_c$ with $a,b,c\in T$ form a basis for $\funcs{\field}{M}^G\otimes\funcs{\field}{M}^G\otimes\funcs{\field}{M}^G$, coassociativity holds if and only if the structure constants satisfy\begin{equation}\label{134}\sum_{v\in T}d_{v,c}^{\,u}\,d_{a,b}^{\,v}\,\,=\sum_{w\in T}d_{a,w}^{\,u}\,d_{b,c}^{\,w}\quad\text{for all}\quad a,b,c\in T.\end{equation}By Lemma~\ref{129}, we expand the left side as follows\begin{equation}\label{132}\sum_{v\in T}d_{v,c}^{\,u}d_{a,b}^{\,v}=\sum_{v\in T}d_{v,c}^{\,u}\Bigg(\sum_{d\in\psf^{-1}(\{(a,b)\})\cap Gv}m_v(d)\Bigg)=\sum_{v\in T}\sum_{d\in\psf^{-1}(\{(a,b)\})\cap Gv}d_{v,c}^{\,u}\,m_v(d).\end{equation}For each $d\in Gv$, there exists $\sigma\in G$ such that $d=\rho_\sigma(v)$, which implies $m_v(d)=\omega(\sigma,v)$. By applying Lemma~\ref{129}, we obtain\[d_{v,c}^{\,u}\,=\,\omega(\sigma,v)^{-1}\omega(1,c)^{-1}\sum_{x\in\psf^{-1}(\{(d,c)\})\cap Gu}m_u(x)\,\stackrel{\eqref{003}}{=}\,\omega(\sigma,v)^{-1}\sum_{x\in\psf^{-1}(\{(d,c)\})\cap Gu}m_u(x).\]Hence, multiplying by $m_v(d)$ cancels the weight $\omega(\sigma,v)$, yielding\begin{equation}\label{133}d_{v,c}^{\,u}\,m_v(d)\,=\sum_{x\in\psf^{-1}(\{(d,c)\})\cap Gu}m_u(x).\end{equation}Note that the orbits $Gv$ with $v\in T$ partition the $\omega$-admissible elements of $M$. If $d\in\psf^{-1}(\{(a,b)\})$ is not $\omega$-admissible, then $(d,c)$ is not $\omega^\otimes$-admissible by Proposition~\ref{092}\eqref{095}; since $\nabla(m_u)$ is $(G\times G)$-invariant, Proposition~\ref{045} implies that\[\sum_{x\in\psf^{-1}(\{(d,c)\})\cap Gu}m_u(x)=\nabla(m_u)(d,c)=0.\]Thus, these vanishing terms can be included in the sum over all $d\in\psf^{-1}(\{(a,b)\})$, and substituting~\eqref{133} back into~\eqref{132} gives\[\sum_{v\in T}d_{v,c}^{\,u}d_{a,b}^{\,v}=\sum_{v\in T}\sum_{d\in\psf^{-1}(\{(a,b)\})\cap Gv}\Bigg(\sum_{x\in\psf^{-1}(\{(d,c)\})\cap Gu}m_u(x)\Bigg)=\sum_{d\in\psf^{-1}(\{(a,b)\})}\,\sum_{x\in\psf^{-1}(\{(d,c)\})\cap Gu}m_u(x).\]Notice that\[\{x\in Gu\mid\psf(x)=(d,c)\text{ with } \psf(d)=(a,b)\text{ for some }d\in M\}\,\,=\,\,\{x\in Gu\mid((\psf\times\id)\circ\psf)(x)=(a,b,c)\}.\]Hence,\[\sum_{v\in T}d_{v,c}^{\,u}d_{a,b}^{\,v}=\sum_{x\in((\psf\times\id)\circ\psf)^{-1}(\{(a,b,c)\})\cap Gu}m_u(x).\]By a symmetric argument, expanding the right side in~\eqref{134} yields\[\sum_{w\in T}d_{a,w}^{\,u}d_{b,c}^{\,w}=\sum_{x\in((\id\times\psf)\circ\psf)^{-1}(\{(a,b,c)\})\cap Gu}m_u(x).\]Since $\psf$ is coassociative, by~\eqref{108} we obtain that~\eqref{134} holds, and therefore $\Delta$ is a coproduct.

Recall that $\rho_\sigma(1)=1$ and $\omega(\sigma,1)=1$ for all $\sigma\in G$. Then, $1\in T$, and $\varepsilon(m_v)=m_v(1)=\delta_{\{1\}}(v)$. Applying Lemma~\ref{129} gives\[\begin{array}{rclclcl}(\varepsilon\otimes\id)(\Delta(m_u))
&=&\displaystyle\sum_{(v,w)\in T\times T}d_{v,w}^{\,u}\,\varepsilon(m_v)\,m_w
&=&\displaystyle\sum_{w\in T}d_{1,w}^{\,u}\,m_w
&=&\displaystyle\sum_{w\in T}\Bigg(\sum_{a\in\psf^{-1}(\{(1,w)\})\cap Gu}\!\!\!m_u(a)\Bigg)m_w\\[0.7cm]
&&&&&\stackrel{\eqref{122}}{=}&\displaystyle\sum_{w\in T}m_u(w)\,m_w\,\,\,=\,\,\,m_u.\end{array}\]Similarly, we obtain $(\id\otimes\varepsilon)(\Delta(m_u))=m_u$. Thus, $\funcs{\field}{M}^G$ is a bigraded bialgebra. Finally, since $\funcs{\field}{M}^G$ is graded, Proposition~\ref{105} implies that $\funcs{\field}{M}^G$ is a Hopf algebra.
\end{proof}

By~\eqref{135} and Lemma~\ref{129}, for every $u\in T$, the coproduct on the monomial basis is explicitly given by:\begin{equation}\label{139}\Delta(m_u)=\sum_{(v,w)\in T\times T}\Bigg(\sum_{a\in\psf^{-1}(\{(v,w)\})\cap Gu}m_u(a)\Bigg)m_v\otimes m_w.\end{equation}

We now show how to transfer the coproduct structure $(\psf,\hsf)$ to the restricted actions defined in Subsection~\ref{171}.

\begin{pro}\label{170}
Let $S\subset G$ be a generating set, and let $\{M_s\mid s\in S\}$ be a family of subsets defining a restricted action $\bar{\rho}$ and a restricted weight $\bar{\omega}$ on $M$, together with a system of representatives $\bar{T}\subset M$ for the $\bar{\omega}$-admissible $\bar{G}$-orbits satisfying Condition~\eqref{051}. Furthermore, suppose that the following conditions hold:
\begin{enumerate}
\item\label{172}For every $s\in S$, we have $\hsf(s,1)\in S$ and $\hsf(1,s)\in S$.
\item\label{173}For every $s\in S$ and every $u\in M$ with $\psf(u)=(u_1,u_2)$, we have $u\in M_{\hsf(s,1)}$ if and only if $u_1\in M_s$, and $u\in M_{\hsf(1,s)}$ if and only if $u_2\in M_s$.
\item\label{174}For every $u\in M\times M\times M$, there exists a bijection $f_u:((\psf\times\id)\circ\psf)^{-1}(\{u\})\to((\id\times\psf)\circ\psf)^{-1}(\{u\})$ such that for all $v\in((\psf\times\id)\circ\psf)^{-1}(\{u\})$, we have $f_u(v)=\bar{\rho}_\sigma(v)$ for some $\sigma\in\bar{G}$ with $\bar{\omega}(\sigma,v)=1$. Furthermore, for every $t\in\bar{T}$, the fibres $\psf^{-1}(\{(1,t)\})$ and $\psf^{-1}(\{(t,1)\})$ are
contained in $\bar{G}t$ and satisfy\[\sum_{v\in\psf^{-1}(\{(1,t)\})}\bar{\omega}(\sigma_v,t)=1=\sum_{v\in\psf^{-1}(\{(t,1)\})}\bar{\omega}(\sigma_v,t),\quad\text{where}\quad \bar{\rho}_{\sigma_v}(t)=v.\]
\end{enumerate}
Then $(\psf,\hsf)$ is a coproduct structure on $M$ with respect to $(\bar{G},\bar{\rho},\bar{\omega})$. Consequently, $\funcs{\field}{M}^{\bar{G}}$ is a connected Hopf algebra, and the canonical inclusion $\funcs{\field}{M}^G\hookrightarrow\funcs{\field}{M}^{\bar{G}}$ is an injective homomorphism of Hopf algebras.
\end{pro}
\begin{proof}
We verify that $(\psf,\hsf)$ satisfies the identities~\eqref{125} of Proposition~\ref{126} and the coassociativity and counitality identities~\eqref{041} and~\eqref{089} with respect to $(\bar{\rho},\bar{\omega})$.

Since $\hsf$ induces a group homomorphism on $\bar{G}\times\bar{G}$ and the weight $\bar{\omega}$ satisfies the $1$-cocycle condition, it suffices to verify the identities on the generators $(s,1)$ and $(1,s)$ of $\bar{G}\times\bar{G}$. Let $u\in M$ with $\psf(u)=(u_1,u_2)$. By~\eqref{172}, we have $\hsf(s,1)\in S$. If $u\in M_{\hsf(s,1)}$, then~\eqref{173} implies $u_1\in M_s$. Thus $\bar{\rho}_{\hsf(s,1)}(u)=\rho_{\hsf(s,1)}(u)$ and $\bar{\rho}_s(u_1)=\rho_s(u_1)$, which yields\begin{gather*}
\psf(\bar{\rho}_{\hsf(s,1)}(u))=\psf(\rho_{\hsf(s,1)}(u))\overset{\eqref{125}}{=} (\rho_s(u_1),u_2)=(\bar{\rho}_s(u_1),u_2)=\bar{\rho}^\otimes_{(s,1)}(\psf(u)),\\
\bar{\omega}(\hsf(s,1),u)=\omega(\hsf(s,1),u)\overset{\eqref{125}}{=}\omega(s,u_1)\,\omega(1,u_2)=\bar{\omega}(s,u_1)\,\bar{\omega}(1,u_2)=\bar{\omega}^\otimes((s,1),\psf(u)).\end{gather*}If $u\notin M_{\hsf(s,1)}$, then~\eqref{173} implies $u_1\notin M_s$. By definition of the restricted action and weight, we obtain trivially\[\psf(\bar{\rho}_{\hsf(s,1)}(u))=\psf(u)=\bar{\rho}^\otimes_{(s,1)}(\psf(u))\quad\text{and}\quad\bar{\omega}(\hsf(s,1),u)=1=\bar{\omega}^\otimes((s,1),\psf(u)).\]The case $(1,s)$ is symmetric.

Let $u\in M\times M\times M$ and $t\in\bar{T}$, and set $F_1=((\psf\times\id)\circ\psf)^{-1}(\{u\})$ and $F_2=((\id\times\psf)\circ\psf)^{-1}(\{u\})$. By~\eqref{174}, there exists a bijection $f_u:F_1\to F_2$ that preserves each $\bar{G}$-orbit $\bar{G}t$, so it restricts to a bijection $f_u:F_1\cap\bar{G}t\to F_2\cap\bar{G}t$. Moreover, for each $v\in F_1\cap\bar{G}t$ with $v=\bar{\rho}_{\sigma_v}(t)$, writing $w=f_u(v)=\bar{\rho}_\sigma(v)=\bar{\rho}_{\sigma\sigma_v}(t)$ with $\bar{\omega}(\sigma,v)=1$, the $1$-cocycle condition yields $\bar{\omega}(\sigma_w,t)=\bar{\omega}(\sigma\sigma_v,t)=\bar{\omega}(\sigma,v)\,\bar{\omega}(\sigma_v,t)=\bar{\omega}(\sigma_v,t)$. Consequently,\[\sum_{v\in F_1\cap\bar{G}t}\bar{\omega}(\sigma_v,t)=\sum_{w\in F_2\cap\bar{G}t}\bar{\omega}(\sigma_w,t),\]which establishes the coassociativity in~\eqref{041}.

Let $t,t'\in\bar{T}$. By~\eqref{174}, the fibre $\psf^{-1}(\{(1,t')\})$ is entirely contained in the restricted orbit $\bar{G}t'$. If $t\neq t'$, the orbits $\bar{G}t$ and $\bar{G}t'$ are disjoint, so $\psf^{-1}(\{(1,t')\})\cap\bar{G}t=\emptyset$, yielding a sum of $0=\delta_{t,t'}$. If $t=t'$, we have $\psf^{-1}(\{(1,t)\})\cap\bar{G}t=\psf^{-1}(\{(1,t)\})$, and the sum equals $1=\delta_{t,t}$ directly by~\eqref{174}. The second equality in~\eqref{089} follows symmetrically. Hence, $(\psf,\hsf)$ is a coproduct structure on $M$ with respect to $(\bar{G},\bar{\rho},\bar{\omega})$, so Theorem~\ref{130} implies that $\funcs{\field}{M}^{\bar{G}}$ is a connected Hopf algebra. Finally, since each unrestricted orbit decomposes into a disjoint union of restricted orbits and both coproducts are induced by the same monoid homomorphism $\psf$, the canonical inclusion $\funcs{\field}{M}^G\hookrightarrow\funcs{\field}{M}^{\bar{G}}$ is an injective homomorphism of Hopf algebras.
\end{proof}

\subsubsection{Duplication}\label{035}

Here, we describe a method to construct a monoid homomorphism $\psf:M\to M\times M$ that can potentially serve as the map needed in Condition~\eqref{127}. Before proceeding, we establish a preparatory proposition.

\begin{pro}\label{036}
Let $\psi:M\to N$ be a monoid isomorphism. Then:
\begin{enumerate}
\item\label{107}The map $\tau^\psi:N\times N\to\field$ defined by $\tau^\psi(u,v)=\tau(\psi^{-1}(u),\psi^{-1}(v))$ constitutes a twisting on the monoid $N$.
\item\label{039}There exists a unique quasi-multiplicative group action $\rho^\psi:G\to\sym(N)$ such that the following diagram commutes for each $\sigma\in G$:\begin{equation}\label{042}\vcdraw{\xymatrix{M\ar[d]_\psi\ar[r]^{\rho_\sigma}&M\ar[d]^\psi\\N\ar[r]_{\rho_\sigma^\psi}&N}}\end{equation}Consequently, every group action $\varrho:G\to\sym(N)$ such that $\psi$ intertwines the actions, that is, $\psi(\rho_\sigma(u))=\varrho_\sigma(\psi(u))$, must coincide with $\rho^\psi$ by uniqueness.
\item\label{040}The map $\omega^\psi:G\times N\to\field^\times$ defined by $\omega^\psi(\sigma,u):=\omega(\sigma,\psi^{-1}(u))$ is a weight function for $\rho^\psi$, such that $(\tau^\psi,\omega^\psi)$ is a multiplicative pair for $\rho^\psi$.
\item\label{037}If $\psi$ is bigraded, there exists a bigraded isomorphism $\funcs{\field}{M}^G\to\funcs{\field}{N}^{G^\psi}$, where $G^\psi$ denotes the group $G$ equipped with the induced action $\rho^\psi$ and the weight function $\omega^\psi$.
\end{enumerate}
\end{pro}
\begin{proof}
\eqref{107}~It follows immediately from the conditions in \eqref{001}. \eqref{039}~Existence and uniqueness follow directly by defining $\rho_\sigma^\psi=\psi\circ\rho_\sigma\circ\psi^{-1}$, which makes the diagram commute by construction and is a well-defined group action since $\psi$ is a bijection. Furthermore, given $\sigma\in G$ and $w\in N$, since $\rho$ is quasi-multiplicative on $\psi^{-1}(w)$, there exist $\alpha,\varsigma\in G$ such that for all $u,v\in N$ with $uv=w$ and $\tau^\psi(u,v)=\tau(\psi^{-1}(u),\psi^{-1}(v))\neq0$, we have\[\rho_\sigma^\psi(uv)=\psi(\rho_\sigma(\psi^{-1}(uv)))=\psi(\rho_\alpha(\psi^{-1}(u))\,\rho_\varsigma(\psi^{-1}(v)))=\psi(\rho_\alpha(\psi^{-1}(u)))\,\psi(\rho_\varsigma(\psi^{-1}(v)))=\rho_\alpha^\psi(u)\,\rho_\varsigma^\psi(v),\]proving that $\rho^\psi$ is quasi-multiplicative with respect to $\tau^\psi$.

\eqref{040}~For any $\sigma,\varsigma\in G$ and $u\in N$, we have\[\begin{array}{rclllll}
\omega^\psi(\sigma\varsigma,u)&=&\omega(\sigma\varsigma,\psi^{-1}(u))&\stackrel{\eqref{003}}{=}&\omega(\sigma,\rho_\varsigma(\psi^{-1}(u)))\,\omega(\varsigma,\psi^{-1}(u))\\
&&&\stackrel{\eqref{042}}{=}&\omega(\sigma,\psi^{-1}(\rho_\varsigma^\psi(u)))\,\omega(\varsigma,\psi^{-1}(u))&=&\omega^\psi(\sigma,\rho_\varsigma^\psi(u))\,\omega^\psi(\varsigma,u).
\end{array}\]Moreover, given $\sigma\in G$ and $w\in N$, let $\alpha,\varsigma\in G$ be the elements associated with $(\sigma,\psi^{-1}(w))$. Since $(\tau,\omega)$ is a multiplicative pair, for all $u,v\in N$ with $uv=w$ and $\tau^\psi(u,v)\neq0$, we have\[\begin{array}{rcl}
\tau^\psi(u,v)\,\omega^\psi(\sigma,uv)&=&\tau(\psi^{-1}(u),\psi^{-1}(v))\,\omega(\sigma,\psi^{-1}(uv))\\
&\stackrel{\eqref{028}}{=}&\tau(\rho_\alpha(\psi^{-1}(u)),\rho_\varsigma(\psi^{-1}(v)))\,\omega(\alpha,\psi^{-1}(u))\,\omega(\varsigma,\psi^{-1}(v))\\
&\stackrel{\eqref{042}}{=}&\tau(\psi^{-1}(\rho_\alpha^\psi(u)),\psi^{-1}(\rho_\varsigma^\psi(v)))\,\omega^\psi(\alpha,u)\,\omega^\psi(\varsigma,v)\\
&=&\tau^\psi(\rho_\alpha^\psi(u),\rho_\varsigma^\psi(v))\,\omega^\psi(\alpha,u)\,\omega^\psi(\varsigma,v).
\end{array}\]Therefore, $(\tau^\psi,\omega^\psi)$ is a multiplicative pair for $\rho^\psi$.

\eqref{037}~By Proposition~\ref{019}\eqref{020}, the map $\psi$ induces a bigraded algebra isomorphism $\Psi:\funcs{\field}{M}\to\funcs{\field}{N}$ defined by $\Psi(g)=g\circ\psi^{-1}$. Since $\Psi$ is already bigraded and its inverse is given by $\Psi^{-1}(g)=g\circ\psi$, it is enough to show that $\Psi$ maps $\funcs{\field}{M}^G$ into $\funcs{\field}{N}^{G^\psi}$; indeed, for any $g\in\funcs{\field}{M}^G$ and $\sigma\in G$, we have\[\Psi(g)(\rho_\sigma^\psi(u))=g(\psi^{-1}(\rho_\sigma^\psi(u)))\stackrel{\eqref{042}}{=}g(\rho_\sigma(\psi^{-1}(u)))\stackrel{\eqref{043}}{=}\omega(\sigma,\psi^{-1}(u))g(\psi^{-1}(u))=\omega^\psi(\sigma,u)\Psi(g)(u).\]By Remark~\ref{044}, it follows that $\Psi(g)\in\funcs{\field}{N}^{G^\psi}$; consequently, \mbox{the restriction $\Psi|_{\funcs{\field}{M}^G}:\funcs{\field}{M}^G\to\funcs{\field}{N}^{G^\psi}$} is a bigraded isomorphism.
\end{proof}

To achieve our goal, we implement a duplication of the monoid $M$ using a congruence $R\subseteq[M,M]$ on the free product $M*M$. We assume there exists a monoid isomorphism $\psi:M\to(M*M)/R$ ensuring that $(M*M)/R$ satisfies Conditions~\eqref{015}~to~\eqref{051}. By Proposition~\ref{036}\eqref{107}, the monoid $(M*M)/R$ can be equipped with the twisting\begin{equation}\label{110}\tau^\psi(u,v)=\tau(\psi^{-1}(u),\psi^{-1}(v))\quad\text{for all}\quad u,v\in(M*M)/R.\end{equation}Proposition~\ref{036} ensures the existence of a quasi-multiplicative action $\rho^\psi:G^\psi\to\sym((M*M)/R)$ together with a weight function $\omega^\psi:G^\psi\times(M*M)/R\to\field^\times$ such that $(\tau^\psi,\omega^\psi)$ is a multiplicative pair; explicitly\begin{equation}\label{138}\rho^\psi_\sigma=\psi\circ\rho_\sigma\circ\psi^{-1}\quad\text{and}\quad\omega^\psi(\sigma,u)=\omega(\sigma,\psi^{-1}(u))\quad\text{for all}\quad\sigma\in G^\psi\quad\text{and}\quad u\in(M*M)/R.\end{equation}

\begin{rem}\label{109}
If $M$ is $\P$-graded with respect to $\tau$ by a degree map $\nu:M\to\P$ and $\psi:M\to N$ is a monoid isomorphism, then $N$ is also $\P$-graded with respect to $\tau^\psi$ by the degree map $\nu^\psi=\nu\circ\psi^{-1}$. Indeed, if $\tau^\psi(x,y)\neq 0$, then $\tau(\psi^{-1}(x),\psi^{-1}(y))\neq 0$, so $\nu^\psi(xy)=\nu(\psi^{-1}(x)\psi^{-1}(y))=\nu(\psi^{-1}(x))\nu(\psi^{-1}(y))=\nu^\psi(x)\nu^\psi(y)$. By definition, $N_x=\psi(M_x)$ for all $x\in\P$, meaning that $\psi$ is a $\P$-graded isomorphism. Furthermore, if $M_1=\{1\}$, then $N_1=\psi(M_1)=\psi(\{1\})=\{1\}$. Therefore, if $M$ is graded, then $N$ is graded as well.
\end{rem}

\begin{rem}\label{136}
Since $\psi$ is a monoid isomorphism, $vw=u$ with $\tau(v,w)\neq 0$ if and only if $\psi(v)\psi(w)=\psi(u)$ with $\tau^\psi(\psi(v),\psi(w))\neq 0$. Therefore, $M$ is locally finite with respect to $\tau$ if and only if $(M*M)/R$ is locally finite with respect to $\tau^\psi$. By this fact and Remark~\ref{109}, the monoid $(M*M)/R$ is locally finite with respect to $\tau^\psi$ and $(\G\times\N)$-graded with respect to $\tau^\psi$, with gradings given by $\{\psi(\ubar{M}_x)\mid x\in\G\}$, and $\{\psi(M^d)\mid d\in\N\}$ with $\psi(M^0)=\{1\}$.

Moreover, since $M^d\cap\ubar{M}_x=\psi^{-1}(\psi(M^d)\cap\psi(\ubar{M}_x))$, each $\psi(M^d)$ intersects only finitely many fibres $\psi(\ubar{M}_x)$. 

By definition, $\rho_\sigma^\psi(\psi(M^d))=\psi(M^d)$ and $\rho_\sigma^\psi(\psi(\ubar{M}_x))=\psi(\ubar{M}_x)$ for all $\sigma\in G$. Therefore, $\rho^\psi$ is bigraded.

Furthermore, it is straightforward to verify that $G^\psi\psi(u)=\psi(Gu)$ and $G^\psi_{(\psi(u))}=G_{(u)}$. This implies that $T^\psi:=\psi(T)$ defines a transversal of the $\omega^\psi$-admissible $G^\psi$-orbits, and each $\psi(M^d)$ intersects only finitely many $\omega^\psi$-admissible $G^\psi$-orbits.
%PROOF: (1) Given $v\in G^\psi\psi(u)$, there is $\sigma\in G$ such that $v=\rho^\psi_\sigma(\psi(u))=\psi(\rho_\sigma(u))\in\psi(Gu)$. The reverse inclusion follows symmetrically because $\psi$ is a bijection, which yields $G^\psi\psi(u)=\psi(Gu)$. Furthermore, by definition, $\rho^\psi_\sigma(\psi(u))=\psi(u)$ if and only if $\psi(\rho_\sigma(u))=\psi(u)$. Therefore, $G^\psi_{(\psi(u))}=G_{(u)}$. (2) Part~(1) guarantees that $\psi$ maps the $G$-orbits of $M$ bijectively onto the $G^\psi$-orbits of $(M*M)/R$. Consequently, $T^\psi$ contains exactly one representative for each orbit in $\{\psi(Gu) \mid u \in T\}$. For any $u\in M$ and $\sigma\in G_{(u)}$, the definition of the weight function yields $\omega^\psi(\sigma,\psi(u))=\omega(\sigma,u)$. By part~(1), we have $\omega^\psi(\sigma,\psi(u))=1$ for all $\sigma\in G^\psi_{(\psi(u))}$ if and only if $\omega(\sigma,u)=1$ for all $\sigma\in G_{(u)}$. Therefore, $T^\psi$ is a system of representatives of all the admissible $G^\psi$-orbits. (3) Let $Gu_1,\ldots,Gu_r$ be the finitely many $G$-orbits that intersect $M^d$. By part~(1), the $G^\psi$-orbits intersecting $\psi(M^d)$ are precisely $G^\psi\psi(u_1),\ldots,G^\psi\psi(u_r)$. Thus, $\psi(M^d)$ intersects only finitely many $G^\psi$-orbits.

Consequently, $(M*M)/R$ satisfies Conditions~\eqref{015}~to~\eqref{051}. Therefore, the twisted algebra of invariants $\funcs{\field}{(M*M)/R}^{G^\psi}$ is well-defined and admits the monomial basis $\{m_t^\psi\mid t\in T^\psi\}$, where $m_t^\psi$ denotes the monomial $G^\psi$-invariant function associated to $t\in T^\psi$.

Thus, by Proposition~\ref{036}\eqref{037}, there exists a bigraded algebra isomorphism $\Psi:\funcs{\field}{M}^G\to\funcs{\field}{(M*M)/R}^{G^\psi}$ defined by $\Psi(g)=g\circ\psi^{-1}$ for all $g\in\funcs{\field}{M}^G$. Furthermore, for any $t\in T$, we have $\Psi(m_t)=m^\psi_{\psi(t)}$.
%PROOF: The existence of $\Psi$ is a direct application of Proposition~\ref{036}\eqref{037}. For the second part, observe that for each $t\in T$, $\Psi(m_t)(u)=m_t(\psi^{-1}(u))=0$ if and only if $\psi^{-1}(u)\notin Gt$, which is equivalent to $u\notin G^\psi\psi(t)$. Hence, $\supp(\Psi(m_t))=G^\psi\psi(t)$. Given $u\in G^\psi\psi(t)$, there exists $\sigma\in G$ such that $u=\rho^\psi_\sigma(\psi(t))$, which is equivalent to $\psi^{-1}(u)=\rho_\sigma(t)$. Consequently, evaluating the series yields\[\Psi(m_t)(u)=m_t(\psi^{-1}(u))=\omega(\sigma,\psi^{-1}(u))=\omega^\psi(\sigma,u)=m^\psi_{\psi(t)}(u).\]Therefore $\Psi(m_t)=m^\psi_{\psi(t)}$, as required.
\end{rem}

As $R\subseteq[M,M]$, the direct product $M\times M$ is a quotient of $(M*M)/R$. Specifically, as explained in Subsubsection~\ref{002}, we have\begin{equation*}M\times M\simeq(M*M)/[M,M]\simeq(\,(M*M)/R\,)\,/\,(\,[M,M]/R\,).\end{equation*}We denote by $\qsf:(M*M)/R\to M\times M$ the canonical projection, defined by\[\qsf([u]_{_R})=[[u]_{_R}]_{_{[M,M]/R}}=[u]_{_{[M,M]}}\quad\text{for all}\quad u\in M*M.\]

\begin{rem}\label{124}
By definition, for any $u\in M*M$, the image $\qsf([u]_{_R})$ coincides with the canonical projection $M*M\to M\times M$ applied to the representative $u$. Recall that any element $u\in M*M$ can be uniquely written as an alternating product; that is, there exist $v_1,\ldots,v_k\in M$ corresponding to elements in the first copy of $M$, and $w_1,\ldots,w_k\in M$ corresponding to elements in the second copy, where at most $v_1$ and $w_k$ could be $1$, such that $u=v_1w_1\cdots v_kw_k$. Therefore, $\qsf([u]_{_R})=(v_1\cdots v_k,w_1\cdots w_k)$.
\end{rem}

With these maps established, we define the required monoid homomorphism as $\psf=\qsf\circ\psi$.

The following proposition shows that several compatibility requirements for the map $\psf$ are directly inherited from corresponding properties on the canonical projection $\qsf$.
\begin{pro}\label{106}
\begin{enumerate}
\item\label{111}If $\qsf$ has finite fibres, then $\psf$ has finite fibres.
\item\label{112}If $\qsf$ is twisting-preserving, then $\psf$ is twisting-preserving.
\item\label{113}Assume that there exists a group homomorphism $\hsf:G\times G\to G$ such that for all $u,v\in(M*M)/R$ and $\sigma,\varsigma\in G$, it satisfies\begin{equation}\label{137}\qsf(\rho^\psi_{\hsf(\sigma,\varsigma)}(u))=\rho^\otimes_{(\sigma,\varsigma)}(\qsf(u))\quad\text{and}\quad\omega^\psi(\hsf(\sigma,\varsigma),u)=\omega^\otimes((\sigma,\varsigma),\qsf(u)).\end{equation}Then, the pair $(\psf,\hsf)$ satisfies the compatibility conditions in~\eqref{125}.
\end{enumerate}
\end{pro}
\begin{proof}
\eqref{111} It is a direct consequence of the fact that $\psi$ is a monoid isomorphism. \eqref{112} Let $u,v\in M$ and $\sigma,\varsigma\in G$. Applying the definition of the twisting $\tau^\psi$ and the fact that $\qsf$ is twisting-preserving, we have\[\tau_\phi(\psf(u),\psf(v))=\tau_\phi(\qsf(\psi(u)),\qsf(\psi(v)))=\tau^\psi(\psi(u),\psi(v))\stackrel{\eqref{110}}{=}\tau(u,v).\]\eqref{113} Similarly, applying the conditions in~\eqref{137} and the definitions of the group action $\rho^\psi$ and the weight function $\omega^\psi$, we have\begin{gather*}\psf(\rho_{\hsf(\sigma,\varsigma)}(u))=\qsf(\psi(\rho_{\hsf(\sigma,\varsigma)}(u)))\stackrel{\eqref{138}}{=}\qsf(\rho^\psi_{\hsf(\sigma,\varsigma)}(\psi(u)))\stackrel{\eqref{137}}{=}\rho^\otimes_{(\sigma,\varsigma)}(\qsf(\psi(u)))=\rho^\otimes_{(\sigma,\varsigma)}(\psf(u)),\quad\text{and}\\\omega(\hsf(\sigma,\varsigma),u)\stackrel{\eqref{138}}{=}\omega^\psi(\hsf(\sigma,\varsigma),\psi(u))\stackrel{\eqref{137}}{=}\omega^\otimes((\sigma,\varsigma),\qsf(\psi(u)))=\omega^\otimes((\sigma,\varsigma),\psf(u)).\end{gather*}Therefore, $(\psf,\hsf)$ satisfies the compatibility conditions required in~\eqref{125}, concluding the proof.
\end{proof}

Consequently, if $M$ admits a duplication $(M*M)/R$ such that the projection $\qsf$ satisfies the hypotheses of Proposition~\ref{106} with $\hsf$ injective, and the resulting map $\psf=\qsf\circ\psi$ is bigraded, coassociative, and counital, then the pair $(\psf,\hsf)$ constitutes a coproduct structure on $M$. In this case, the monoid $M$ fulfills Condition~\eqref{127}, and by Theorem~\ref{130}, the twisted algebra of invariants $\funcs{\field}{M}^G$ is endowed with the structure of a Hopf algebra.

\subsection{Quotients}\label{176}

Having established the conditions under which the twisted algebra of invariants $\funcs{\field}{M}^G$ admits a Hopf structure, we now investigate how this structure behaves under monoid quotients. Specifically, our goal is to construct a quotient Hopf algebra by means of a congruence on $M$. For this purpose, we further assume that $M$ satisfies Condition~\eqref{127}, which ensures that $\funcs{\field}{M}^G$ carries a Hopf structure as established in Theorem~\ref{130}.\medskip

A congruence $S$ on $M$ is said to be \emph{stable} under $\rho$, or that $\rho$ \emph{preserves} $S$, if $(\rho_\sigma(u),\rho_\sigma(v))\in S$ whenever $(u,v)\in S$, for all $\sigma\in G$.

\begin{pro}\label{073}
Let $S\subseteq M\times M$ be a congruence on $M$.
\begin{enumerate}
\item\label{074}If $\tau(u,v)=\tau(u',v')$ for all $(u,u'),(v,v')\in S$, then the map $\tau^{_S}:M/S\times M/S\to\field$ given by $\tau^{_S}([u],[v])=\tau(u,v)$ defines a twisting on $M/S$.
\item\label{119}If $S$ is stable under $\rho$, then the map $\rho^{_S}:G\to\sym(M/S)$ defined by $\rho^{_S}_{\sigma}([u])=[\rho_{\sigma}(u)]$, induces a well-defined group action. Furthermore, if $\tau$ is compatible with $S$ as in~\eqref{074} and $\rho$ is multiplicative with respect to $\tau$, then the induced action $\rho^{_S}$ is multiplicative with respect to $\tau^{_S}$.
\item\label{145}If the action $\rho:G\to\sym(M)$ is $\P$-graded and $S$ is $\P$-graded and stable under $\rho$, then the induced action $\rho^{_S}$ on $M/S$ is also $\P$-graded.
\end{enumerate}
\end{pro}
\begin{proof}
\eqref{074}~It follows immediately from the conditions in \eqref{001}. \eqref{119}~Since $S$ is stable under $\rho$, for each $\sigma\in G$, the map $\rho^{_S}_\sigma:M/S\to M/S$ given by $\rho^{_S}_\sigma([u])=[\rho_\sigma(u)]$ is well-defined and bijective. If $\rho$ is multiplicative with respect to $\tau$, then for any $[u],[v]\in M/S$ with $\tau^{_S}([u],[v])=\tau(u,v)\neq0$, we have $\rho_\sigma^{_S}([u][v])=[\rho_\sigma(uv)]=[\rho_\sigma(u)\rho_\sigma(v)]=\rho_\sigma^{_S}([u])\rho_\sigma^{_S}([v])$, so $\rho^{_S}$ is multiplicative with respect to $\tau^{_S}$. Moreover, for any $\sigma,\varsigma\in G$, we have $\rho^{_S}_{\sigma\varsigma}([u])=[\rho_{\sigma\varsigma}(u)]=\rho^{_S}_\sigma([\rho_\varsigma(u)])=(\rho^{_S}_\sigma\circ\rho^{_S}_\varsigma)([u])$, confirming that $\rho^{_S}$ is a group action on $M/S$.

\eqref{145}~Let $\nu:M\to\P$ be the degree map of $M$. Since $S$ is $\P$-graded and stable under $\rho$, the induced degree map $\nu^{_S}$ on the quotient is well-defined and satisfies $\nu^{_S}(\rho^{_S}_\sigma([u]))=\nu^{_S}([\rho_\sigma(u)])=\nu(\rho_\sigma(u))=\nu(u)=\nu^{_S}([u])$, where the third equality holds because $\rho$ is $\P$-graded. Therefore, $\rho^{_S}$ is a $\P$-graded action.
\end{proof}

If $M$ is $\P$-graded with respect to $\tau$ by a degree map $\nu:M\to\P$, a congruence $S$ on $M$ is called \emph{$\P$-graded} if it is stable under $\nu$; that is, $(u,v)\in S$ implies $\nu(u)=\nu(v)$. This stability guarantees the existence of a well-defined induced map $\nu^{_S}:M/S\to\P$ given by $\nu^{_S}([u]_{_S})=\nu(u)$. Consequently, the quotient monoid $M/S$ inherits a natural $\P$-grading with respect to the induced twisting $\tau^{_S}$. In particular, if $M$ is $\N$-graded with $M_0=\{1\}$, then $(M/S)_0=\{[1]_{_S}\}$.

\begin{pro}\label{146}
Let $S$ be a congruence on $M$ that is stable under $\rho$ and compatible with $\omega$; that is, $\omega(\sigma,u)=\omega(\sigma,v)$ for all $(u,v)\in S$ and $\sigma\in G$. Then:
\begin{enumerate}
\item\label{147}There exists a well-defined weight function $\omega^{_S}:G\times M/S\to\field^\times$ for the induced action $\rho^{_S}$, defined by $\omega^{_S}(\sigma,[u])=\omega(\sigma,u)$ for all $u\in M$ and $\sigma\in G$.
\item\label{148} If $\tau:M\times M\to\field$ is a twisting on $M$ compatible with $S$ and $\rho$ is multiplicative with respect to $\tau$ such that $(\tau,\omega)$ is a multiplicative pair for $\rho$, then $(\tau^{_S},\omega^{_S})$ is a multiplicative pair for $\rho^{_S}$.
\end{enumerate}
\end{pro}
\begin{proof}
\eqref{147}~Because $S$ is compatible with $\omega$, the map $\omega^{_S}:G\times M/S\to\field$ is immediately well-defined. Given $\sigma,\varsigma\in G$ and $u\in M$, we have $\omega^{_S}(1,[u])=\omega(1,u)=1$ and\[\omega^{_S}(\sigma\varsigma,[u])=\omega(\sigma\varsigma,u)=\omega(\sigma,\rho_\varsigma(u))\,\omega(\varsigma,u)=\omega^{_S}(\sigma,[\rho_\varsigma(u)])\,\omega^{_S}(\varsigma,[u])=\omega^{_S}(\sigma,\rho_\varsigma^{_S}([u]))\,\omega^{_S}(\varsigma,[u]).\]Therefore $\omega^{_S}$ is a well-defined weight function. \eqref{148}~Assume now that $\tau$ is compatible with $S$, $\rho$ is multiplicative with respect to $\tau$, and $(\tau,\omega)$ is a multiplicative pair. For any $[u],[v]\in M/S$ with $\tau^{_S}([u],[v])=\tau(u,v)\neq0$, we have\[\tau^{_S}([u],[v])\,\omega^{_S}(\sigma,[u][v])=\tau(u,v)\,\omega(\sigma,uv)=\tau(\rho_\sigma(u),\rho_\sigma(v))\,\omega(\sigma,u)\,\omega(\sigma,v)=\tau^{_S}(\rho_\sigma^{_S}([u]),\rho_\sigma^{_S}([v]))\,\omega^{_S}(\sigma,[u])\,\omega^{_S}(\sigma,[v]).\]Therefore $(\tau^{_S},\omega^{_S})$ is a multiplicative pair.
\end{proof}

\begin{pro}\label{141}
Let $S\subseteq M\times M$ be a congruence on $M$ whose equivalence classes are finite, and let $U$ be a system of representatives for these classes. Assume that $S$ is bigraded and compatible with the twisting; that is, $\bdg(u)=\bdg(u')$ and $\tau(u,v)=\tau(u',v')$ for all $(u,u'),(v,v')\in S$. Then, the assignment $g\mapsto g_{_U}$ defines a bigraded linear endomorphism on $\funcs{\field}{M}$, where $g_{_U}:M\to\field$ is the series defined by
\begin{gather}\label{140}g_{_U}=\sum_{t\in U}\Bigg(\sum_{u\in[t]}g(u)\Bigg)t.\end{gather}
\end{pro}
\begin{proof}
Let $g,h\in\funcs{\field}{M}$. First, since each equivalence class $[t]$ is finite, the inner sum in~\eqref{140} belongs to $\field$; thus, $g_{_U}$ is a well-defined formal power series. Since $g$ is of bounded degree, there exists $D\in\N$ such that $g(u)=0$ for all $u\in M$ with $\deg(u)>D$. Let $t\in U$ with $\deg(t)>D$. Since $S$ preserves the bigrading, every $u\in[t]$ satisfies $\deg(u)=\deg(t)>D$, which forces $g_{_U}(t)=0$. Hence, $g_{_U}$ is of bounded degree.

To check linearity, let $a\in\field$ and $t\in U$. Then\[(ag+h)_{_U}(t)=\sum_{u\in[t]}(ag+h)(u)=a\sum_{u\in[t]}g(u)+\sum_{u\in[t]}h(u)=ag_{_U}(t)+h_{_U}(t).\]Assume now that $\supp(g)\subseteq\ubar{M}_x\cap M^d$ for some $(x,d)\in\G\times\N$. Since $S$ preserves the bigrading, every $t\in U$ satisfies $\bdg(u)=\bdg(t)$ for all $u\in[t]$. Consequently, if $\bdg(t)\neq(x,d)$, then $g(u)=0$ for all $u\in[t]$, and hence $g_{_U}(t)=0$. Therefore, $\supp(g_{_U})\subseteq U\cap\ubar{M}_x\cap M^d$, which implies that the assignment is bigraded.
\end{proof}

Note that if $g\in\funcs{\field}{M}$ is $G$-invariant, its image $g_{_U}$ is not necessarily $G$-invariant. Indeed, since the action $\rho$ does not preserve $U$ in general, for a given $t\in U$ with $g_{_U}(t)\neq0$, there may exist $\sigma\in G$ such that $\rho_\sigma(t)\notin U$. In such a case, the definition forces $g_{_U}(\rho_\sigma(t))=0$, which violates the $G$-invariance condition in Remark~\ref{044}. Consequently, the linear map $g\mapsto g_{_U}$ does not necessarily restrict to an endomorphism on the subspace $\funcs{\field}{M}^G$.

A pair $(S,U)$ satisfying the conditions of Proposition~\ref{141} will be called a \emph{retraction pair} on $M$.

\begin{lem}\label{142}
Let $(S,U)$ be a retraction pair on $M$. Then, Conditions~\eqref{015}~to~\eqref{017} hold for the quotient $M/S$.
\end{lem}
\begin{proof}
Let $\ssf:M\to M/S$ be the canonical projection defined by $\ssf(u)=[u]$. Since the twisting $\tau$ is constant on the congruence classes modulo $S$, Proposition~\ref{073}\eqref{074} ensures the existence of a twisting $\tau^{_S}:M/S\times M/S\to\field$ defined by $\tau^{_S}([u],[v])=\tau(u,v)$ for all $u,v\in M$. As explained above, since $S$ is $\G$-graded, $M/S$ inherits a $\G$-grading $\{\ssf(\ubar{M}_x)\mid x\in\G\}$ such that $|[u]|=|u|$. For $u\in M$, we have $[v][w]=[u]$ in $M/S$ with $\tau^{_S}([v],[w])\neq0$ if and only if $vw\in[u]$ with $\tau(v,w)\neq0$. Because each congruence class $[u]$ is finite and $M$ is locally finite with respect to $\tau$, each element $z\in[u]$ admits only finitely many factorizations $(v,w)\in M\times M$ such that $vw=z$ and $\tau(v,w)\neq0$. Consequently, the set of all factorizations $(v,w)\in M\times M$ such that $vw\in[u]$ and $\tau(v,w)\neq0$ is a finite union of finite sets, and thus finite. Since the set of pairs $\{([v],[w])\in M/S\times M/S\mid [v][w]=[u]\text{ and }\tau^{_S}([v],[w])\neq0\}$ is the image of this finite set under $\ssf\times\ssf$, it is also finite. Therefore, the quotient $M/S$ is locally finite with respect to $\tau^{_S}$. Similarly, since $S$ is $\N$-graded, there is an induced $\N$-grading $\{\ssf(M^d)\mid d\in\N\}$ such that $\deg([u])=\deg(u)$ and $\ssf(M^0)=\{[1]\}$. Because $S$ is bigraded, we have $\ssf(\ubar{M}_x)\cap\ssf(M^d)=\ssf(\ubar{M}_x\cap M^d)$; indeed, if $\ssf(u)=\ssf(v)$ for some $u\in\ubar{M}_x$ and $v\in M^d$, then $(u,v)\in S$, which implies $u,v\in\ubar{M}_x\cap M^d$. In particular, two elements of different $\G$-degree can never be identified under $S$. Consequently, the component $\ssf(M^d)$ intersects the exact same number of fibres $\ssf(\ubar{M}_x)$ as $M^d$ does with $\ubar{M}_x$. Thus, Conditions~\eqref{015}~to~\eqref{017} hold.
\end{proof}

\begin{lem}\label{156}
Let $(S,U)$ be a retraction pair on $M$, and let $\ssf:M\to M/S$ be the canonical projection defined by $\ssf(u)=[u]$ for all $u\in M$. Assume that $S$ is stable under the action and compatible with the weight function; that is, $(\rho_\sigma(u),\rho_\sigma(v))\in S$ and $\omega(\sigma,u)=\omega(\sigma,v)$ for all $(u,v)\in S$ and $\sigma\in G$. Then
\begin{enumerate}
\item\label{143}The induced action $\rho^{_S}$ is bigraded with weight function $\omega^{_S}$, and Condition~\eqref{051} holds for $M/S$. Moreover, if $\rho$ is multiplicative with respect to $\tau$, then Condition~\eqref{018} holds on $M/S$.
\item\label{150}There is a surjective bigraded algebra homomorphism $[\,]:\funcs{\field}{M}\to\funcs{\field}{M/S}$ such that\[[\funcs{\field}{M}^G]\subseteq\funcs{\field}{M/S}^G.\]
\item\label{151}For every $g\in\funcs{\field}{M/S}^G$ there exists $\bar{g}\in\funcs{\field}{M}^G$ such that $g([u])=|[u]|\,\bar{g}(u)$ for all $u\in M$. Consequently, $\funcs{\field}{M/S}^G=[\funcs{\field}{M}^G]$ is a bigraded subalgebra of $\funcs{\field}{M/S}$.
\item\label{144}Assume, additionally, that $S$ is compatible with the coproduct structure; that is, $(\ssf\times\ssf)(\psf(u))=(\ssf\times\ssf)(\psf(v))$ for all $(u,v)\in S$. Then, Condition~\eqref{127} holds.
\end{enumerate}
In consequence, $\funcs{\field}{M/S}^G$ is a connected bialgebra with monomial basis $\{m_{[t]}\mid [t]\in T^{_S}\}$, where $T^{_S}$ is a system of representatives of the admissible $G$-orbits.
\end{lem}
\begin{proof}
\eqref{143}~Since $S$ is stable under the action of $G$ and compatible with $\omega$, Propositions~\ref{073} and~\ref{146} imply that $\rho^{_S}:G\to\sym(M/S)$ is a bigraded group action with weight function $\omega^{_S}:G\times M/S\to\field^\times$, and that Condition~\eqref{018} holds on $M/S$ whenever $\rho$ is multiplicative with respect to $\tau$. Let $G[u]$ be an $\omega^{_S}$-admissible orbit in $M/S$, then $\omega(\sigma,u)=\omega^{_S}(\sigma,[u]_S)=1$ for all $\sigma\in G_{([u])}$. Since $G_{(u)}$ is a subgroup of $G_{([u])}$ we conclude that $Gu$ is $\omega$-admissible. Furthermore, $\ssf(Gu)=\{[\rho_\sigma(u)]\mid\sigma\in G\}=\{\rho^{_S}_\sigma([u])\mid\sigma\in G\}=G[u]$. Consequently, every $\omega^{_S}$-admissible $G$-orbit is the projection under $\ssf$ of an $\omega$-admissible $G$-orbit. Notice that $G[u]\cap\ssf(M^d)=\ssf(Gu\cap M^d)$; indeed, if $[v]\in G[u]\cap\ssf(M^d)$, then $[v]=\rho^{_S}_\sigma([u])=[\rho_\sigma(u)]$ for some $\sigma\in G$, and $[v]=\ssf(w)=[w]$ for some $w\in M^d$. Thus, $[\rho_\sigma(u)]=[w]$, meaning that $(\rho_\sigma(u),w)\in S$. Since $S$ is $\N$-graded and $w\in M^d$, it follows that $\rho_\sigma(u)\in M^d$. Therefore, $\rho_\sigma(u)\in Gu\cap M^d$, which proves that $[v]=\ssf(\rho_\sigma(u))\in\ssf(Gu\cap M^d)$. Since $M^d$ intersects only finitely many $\omega$-admissible $G$-orbits, we obtain that $\ssf(M^d)$ intersects only finitely many $\omega^{_S}$-admissible $G$-orbits. Thus, Condition~\eqref{051} holds.

\eqref{150}~Since $\ssf$ is a monoid homomorphism, Proposition~\ref{019}\eqref{022} ensures the existence of a surjective bigraded algebra homomorphism $[\,]:\funcs{\field}{M}\to\funcs{\field}{M/S}$ defined for any $h\in\funcs{\field}{M}$ by\begin{equation}\label{152}[h]=\sum_{[u]\in M/S}\Bigg(\sum_{\ssf(v)=[u]}h(v)\Bigg)[u]=\sum_{u\in U}\Bigg(\sum_{v\in[u]}h(v)\Bigg)[u].\end{equation}Furthermore, $[\funcs{\field}{M}^G]\subseteq\funcs{\field}{M/S}^G$; indeed for $h\in\funcs{\field}{M}^G$, we have\[\begin{array}{rclclcl}
\bm{\rho}_\sigma^{_S}([h])([u])&=&\omega^{_S}(\sigma,\rho^{_S}_{\sigma^{-1}}([u]))\,[h](\rho^{_S}_{\sigma^{-1}}([u]))\\
&=&\omega^{_S}(\sigma,[\rho_{\sigma^{-1}}(u)])\,[h]([\rho_{\sigma^{-1}}(u)])&=&{\displaystyle\omega^{_S}(\sigma,[\rho_{\sigma^{-1}}(u)])\!\!\sum_{v\in[\rho_{\sigma^{-1}}(u)]}h(v)}\\
&&&=&{\displaystyle\sum_{w\in[u]} \omega(\sigma,\rho_{\sigma^{-1}}(w))\,h(\rho_{\sigma^{-1}}(w))}\\[0.5cm]
&&&=&{\displaystyle\sum_{w\in[u]} \bm{\rho}_\sigma(h)(w)}\,\,\,=\,\,\,[\bm{\rho}_\sigma(h)]([u])\,\,\,=\,\,\,[h]([u]).
\end{array}\]

\eqref{151}~Since $(S,U)$ is a retraction pair, each $[u]\in M/S$ is finite and non-empty, thus, for every $g\in\funcs{\field}{M/S}^G$ we define a formal power series $\bar{g}:M\to\field$ by $\bar{g}(u)=g([u])/|[u]|$. Note that $|[\rho_{\sigma^{-1}}(u)]|=|[u]|$ for all $\sigma\in G$, because $\rho_{\sigma^{-1}}$ is a bijection. Now, since $g$ is $G$-invariant, we have\[\begin{array}{rclcl}
\bm{\rho}_\sigma(\bar{g})(u)&\!\!=\!\!&\omega(\sigma,\rho_{\sigma^{-1}}(u))\,\bar{g}(\rho_{\sigma^{-1}}(u))&\!\!=\!\!&\omega^{_S}(\sigma,[\rho_{\sigma^{-1}}(u)])\,g([\rho_{\sigma^{-1}}(u)])/|[\rho_{\sigma^{-1}}(u)]|\\
&&&\!\!=\!\!&\omega^{_S}(\sigma,\rho^{_S}_{\sigma^{-1}}([u]))\,g(\rho^{_S}_{\sigma^{-1}}([u]))/|[u]|\,\stackrel{\eqref{043}}{=}\,g([u])/|[u]|\,=\,\bar{g}(u).\end{array}\]Finally, applying $[\,]$, we obtain\[[\bar{g}]([u])=\sum_{v\in[u]}\bar{g}(v)=\sum_{v\in[u]}g([v])/|[v]|=|[u]|\,g([u])/|[u]|=g([u]).\]Together with~\eqref{150}, this shows that $\funcs{\field}{M/S}^G=[\funcs{\field}{M}^G]$ is a bigraded subalgebra of $\funcs{\field}{M/S}$.

\eqref{144}~By assumption, the map $\psf^{_S}:M/S\to M/S\times M/S$ given by $\psf^{_S}([u])=(\ssf\times\ssf)(\psf(u))$ is well-defined. Furthermore, since $\ssf$ is a monoid homomorphism, we have\begin{gather*}\psf^{_S}([u][v])=\psf^{_S}([uv])=(\ssf\times\ssf)(\psf(uv))=(\ssf\times\ssf)(\psf(u)\psf(v))=(\ssf\times\ssf)(\psf(u))\,(\ssf\times\ssf)(\psf(v))=\psf^{_S}([u])\,\psf^{_S}([v])\\\text{and}\quad\psf^{_S}([1])=(\ssf\times\ssf)(\psf(1))=(\ssf\times\ssf)(1,1)=([1],[1]).\end{gather*}Thus, $\psf^{_S}$ is a monoid homomorphism. To verify that $\psf^{_S}$ is twisting-preserving and bigraded, let $u,v\in M$ such that $\psf(u)=(u_1,u_2)$ and $\psf(v)=(v_1,v_2)$. Note that $\psf^{_S}([u])=([u_1],[u_2])$ and $\psf^{_S}([v])=([v_1],[v_2])$. Since $S$ is bigraded and $\psf$ is twisting-preserving and bigraded, we have\[\begin{array}{rclcl}
\tau^{_S}_\phi(\psf^{_S}([u]),\psf^{_S}([v]))
&=&\phi(|[u_2]|,|[v_1]|)\,\tau^{_S}([u_1],[v_1])\,\tau^{_S}([u_2],[v_2])\\
&=&\phi(|u_2|,|v_1|)\,\tau(u_1,v_1)\,\tau(u_2,v_2)\,\,\,=\,\,\,\tau_\phi(\psf(u),\psf(v))\,\,\,=\,\,\,\tau(u,v)\,\,\,=\,\,\,\tau^{_S}([u],[v])
\end{array}\]\[\text{and}\quad\bdg(\psf^{_S}([u]))=\bdg([u_1],[u_2])=\bdg(u_1,u_2)=\bdg(\psf(u))=\bdg(u)=\bdg([u]).\]To verify that $\psf^{_S}$ is coassociative, first observe that\begin{equation}\label{079}(\psf^{_S}\times\id)\circ\psf^{_S}\circ\ssf=(\ssf\times\ssf\times\ssf)\circ(\psf\times\id)\circ\psf\quad\text{and}\quad(\id\times\psf^{_S})\circ\psf^{_S}\circ\ssf=(\ssf\times\ssf\times\ssf)\circ(\id\times\psf)\circ\psf.\end{equation}Let $g\in\funcs{\field}{M/S}^G$ and $([a],[b],[c])\in M/S\times M/S\times M/S$. Since $\bar{g}(u)=\bar{g}(v)$ for all $(u,v)\in S$, we have\[\sum_{[v]\in((\psf^{_S}\times\id)\circ\psf^{_S})^{-1}(\{([a],[b],[c])\})}g([v])\,\,=\,\,\sum_{[v]\in((\psf^{_S}\times\id)\circ\psf^{_S})^{-1}(\{([a],[b],[c])\})}|[v]|\bar{g}(v)\,\,\refeq{079}\,\,\sum_{w\in((\psf\times\id)\circ\psf)^{-1}([a]\times[b]\times[c])}\bar{g}(w),\]and similarly,\[\sum_{[v]\in((\id\times\psf^{_S})\circ\psf^{_S})^{-1}(\{([a],[b],[c])\})}g([v])\,\,=\,\,\sum_{[v]\in((\id\times\psf^{_S})\circ\psf^{_S})^{-1}(\{([a],[b],[c])\})}|[v]|\bar{g}(v)\,\,\refeq{079}\,\,\sum_{w\in((\id\times\psf)\circ\psf)^{-1}([a]\times[b]\times[c])}\bar{g}(w).\]Since $\psf$ is coassociative, Lemma~\ref{104}\eqref{116} guarantees that both equations yield the same result, which implies that $\psf^{_S}$ is coassociative.

Note that $\psf^{_S}([u])=(\ssf\times\ssf)(\psf(u))=([v],[w])$ if and only if $\psf(u)\in(\ssf\times\ssf)^{-1}(\{([v],[w])\})=[v]\times[w]$, then\begin{equation}\label{157}(\psf^{_S})^{-1}(\{([v],[w])\})=\ssf(\psf^{-1}([v]\times[w])).\end{equation}Since $[v],[w]$ are finite sets, the product $[v]\times[w]$ is also finite. Consequently, $\ssf(\psf^{-1}([v]\times[w]))$ is finite because $\psf$ has finite fibres, proving that $\psf^{_S}$ has finite fibres as well. Since $(\psf,\hsf)$ is a coproduct structure, for $\sigma,\varsigma\in G$, we have\[\begin{array}{rclclcl}\psf^{_S}(\rho^{_S}_{\hsf(\sigma,\varsigma)}([u]))&\!\!=\!\!&\psf^{_S}([\rho_{\hsf(\sigma,\varsigma)}(u)])&\!\!=\!\!&(\ssf\times\ssf)(\psf(\rho_{\hsf(\sigma,\varsigma)}(u)))\\&&&\!\!=\!\!&(\ssf\times\ssf)(\rho^\otimes_{(\sigma,\varsigma)}(\psf(u)))\\
&&&\!\!=\!\!&(\ssf\times\ssf)(\rho_\sigma(u_1),\rho_\varsigma(u_2))\,=\,(\rho^{_S}_\sigma([u_1]),\rho^{_S}_\varsigma([u_2]))\,=\,(\rho^{_S})^\otimes_{(\sigma,\varsigma)}(\psf^{_S}([u])).
\end{array}\]
\[\begin{array}{rclclcl}\omega^{_S}(\hsf(\sigma,\varsigma),[u])&=&\omega(\hsf(\sigma,\varsigma),u)&=&\omega^\otimes((\sigma,\varsigma),\psf(u))\\
&&&=&\omega(\sigma,u_1)\,\omega(\varsigma,u_2)\\
&&&=&\omega^{_S}(\sigma,[u_1])\,\omega^{_S}(\varsigma,[u_2])&=&(\omega^{_S})^\otimes((\sigma,\varsigma),\psf^{_S}([u])).\end{array}\]Finally, given $g\in\funcs{\field}{M/S}^G$, since $\bar{g}(u)=\bar{g}(v)$ for all $(u,v)\in S$, we have
\[\begin{array}{rclcl}
{\displaystyle\sum_{[v]\in(\psf^{_S})^{-1}(\{([1],[u])\})}\!\!\!\!\!g([v])}&=&{\displaystyle\sum_{[v]\in(\psf^{_S})^{-1}(\{([1],[u])\})}\!\!\!|[v]|\,\bar{g}(v)}\\
&=&{\displaystyle\sum_{[v]\in(\psf^{_S})^{-1}(\{([1],[u])\})}\,\,\sum_{w\in[v]}\bar{g}(w)} \\[0.5cm]
&=&{\displaystyle\sum_{w\in\psf^{-1}([1]\times[u])}\!\!\bar{g}(w)}\\[0.5cm]
&=&{\displaystyle\sum_{y\in[u]}\,\,\sum_{w\in\psf^{-1}(\{(1,y)\})}\!\!\!\!\bar{g}(w)}\,\,\,\refeq{122}\,\,\,{\displaystyle\sum_{y\in[u]}\bar{g}(y)}\,\,\,=\,\,\,|[u]|\,\bar{g}(u)\,\,\,=\,\,\, g([u]).
\end{array}\]Similarly we show that\[\sum_{[v]\in(\psf^{_S})^{-1}(\{([u],[1])\})}\!\!\! g([v])=g([u]),\]which proves~\eqref{122} for $\psf^{_S}$. Therefore, Condition~\eqref{127} holds.
\end{proof}

\begin{lem}\label{160}
Let $(S,U)$ be a retraction pair on $M$ such that $S$ is stable under the action of $G$ and compatible with the weight function, and let $t\in T$. If $[t]$ is $\omega^{_S}$-admissible in $M/S$, then $m_t(v)=1$ for all $v\in[t]\cap Gt$, and\[[m_t]_{_G}=c_t\,m_{[t]}\,,\quad\text{where}\quad c_t=|[t]\cap Gt|\geq1.\]Otherwise, if $[t]$ is not $\omega^{_S}$-admissible in $M/S$, then $[m_t]_G=0$ and we set $c_t=0$.
\end{lem}
\begin{proof}
By definition, $[m_t]_{_G}([u])$ is the sum of all $m_t(v)$ with $v\in[u]$. Since $S$ is stable under the action of $G$, it follows that $[Gt]=G[t]$ for all $t\in T$. If $[u]\notin G[t]$, then $[u]\cap Gt=\emptyset$. Consequently, $m_t(v)=0$ for all $v\in [u]$, forcing $[m_t]_{_G}([u])=0$. In particular, $\supp([m_t]_{_G})\subseteq G[t]$. Furthermore, since $m_t\in\funcs{\field}{M}^G$, Lemma~\ref{156}\eqref{150} implies that $[m_t]_{_G}\in\funcs{\field}{M/S}^G$. If $[t]$ is not $\omega^{_S}$-admissible, Proposition~\ref{045} implies that $[m_t]_{_G}=0$. Assume now that $[t]$ is $\omega^{_S}$-admissible. For any $v=\rho_\sigma(t)\in[t]\cap Gt$, we have $\rho_\sigma^{_S}([t])=[v]=[t]$, which means $\sigma\in G_{([t])}$, and thus $m_t(v)=\omega(\sigma,t)=\omega^{_S}(\sigma,[t])=1$. Evaluating at the representative class $[t]$, since $\supp(m_t)=Gt$, we have\[[m_t]_{_G}([t])=\sum_{v\in[t]}m_t(v)=\sum_{v\in[t]\cap Gt}1=|[t]\cap Gt|=c_t=c_t\,m_{[t]}([t]).\]Then, by Remark~\ref{044}, for any class $[u]=\rho_\sigma^{_S}([t])\in G[t]$, we have\[[m_t]_{_G}([u])=[m_t]_{_G}(\rho_\sigma^{_S}([t]))\refeq{043}\omega^{_S}(\sigma,[t])\,[m_t]_{_G}([t])=c_t\,\omega^{_S}(\sigma,[t])=c_t\,m_{[t]}([u]).\]We conclude that $[m_t]_{_G}=c_t\,m_{[t]}$, completing the proof.
\end{proof}

Let $I\subseteq\funcs{\field}{M}$ be the kernel of the linear endomorphism $g\mapsto g_{_U}$; that is,\[I=\{g\in\funcs{\field}{M}\mid g_{_U}=0\}=\{g-h\mid g_{_U}=h_{_U}\}.\] 
Notice that $g_{_U}=0$ if and only if $[g]=0$, which means that $I=\ker([\,])$. Consequently, since $[\,]:\funcs{\field}{M}\to\funcs{\field}{M/S}$ is an algebra homomorphism by Lemma~\ref{156}\eqref{150}, $I$ is a two-sided ideal of $\funcs{\field}{M}$.

\begin{lem}
Let $(S,U)$ be a retraction pair on $M$. Then, the map $\dbrack{\,}:\funcs{\field}{M}/I\to\funcs{\field}{M/S}$ defined by $\dbrack{g+I}=[g]$ is a bigraded algebra isomorphism, where $[\,]$ is defined in~\eqref{152}.
\end{lem}
\begin{proof}
By Lemma~\ref{156}\eqref{150}, the map $[\,]:\funcs{\field}{M}\to\funcs{\field}{M/S}$ is a surjective bigraded algebra homomorphism. Since $\ker([\,])=I$, by the first isomorphism theorem for algebras, the induced map $\dbrack{\,}:\funcs{\field}{M}/I\to\funcs{\field}{M/S}$ given by $\dbrack{g+I}=[g]$ is a well-defined algebra isomorphism. Finally, to see that $\dbrack{\,}$ is bigraded, recall that since $\funcs{\field}{M}$ is $(\G\times\N)$-graded, every function $g$ can be uniquely expressed as a sum of functions $g_x^d\in\funcs{\field}{M_x^d}$. Notice that $(g_x^d)_{_U}\in\funcs{\field}{M_x^d}$ because $\supp((g_x^d)_{_U})\subseteq U\cap M_x^d$. Furthermore, since $S$ is bigraded, if $t\in U$ with $\bdg(t)=(x,d)$, we have\[g_{_U}(t)=\sum_{u\in[t]}g(u)=\sum_{u\in[t]}g_x^d(u)=(g_x^d)_{_U}(t)=\sum_{(y,e)\in\G\times\N}(g_y^e)_{_U}(t).\]Now, if $g\in I$, then $0=g_{_U}=(g_x^d)_{_U}$, that is, $g_x^d\in I$ for all $(x,d)\in\G\times\N$, which implies that
\[\funcs{\field}{M}/I=\bigoplus_{(x,d)\in\G\times\N}\funcs{\field}{M_x^d}/I_x^d,\quad\text{where}\quad I_x^d=I\cap\funcs{\field}{M_x^d}.\]Since $S$ is bigraded, given $g\in\funcs{\field}{M_x^d}$ and $t\in M_y^e$ with $(y,e)\neq(x,d)$, we have
\[\dbrack{g+I}([t])=[g]([t])=\sum_{u\in[t]}g(u)=0,\]which implies that $\dbrack{g+I}\in\funcs{\field}{(M/S)_x^d}$. Therefore, $\dbrack{\,}$ is bigraded, concluding the proof.
\end{proof}

A retraction pair $(S,U)$ satisfying all the conditions of Lemma~\ref{156} will be called an \emph{invariant-retraction pair}. Explicitly, this means that the following conditions hold:
\begin{enumerate}
\item $S$ is a congruence on $M$ whose equivalence classes are finite, and $U$ is a system of representatives for these classes.
\item $S$ is bigraded and compatible with the twisting; that is, $\bdg(u)=\bdg(u')$ for all $(u,u')\in S$, and $\tau(u,v)=\tau(u',v')$ for all $(u,u'),(v,v')\in S$.
\item $S$ is stable under the action of $G$ and compatible with the weight function; that is, $(\rho_\sigma(u),\rho_\sigma(v))\in S$ and $\omega(\sigma,u)=\omega(\sigma,v)$ for all $(u,v)\in S$ and $\sigma\in G$.
\item $S$ is compatible with the coproduct structure; that is, $(\ssf\times\ssf)(\psf(u))=(\ssf\times\ssf)(\psf(v))$ for all $(u,v)\in S$.
\end{enumerate}

Throughout the remainder of this subsection, we assume the following standing structural condition.
\stepcounter{condition}
\begin{enumerate}[start=\thecondition,label={(\condition)},ref={\condition}]
\item\label{161}\conditions{6}
\end{enumerate}

\begin{pro}\label{154}
The algebra homomorphism $[\,]$ restricts to a surjective bigraded bialgebra homomorphism\[[\,]_{_G}:\funcs{\field}{M}^G\to\funcs{\field}{M/S}^G.\]
\end{pro}
\begin{proof}
By Lemma~\ref{156}\eqref{150}, we have $[\funcs{\field}{M}^G]=\funcs{\field}{M/S}^G$. Since $[\,]$ is a bigraded algebra homomorphism, its restriction $[\,]_{_G}:\funcs{\field}{M}^G\to\funcs{\field}{M/S}^G$ is a well-defined surjective bigraded algebra homomorphism. Hence, it remains only to show that $[\,]_{_G}$ preserves the counit and the coproduct.

Let $\varepsilon$ and $\varepsilon^{_S}$ denote the counits of $\funcs{\field}{M}^G$ and $\funcs{\field}{M/S}^G$, respectively. Given $g\in\funcs{\field}{M}^G$, recall that $\varepsilon(g)=g(1)$ and $\varepsilon^{_S}([g])=[g]([1])$. Since $M$ is graded with $M^0=\{1\}$ and $S$ is a bigraded congruence, the identity element $1$ is the unique element of degree $0$ in $M$, which implies that its congruence class is the singleton $[1]=\{1\}$. Consequently, we have
\[\varepsilon^{_S}([g])=[g]([1])=\sum_{u\in[1]}g(u)=g(1)=\varepsilon(g).\]

To verify that $[\,]_{_G}$ preserves the coproduct, recall from Theorem~\ref{130} that the coproducts are given by $\Delta=\Phi_2\circ\nabla$ and $\Delta^{_S}=\Phi^{_S}_2\circ\nabla^{_S}$, where $\Phi_2$ and $\Phi^{_S}_2$ are the algebra isomorphisms from Corollary~\ref{097}. By Proposition~\ref{019}\eqref{022}, the surjective monoid homomorphism $\ssf\times\ssf$ induces a surjective algebra homomorphism\[[\,]^\otimes:\funcs{\field}{M\times M}^{G\times G}\to\funcs{\field}{(M/S)\times(M/S)}^{G\times G},\quad[g]^\otimes=\sum_{[u],[v]\in M/S}\Bigg(\sum_{(a,b)\in[u]\times[v]}\!\!\!\!g(a,b)\Bigg)([u],[v]).\]

Given $([v],[w])\in(M/S)\times(M/S)$, notice that\begin{equation}\label{155}\bigsqcup_{(a,b)\in[v]\times[w]}\psf^{-1}(\{(a,b)\})=\psf^{-1}([v]\times[w]).\end{equation}We compute\[[\nabla(g)]^\otimes([v],[w])\,\,\,=\sum_{(a,b)\in[v]\times[w]}\nabla(g)(a,b)\,\,\,\refeq{128}\sum_{(a,b)\in[v]\times[w]}\,\sum_{u\in\psf^{-1}(\{(a,b)\})}g(u)\,\,\,\refeq{155}\sum_{u\in\psf^{-1}([v]\times[w])}g(u).\]By assumption, for $u\in\psf^{-1}([v]\times[w])$ and $u'\in[u]$, we have $(\ssf\times\ssf)(\psf(u'))=(\ssf\times\ssf)(\psf(u))=([v],[w])$, which implies that $[u]\subseteq\psf^{-1}([v]\times[w])$. Hence, $u\in\psf^{-1}([v]\times[w])$ if and only if $[u]\in\ssf(\psf^{-1}([v]\times[w]))$. Then\begin{equation}\label{163}[\nabla(g)]^\otimes([v],[w])=\!\!\!\!\sum_{[u]\in\ssf(\psf^{-1}([v]\times[w]))}\,\sum_{u'\in[u]}g(u')\refeq{157}\!\!\!\!\sum_{[u]\in(\psf^{_S})^{-1}(\{([v],[w])\})}\!\!\!\![g]([u])=\nabla^{_S}([g])([v],[w]).\end{equation}Hence, $\Delta^{_S}([g])=\Phi_2^{_S}(\nabla^{_S}([g]))=\Phi_2^{_S}([\nabla(g)]^\otimes)$. We need to show that $\Phi_2^{_S}\circ[\,]^\otimes=([\,]_{_G}\otimes_{\phi}[\,]_{_G})\circ\Phi_2$. We proceed by applying each side on the monomial basis. On one hand, we have\begin{equation}\label{162}([\,]_{_G}\otimes_{\phi}[\,]_{_G})(\Phi_2(m_{(u,v)}^\otimes))=([\,]_{_G}\otimes[\,]_{_G})(m_u\otimes m_v)=[m_u]_{_G}\otimes_{\phi}[m_v]_{_G}.\end{equation}Setting $T_S=\{t\in T\mid [t]\text{ is }\omega^S\text{-admissible}\}$, if $u\notin T_S$ or $v\notin T_S$, Lemma~\ref{160} implies that $[m_u]_{_G}=0$ or $[m_v]_{_G}=0$, so both $[m_{(u,v)}^\otimes]^\otimes$ and $([\ ]_{_G}\otimes_\phi[\ ]_{_G})(\Phi_2(m_{(u,v)}^\otimes))$ vanish. Assume now that $u,v\in T_S$. On the other hand, by applying Lemma~\ref{160} and Proposition~\ref{159}, we have\[\begin{array}{rclcl}
[m_{(u,v)}^\otimes]^\otimes([a],[b])&=&\displaystyle\sum_{(x,y)\in[a]\times[b]}m_{(u,v)}^\otimes(x,y)&=&\displaystyle\sum_{(x,y)\in[a]\times[b]}m_u(x)\,m_v(y)\\[0.5cm]
&&&=&\displaystyle\sum_{x\in[a]}m_u(x)\,\sum_{y\in[b]}m_v(y)\\[0.4cm]
&&&=&[m_u]_{_G}([a])\,[m_v]_{_G}([b])\\
&&&=&c_u\,m_{[u]}([a])\,c_v\,m_{[v]}([b])\,\,\,=\,\,\,c_u\,c_v\,m^\otimes_{([u],[v])}([a],[b]).
\end{array}\]Thus, we obtain $[m_{(u,v)}^\otimes]^\otimes=c_u\,c_v\,m_{([u],[v])}^\otimes$. Applying the isomorphism $\Phi_2^{_S}$ and using Lemma~\ref{160} yields\begin{equation}\label{164}\Phi_2^{_S}\big([m_{(u,v)}^\otimes]^\otimes\big)=c_u\,c_v\Phi_2^{_S}\big(m_{([u],[v])}^\otimes\big)=c_u\,m_{[u]}\otimes c_v\,m_{[v]}=[m_u]_{_G}\otimes[m_v]_{_G}\!\!\!\refeq{162}\!([\,]_{_G}\otimes_{\phi}[\,]_{_G})(\Phi_2(m_{(u,v)}^\otimes)).\end{equation}Finally, applying this identity to $\nabla(g)$, we conclude
\[\Delta^{_S}([g])\refeq{163}\Phi_2^{_S}\big([\nabla(g)]^\otimes\big)\refeq{164}([\,]_{_G}\otimes_{\phi}[\,]_{_G})\big(\Phi_2(\nabla(g))\big)=([\,]_{_G}\otimes_{\phi}[\,]_{_G})(\Delta(g)),\]
which proves that $[\,]_{_G}$ preserves the coproduct, concluding the proof.
\end{proof}

\begin{thm}\label{153}
The intersection $J=I\cap\funcs{\field}{M}^G$ is a graded bi-ideal of $\funcs{\field}{M}^G$, such that $\funcs{\field}{M}^G/J$ and $\funcs{\field}{M/S}^G$ are isomorphic connected bialgebras.
\end{thm}
\begin{proof}
By Proposition~\ref{154}, the restriction $[\,]_{_G}:\funcs{\field}{M}^G\to\funcs{\field}{M/S}^G$ is a well-defined surjective bigraded bialgebra homomorphism. Observe that the kernel of this homomorphism is precisely $J$; indeed\[\ker([\,]_{_G})=\ker([\,])\cap\funcs{\field}{M}^G=I\cap\funcs{\field}{M}^G=J.\]Since $[\,]_{_G}$ is a bigraded bialgebra homomorphism, Proposition~\ref{165} implies that $J$ is a graded bi-ideal of $\funcs{\field}{M}^G$ and that $\funcs{\field}{M}^G/J$ and $\funcs{\field}{M/S}^G$ are isomorphic connected bialgebras.
\end{proof}

\begin{pro}\label{166}
Let $T_S=\{t\in T\mid [t]\text{ is }\omega^{_S}\text{-admissible}\}$. The graded bi-ideal $J$ is generated by the set\[\{m_u\mid u\in T\setminus T_S\}\cup\{c_v\,m_u-c_u\,\omega^{_S}(\sigma,[u])\,m_v\mid u,v\in T_S\text{ and }\sigma\in G\text{ such that }\rho_\sigma^{_S}([u])=[v]\}.\]
\end{pro}
\begin{proof}
Let $K$ be the ideal of $\funcs{\field}{M}^G$ generated by the proposed set. By Lemma~\ref{160}, we have $[m_u]_{_G}=0$ for all $u\in T\setminus T_S$. Moreover, given $u,v\in T_S$ and $\sigma\in G$ satisfying $\rho_\sigma^{_S}([u])=[v]$, Lemma~\ref{046} applied to $M/S$ yields $m_{[u]}=\omega^{_S}(\sigma,[u])\,m_{[v]}$; hence Proposition~\ref{154} and Lemma~\ref{160} imply that\[[c_v\,m_u-c_u\,\omega^{_S}(\sigma,[u])\,m_v]_{_G}=c_v[m_u]_{_G}-c_u\,\omega^{_S}(\sigma,[u])[m_v]_{_G}=c_uc_v\,m_{[u]}-c_uc_v\,\omega^{_S}(\sigma,[u])\,m_{[v]}=0.\]Hence, $K\subseteq J$. To prove the reverse inclusion, consider $g\in J$. Since the monomial invariant functions form a linear basis for $\funcs{\field}{M}^G$, we can express $g$ as\[g=g'+\sum_{t\in T_S}a_t\,m_t\quad\text{with}\quad g'=\sum_{t\in T\setminus T_S}a_t\,m_t\,,\]where all but finitely many coefficients $a_t$ are zero. Let $\sim$ be the equivalence relation on $T_S$ defined by $u\sim v$ if and only if $G[u]=G[v]$ in $M/S$. For each equivalence class $E\in T_S/{\sim}$, fix a representative $t_E\in E$, which satisfies $c_{t_E}\neq0$ by Lemma~\ref{160}, and for each $t\in E$ choose $\sigma_t\in G$ such that $\rho_{\sigma_t}^{_S}([t])=[t_E]$, so that $m_{[t]}=\lambda_t\,m_{[t_E]}$ with $\lambda_t:=\omega^{_S}(\sigma_t,[t])\in\field^\times$. We rewrite $g$ as
\[g=g'+\sum_{E\in T_S/{\sim}}\,\sum_{t\in E}\frac{a_t}{c_{t_E}}(c_{t_E}\,m_t-c_t\,\lambda_t\,m_{t_E})+\sum_{E\in T_S/{\sim}}\Bigg(\sum_{t\in E}a_t\frac{c_t\,\lambda_t}{c_{t_E}}\Bigg)m_{t_E}.\]Since $[g']_{_G}=0$ and $[c_{t_E}\,m_t-c_t\,\lambda_t\,m_{t_E}]_{_G}=0$, applying $[\ ]_{_G}$ yields\[0=[g]_{_G}=\sum_{E\in T_S/{\sim}}\Bigg(\sum_{t\in E}a_t\frac{c_t\,\lambda_t}{c_{t_E}}\Bigg)[m_{t_E}]_{_G}=\sum_{E\in T_S/{\sim}}\Bigg(\sum_{t\in E}a_tc_t\lambda_t\Bigg)m_{[t_E]}\,.\]Since distinct classes $E\in T_S/{\sim}$ correspond to distinct $\omega^{_S}$-admissible $G$-orbits in $M/S$, the set $\{m_{[t_E]}\mid E\in T_S/{\sim}\}$ is linearly independent by Proposition~\ref{048}, forcing $\sum_{t\in E}a_tc_t\lambda_t=0$ for every $E\in T_S/{\sim}$. Thus,\[g=g'+\sum_{E\in T_S/{\sim}}\,\sum_{t\in E}\frac{a_t}{c_{t_E}}(c_{t_E}\,m_t-c_t\,\lambda_t\,m_{t_E})\in K.\]Therefore, $J\subseteq K$, which concludes the proof.
\end{proof}

\section{Symmetric functions}\label{169}

Let $X=(x_r\mid r\in\Q)$ be an infinite alphabet linearly ordered by the natural order of its indices, and let $X^*$ be the free monoid generated by $X$.

\subsection{Classical symmetric functions}\label{179}

For every $q\in\field^\times$, let $\phi_q:\N\times\N\to\field^\times$ be the bicharacter defined by $\phi_q(n,m)=q^{nm}$.

\subsubsection{Symmetric functions in noncommuting variables}\label{123}

Given a word $u=x_{r_1}\cdots x_{r_k}\in X^*$, an \emph{inversion} of $u$ is a pair of indices $(i,j)$ such that $i<j$ and $r_i>r_j$. The total number of inversions of $u$ is denoted by $\inv(u)$. We denote by $\sort(u)$ the unique word obtained by rearranging the letters of $u$ such that its indices are in nondecreasing order. Observe that $\sort(u)$ is the result of repeatedly exchanging letters $x_{r_i},x_{r_j}$ whenever $(i,j)$ is an inversion until no more inversions remain; consequently, $\inv(\sort(u))=0$. For any $u,v\in X^*$, a \emph{cross-inversion} of $(u,v)$ is an inversion of the word $\sort(u)\,\sort(v)$. The number of cross-inversions of $(u,v)$ is denoted by $\inv(u,v)$. Furthermore, note that for any $u,v,w\in X^*$, we have
\begin{gather}
\label{024}\inv(uv,w)=\inv(u,w)+\inv(v,w)\quad\text{and}\quad\inv(u,vw)=\inv(u,v)+\inv(u,w),\\
\label{026}\inv(uv)=\inv(u)+\inv(v)+\inv(u,v).
\end{gather}

Let $\tau_q:X^*\times X^*\to\field^\times$ be defined by $\tau_q(u,v)=q^{\inv(u,v)}$. Note that $\tau_q$ defines a twisting on $X^*$, called \emph{$q$-twisting}; indeed, by~\eqref{024}, it immediately follows that
\begin{equation}\label{025}
\inv(u,v)+\inv(uv,w)=\inv(u,v)+\inv(u,w)+\inv(v,w)=\inv(u,vw)+\inv(v,w),
\end{equation}
which implies the 2-cocycle condition $\tau_q(u,v)\,\tau_q(uv,w)=\tau_q(u,vw)\,\tau_q(v,w)$. Moreover, the normalization condition $\tau_q(u,1)=1=\tau_q(1,u)$ holds trivially since the empty word produces no inversions. By~\eqref{025}, the assignment $(u,v)\mapsto\inv(u,v)$ defines a \emph{twisting} in the sense of~\cite[Definition~3.1]{EaGrMuRu26}.

In addition, the free monoid $X^*$ is locally finite for any twisting because every word $w\in X^*$ has finite length $\ell(w)$, yielding only finitely many factorizations $uv=w$. We equip $X^*$ with the $\N$-grading given by word length, setting $\ubar{X}^*_d=\{u\in X^*\mid\ell(u)=d\}$ for each $d\in\N$, which ensures $\ubar{X}^*_0=\{1\}$. The second grading is defined identically, setting $(X^*)^d=\ubar{X}^*_d$ for all $d\in\N$. Consequently, Conditions~\eqref{015},~\eqref{016}~and~\eqref{017} hold for $X^*$.

Let $\sym_\Q$ be the \emph{finitary symmetric group} on $\Q$, that is, the group of permutations of $\Q$ that fix all but finitely many elements. The \emph{symmetrizing action} of $\sym_\Q$ on $X^*$ is the multiplicative group action $\rho:\sym_\Q\to\Aut(X^*)$ defined by $\rho_\sigma(u)=x_{\sigma(r_1)}\cdots x_{\sigma(r_k)}$ for $u=x_{r_1}\cdots x_{r_k}\in X^*$, which is clearly bigraded because it preserves word lengths. We equip this action with the map $\omega_q:\sym_\Q\times X^*\to\field^\times$, called \emph{$q$-weight function}, defined by $\omega_q(\sigma,u)=q^{\inv(\rho_\sigma(u))-\inv(u)}$. Note that $\omega_q$ is a well-defined weight function; indeed, for any $\sigma,\varsigma\in\sym_\Q$ and $u\in X^*$, we immediately have
\[\inv(\rho_{\sigma\varsigma}(u))-\inv(u)=\inv(\rho_\sigma(\rho_\varsigma(u)))-\inv(\rho_\varsigma(u))+\inv(\rho_\varsigma(u))-\inv(u),\]
which establishes the 1-cocycle condition $\omega_q(\sigma\varsigma,u)=\omega_q(\sigma,\rho_\varsigma(u))\,\omega_q(\varsigma,u)$, and also $\omega_q(1,u)=q^{\inv(u)-\inv(u)}=1$. Moreover, we can verify that $(\tau_q,\omega_q)$ forms a multiplicative pair; indeed,
\begin{gather*}
\inv(\rho_\sigma(uv))-\inv(u)-\inv(v)=\inv(\rho_\sigma(uv))-\inv(uv)+\inv(uv)-\inv(u)-\inv(v),\quad\text{and}\\
\inv(\rho_\sigma(uv))-\inv(u)-\inv(v)=\inv(\rho_\sigma(uv))-\inv(\rho_\sigma(u))-\inv(\rho_\sigma(v))+\inv(\rho_\sigma(u))-\inv(u)+\inv(\rho_\sigma(v))-\inv(v),
\end{gather*}
which implies $\omega_q(\sigma,uv)\,\tau_q(u,v)=\tau_q(\rho_\sigma(u),\rho_\sigma(v))\,\omega_q(\sigma,u)\,\omega_q(\sigma,v)$. Consequently, Condition~\eqref{018} holds.\medskip

Together with these ingredients, we set $\NCSym_q=\funcs{\field}{X^*}^{\sym_\Q}$, the \emph{algebra of $q$-symmetric functions in noncommuting variables}. In particular, $\NCSym:=\NCSym_1$ is the classical algebra of symmetric functions in noncommuting variables~\cite{Wol36,RoSa06}.\medskip

Given $u\in X^*$ and $\sigma\in (\sym_\Q)_{(u)}$, we have $\inv(\rho_\sigma(u))-\inv(u)=\inv(u)-\inv(u)=0$, which immediately yields $\omega_q(\sigma,u)=1$, and thus every $\sym_\Q$-orbit is $\omega_q$-admissible. We choose as our system of representatives the set $T_q\subset\{x_1,x_2,\ldots\}^*\subset X^*$ of all \emph{restricted growth words}~\cite[Subsection~6.1]{BerLau10}, where a word $u=x_{r_1}\cdots x_{r_k}$ is \emph{restricted growth} if either $u=1$, or $r_1=1$ and $r_i\leq\max\{r_1,\ldots,r_{i-1}\}+1$ for all $i\in[2,k]$. Since the symmetrizing action preserves word length, each homogeneous component $(X^*)^d$ is a disjoint union of $\sym_\Q$-orbits. Thus, the number of $\sym_\Q$-orbits intersecting $(X^*)^d$ is precisely the number of restricted growth words of length $d$, which is known to be the $d$th Bell number~\cite{Rio79}~\cite[\href{https://oeis.org/A000110}{A000110}]{OEIS}. Consequently, Condition~\eqref{051} holds, and therefore, by Theorem~\ref{050}, the algebra $\NCSym_q$ admits the monomial basis $\{m_u^q\mid u\in T_q\}$.

To equip the free monoid $X^*$ with a coproduct structure, we apply the method of duplication in Subsubsection~\ref{035}. Note that the linear order of $\Q$ can be extended to $\Q\sqcup\Q$ by declaring that the elements of the first copy are strictly less than the elements of the second copy. With this order, $\Q\sqcup\Q$ is a countable unbounded dense linearly ordered set~\cite[Definition~4.2]{Jec03}. Then, by Cantor's theorem~\cite[Theorem~4.3]{Jec03}, $\Q$ and $\Q\sqcup\Q$ are order-isomorphic, say by $\iota:\Q\to\Q\sqcup\Q$. This guarantees the existence of an order-isomorphism $X\to X\sqcup Y$, where $Y$ is a disjoint ordered copy of $X$ such that $x<y$ for all $(x,y)\in X\times Y$. By the universal property of free monoids~\cite[Subsection~1.6]{How95}, this uniquely extends to a monoid isomorphism\[\psi:X^*\to(X\sqcup Y)^*=X^**Y^*=(X^**Y^*)/R,\]where $R=\{(u,u)\mid u\in (X\sqcup Y)^*\}$ is the trivial congruence on $(X\sqcup Y)^*$. Then, we define $\psf=\qsf\circ\psi$, where $\qsf:(X\sqcup Y)^*\to X^*\times Y^*$ is the canonical projection. Furthermore, we define the group homomorphism $\hsf:\sym_\Q\times\sym_\Q\to\sym_\Q$ by $\hsf(\sigma,\varsigma)=\iota^{-1}\circ(\sigma\sqcup\varsigma)\circ\iota$ for all $\sigma,\varsigma\in\sym_\Q$, where $\sigma\sqcup\varsigma$ acts on $\Q\sqcup\Q$ by applying $\sigma$ to the elements of the first copy and $\varsigma$ to the elements of the second copy.

To verify that $(\psf,\hsf)$ is a coproduct structure, we first introduce some terminology. For a word $u=x_{r_1}\cdots x_{r_k}\in X^*$, we set $\alph(u)=\{x_{r_1},\ldots,x_{r_k}\}$, and for any subset $A\subset X$, we write $u|_A$ for the subword obtained by removing the letters from $u$ that do not belong to $A$. Given $u\in X^*$ and $v\in Y^*$, a \emph{$(u,v)$-shuffle} is a word $w\in(X\sqcup Y)^*$ such that $\alph(w)=\alph(u)\cup\alph(v)$ with $w|_X=u$ and $w|_Y=v$. The set of all such shuffles corresponds to the summands of the classical shuffle product of words~\cite[Section~1.4]{Reu93}, which is known to be finite~\cite[\href{https://oeis.org/A007318}{A007318}]{OEIS}. By definition, the preimage of any pair $(u,v)\in X^*\times Y^*$ under $\qsf$ consists exactly of the $(u,v)$-shuffles, which implies that the canonical projection $\qsf$ has finite fibres.

Recall that the monoid $(X\sqcup Y)^*$ is equipped with the twisting $\tau_q^\psi(u,v)=\tau_q(\psi^{-1}(u),\psi^{-1}(v))$ in Proposition~\ref{036}\eqref{107}, while the direct product $X^*\times Y^*$ is equipped with the combined twisting $(\tau_q)_\phi=\tau_q^\times\tau_q^\phi$ given in Propositions~\ref{083} and~\ref{075}\eqref{009}. Because the isomorphism $\psi^{-1}$ is order-preserving, the letters in $\psi^{-1}(Y)$ are greater than the letters in $\psi^{-1}(X)$. Thus, for any $u,v\in (X\sqcup Y)^*$ with canonical projections $\qsf(u)=(u_1,u_2)$ and $\qsf(v)=(v_1,v_2)$, the cross-inversions between the components are exactly $\inv(u_1,v_2)=0$ and $\inv(u_2,v_1)=\ell(u_2)\ell(v_1)$. Consequently, by expanding the cross-inversions according to~\eqref{024}, the number of inversions expands as\[\begin{array}{rcl}
\inv(\psi^{-1}(u),\psi^{-1}(v))&=&\inv(u_1,v_1)+\inv(u_1,v_2)+\inv(u_2,v_1)+\inv(u_2,v_2)\\
&=&\inv(u_1,v_1)+0+\ell(u_2)\ell(v_1)+\inv(u_2,v_2)\,\,\,\,=\,\,\,\,\inv(u_1,v_1)+\inv(u_2,v_2)+\ell(u_2)\ell(v_1).
\end{array}\]Evaluating the combined twisting on the product space, we obtain\[\begin{array}{rclcl}
(\tau_q)_\phi(\qsf(u),\qsf(v))&=&\tau_q^\times((u_1,u_2),(v_1,v_2))\,\phi_q(\ell(u_2),\ell(v_1))\kern-3em\\[0.1cm]
&=&\tau_q(u_1,v_1)\,\tau_q(u_2,v_2)\,q^{\ell(u_2)\ell(v_1)}&=&q^{\inv(u_1,v_1)+\inv(u_2,v_2)+\ell(u_2)\ell(v_1)}\\[0.1cm]
&&&=&\tau_q(\psi^{-1}(u),\psi^{-1}(v))\,\,\,\,=\,\,\,\,\tau_q^\psi(u,v).
\end{array}\]This shows that $\qsf$ is twisting-preserving. To verify that $(\qsf,\hsf)$ satisfies~\eqref{137}, let $u\in(X\sqcup Y)^*$ with canonical projection $\qsf(u)=(u_1,u_2)$, and let $\sigma,\varsigma\in\sym_\Q$. By definition, $\rho^\psi_{\hsf(\sigma,\varsigma)}$ evaluates as $\rho_\sigma$ on the letters of $X$ and as $\rho_\varsigma$ on the letters of $Y$; that is, it stabilizes the alphabets $X$ and $Y$, which implies\[\qsf(\rho^\psi_{\hsf(\sigma,\varsigma)}(u))=(\rho_\sigma(u_1),\rho_\varsigma(u_2))=\rho^\otimes_{(\sigma,\varsigma)}(\qsf(u)).\]On the other hand, since the letters in $Y$ are greater than the letters in $X$, a cross-inversion in $u$ between a letter $y_i$ and a letter $x_j$ occurs if and only if $y_i$ precedes $x_j$ in $u$. Because $\rho^\psi_{\hsf(\sigma,\varsigma)}$ stabilizes both alphabets $X$ and $Y$, the number of cross-inversions between $y_i$ and $x_j$ in $u$ is invariant under the action of $\hsf(\sigma,\varsigma)$, hence, applying~\eqref{026},\[\inv(\rho^\psi_{\hsf(\sigma,\varsigma)}(u))-\inv(u)=\big(\inv(\rho_\sigma(u_1))-\inv(u_1)\big)+\big(\inv(\rho_\varsigma(u_2))-\inv(u_2)\big).\]Evaluating the weight function $\omega_q^\psi$, we obtain\[\omega_q^\psi(\hsf(\sigma,\varsigma),u)=q^{\inv(\rho^\psi_{\hsf(\sigma,\varsigma)}(u))-\inv(u)}=q^{\inv(\rho_\sigma(u_1))-\inv(u_1)}\,q^{\inv(\rho_\varsigma(u_2))-\inv(u_2)}=\omega_q(\sigma,u_1)\,\omega_q(\varsigma,u_2)=\omega_q^\otimes((\sigma,\varsigma),\qsf(u)).\]Thus, by Proposition~\ref{106}, $\psf=\qsf\circ\psi$ is twisting-preserving with finite fibres, and $(\psf,\hsf)$ satisfies the conditions in~\eqref{125}. Because $\psi$ is induced by a bijection of alphabets, it preserves the length of words, and $\qsf$ partitions the letters of a word into two subwords without altering the total number of letters, for any $u\in X^*$ with $\psf(u)=(u_1,u_2)$, we have $\ell(u)=\ell(u_1)+\ell(u_2)$, which implies that $\psf$ is bigraded.

To verify that $\psf$ is counital, let $t,t'\in T_q$. By definition, $v\in\psf^{-1}(\{(1,t')\})$ if and only if $\qsf(\psi(v))=(1,t')$. This means that $\psi(v)=t'\in Y^*\subset(X\sqcup Y)^*$; hence, $v=\psi^{-1}(t')$ is the unique word in $X^*$ satisfying $\psf(v)=(1,t')$, which implies $\psf^{-1}(\{(1,t')\})=\{v\}$. Furthermore, since $\psi$ is induced by an order-isomorphism, $v$ shares the same relative order of letters as $t'$, and thus $v\in Gt'$. Consequently, $\psf^{-1}(\{(1,t')\})\cap Gt$ is empty whenever $t\neq t'$. On the other hand, if $t=t'$, then $\psf^{-1}(\{(1,t')\})\cap Gt=\{v\}$. Since $v$ and $t$ share the same relative order of letters, we obtain $\inv(v)=\inv(t)$. Then, $\rho_{\sigma_v}(t)=v$ implies $\omega_q(\sigma_v,t)=q^{\inv(v)-\inv(t)}=q^0=1$, hence\[\sum_{v\in\psf^{-1}(\{(1,t')\})\cap Gt}\!\!\omega_q(\sigma_v,t)\,\,\,=\,\,\,\delta_{t,t'}.\]A symmetric evaluation holds for $\psf^{-1}(\{(t',1)\})\cap Gt$, which proves that~\eqref{089} is satisfied. 

To verify that $\psf$ is coassociative, let $u=(u_1,u_2,u_3)\in X^*\times X^*\times X^*$, and consider the preimages\begin{gather}\label{175}A=((\psf\times\id)\circ\psf)^{-1}(\{u\})\quad\text{and}\quad B=((\id\times\psf)\circ\psf)^{-1}(\{u\}).\end{gather}Both composed maps are formed by applying the order-isomorphism $\psi$ to partition the letters of $X$ into three disjoint ordered alphabets $X\sqcup Y\sqcup Z$, where $Z$ is another copy of the alphabet $X$. These compositions respect the order of these components, meaning that letters mapped to $X$ are strictly less than those in $Y$, which in turn are strictly less than those in $Z$. As a result, every word $v\in A$ corresponds to a unique word $w\in B$. This correspondence defines a bijection $f:A\to B$ such that $w=f(v)$ preserves the relative order of the letters. Consequently, $v\in Gw$ and $\inv(v)=\inv(w)$. For any $t\in T_q$, if $v\in Gt$, then $w\in Gt$, which implies $f(A\cap Gt)=B\cap Gt$. We have\[\omega_q(\sigma_v,t)=q^{\inv(v)-\inv(t)}=q^{\inv(w)-\inv(t)}=\omega_q(\sigma_w,t).\]Summing over the bijection $f$, we obtain\[\sum_{v\in((\psf\times\id)\circ\psf)^{-1}(\{u\})\cap Gt}\!\!\omega_q(\sigma_v,t)=\sum_{v\in A\cap Gt}\!\!q^{\inv(v)-\inv(t)}=\sum_{w\in B\cap Gt}\!\!q^{\inv(w)-\inv(t)}=\sum_{w\in((\id\times\psf)\circ\psf)^{-1}(\{u\})\cap Gt}\!\!\omega_q(\sigma_w,t).\]This satisfies~\eqref{041}. Therefore, $(\psf,\hsf)$ forms a coproduct structure on $X^*$. By Theorem~\ref{130}, $\NCSym_q$ admits a Hopf structure.

\subsubsection{Quasisymmetric functions in noncommuting variables}\label{182}

Here we apply the method of restricted actions introduced in Subsection~\ref{056}. Let $S$ be the set of all transpositions $s_{a,b}=(a\,\,b)$ with $a,b\in\Q$ and $a<b$, which forms a generating set for $\sym_\Q$. For each $s=s_{a,b}\in S$, we consider
\[X^*_s=\{u\in X^*\mid |\alph(u)\cap[x_a,x_b]|\leq1\},\quad\text{where}\quad [x_a,x_b]=\{x_r\mid a\leq r\leq b\}.\]Since $\rho_s([x_a,x_b])=[x_a,x_b]$, we have $|\alph(\rho_s(u))\cap[x_a,x_b]|=|\alph(u)\cap[x_a,x_b]|$, which implies $\rho_s(X^*_s)=X^*_s$. Moreover, because $\alph(uv)=\alph(u)\cup\alph(v)$ for all $u,v\in X^*$, the condition $|\alph(uv)\cap[x_a,x_b]|\leq 1$ forces both $u,v\in X^*_s$, ensuring that the complement $X^*\setminus X^*_s$ is a two-sided ideal of $X^*$. Furthermore, $\rho_s$ is multiplicative on $X^*_s$ because the symmetrizing action $\rho$ is multiplicative. Therefore, by Propositions~\ref{054}~and~\ref{055}, the family $\{X^*_s\mid s\in S\}$ induces a bigraded quasi-multiplicative restricted action $\bar{\rho}:\bar{\sym}_\Q\to\sym(X^*)$ and a restricted weight function $\bar{\omega}_q:\bar{\sym}_\Q\times X^*\to\field^\times$ such that $(\tau_q,\bar{\omega}_q)$ is a multiplicative pair. Hence, Condition~\eqref{018} holds for $\bar{\sym}_\Q$.

We set $\NCQSym_q=\funcs{\field}{X^*}^{\bar{\sym}_\Q}$, the \emph{algebra of $q$-quasisymmetric functions in noncommuting variables}. In particular, $\NCQSym:=\NCQSym_1$ is the algebra of quasisymmetric functions in noncommuting variables~\cite{BeZa09,NoThWi10}.

A word $u\in X^*$ is called \emph{packed} if either $u=1$, or $\alph(u)=\{x_1,\ldots,x_m\}$ for some $m\geq1$~\cite{NoTh06}. The \emph{packing} of a word $u$, denoted by $\pack(u)$, is the unique packed word obtained by replacing the $i$th smallest letter occurring in $u$ by $x_i$ for each $i$; for instance, $\pack(x_{\frac{7}{2}}x_{-1}x_{\frac{7}{2}}x_4)=x_2x_1x_2x_3$. Since each $s_{a,b}$ acts nontrivially on $u$ under $\bar{\rho}$ only when $u$ contains exactly one letter from $\{x_a,x_b\}$ and no letters in $(x_a,x_b)$, the restricted action $\bar{\rho}$ preserves the relative order of the letters of $u$. Moreover, any two words sharing the same packing can be transformed into one another by a finite sequence of such order-preserving reindexings. Consequently, two words belong to the same $\bar{\sym}_\Q$-orbit if and only if they have the same packing.

By definition, for any word $u\in X^*$ and any $\sigma\in(\bar{\sym}_\Q)_{(u)}$, we have $\inv(\bar{\rho}_\sigma(u))-\inv(u)=0$ and thus $\bar{\omega}_q(\sigma,u)=q^0=1$. Therefore, every $\bar{\sym}_\Q$-orbit is $\bar{\omega}_q$-admissible. We choose as our system of representatives the set $\bar{T}_q\subset X^*$ of all packed words. Since the restricted action preserves word length, the number of $\bar{\sym}_\Q$-orbits intersecting $(X^*)^d$ is the number of packed words of length $d$, which equals the $d$th Fubini number~\cite{NoTh06}, \cite[\href{https://oeis.org/A000670}{A000670}]{OEIS}. Thus, Condition~\eqref{051} holds, and by Theorem~\ref{057}, the algebra $\NCQSym_q$ admits the monomial basis $\{\bar{m}_u^q\mid u\in\bar{T}_q\}$, such that\[m_t^q=\sum_{u\in T_t}\omega_q(\sigma_u,t)\,\bar{m}_u^q\,,\]where $T_t$ is the set of packed words in the orbit $\sym_\Q \cdot t$, and $\sigma_u\in\sym_\Q$ satisfies $\rho_{\sigma_u}(t)=u$.

Finally, we show that the coproduct structure $(\psf,\hsf)$ on $X^*$ defined in Subsubsection~\ref{123} transfers to $\NCQSym_q$ by verifying the conditions of Proposition~\ref{170}. Let $s=s_{a,b}\in S$ with $a<b$ in $\Q$. By definition, $s\sqcup1$ transposes $a$ with $b$ in the first copy of $\Q$ within $\Q\sqcup\Q$ while fixing all other elements. Because $\iota:\Q\to\Q\sqcup\Q$ is an order-isomorphism, the conjugation $\hsf(s,1)=\iota^{-1}\circ(s\sqcup 1)\circ\iota$ is precisely the transposition in $\sym_\Q$ exchanging $\iota^{-1}(a)$ and $\iota^{-1}(b)$, so $\hsf(s,1)\in S$. Symmetrically, $\hsf(1,s)\in S$. Now, let $u\in X^*$ with $\psf(u)=\qsf(\psi(u))=(u_1,u_2)$, and let $s=s_{a,b}\in S$. In the direct product, the generator $(s,1)$ acts by transposing $a$ and $b$ within the first copy of $\Q$ inside $\Q\sqcup\Q$. Since $\iota:\Q\to\Q\sqcup\Q$ is an order-isomorphism and all elements of the first copy are strictly less than those of the second copy, the closed interval $[\iota^{-1}(a),\iota^{-1}(b)]$ is entirely contained in the preimage of the first copy, and is therefore disjoint from the indices contributing to $u_2$. Because $\psi$ induces an order-isomorphism between $\alph(u)\cap[x_{\iota^{-1}(a)}, x_{\iota^{-1}(b)}]$ and $\alph(u_1)\cap[x_a,x_b]$, we obtain $|\alph(u)\cap[x_{\iota^{-1}(a)},x_{\iota^{-1}(b)}]|=|\alph(u_1)\cap[x_a,x_b]|$, which implies that $u\in M_{\hsf(s,1)}$ if and only if $u_1\in X^*_s$. An analogous argument yields $u\in M_{\hsf(1,s)}$ if and only if $u_2\in X^*_s$.

Let $u=(u_1,u_2,u_3)\in X^*\times X^*\times X^*$, and consider the fibres $A$ and $B$ in~\eqref{175}. As established in Subsubsection~\ref{123}, there is a bijection $f:A\to B$ induced by the order-preserving reindexing between the two partitions of $X$ into three ordered copies $X<Y<Z$. Because $f$ preserves the relative order of the letters of $v$ for every $v\in A$, the words $v$ and $w=f(v)$ have the exact same packed word, meaning $w=\bar{\rho}_\sigma(v)$ for some $\sigma\in\bar{\sym}_\Q$, and $\inv(w)=\inv(v)$, so $\bar{\omega}_q(\sigma,v)=q^{\inv(w)-\inv(v)}=1$. Furthermore, for any $t'\in\bar{T}_q$, if $\psf(v)=(1,t')$, then $v=\psi^{-1}(t')$. Since $\psi^{-1}$ is induced by an order-isomorphism, $v$ preserves the relative order of the letters of $t'$, so $\pack(v)=t'$ and $\psf^{-1}(\{(1,t')\})=\{v\}\subset\bar{\sym}_\Q t'$ with $\bar{\omega}_q(\sigma_v,t')=1$. Symmetrically, $\psf^{-1}(\{(t',1)\})$ is a singleton in $\bar{\sym}_\Q t'$ with restricted weight $1$. Therefore, $(\psf,\hsf)$ is a coproduct structure on $X^*$ with respect to $(\bar{\sym}_\Q,\bar{\rho},\bar{\omega}_q)$. Consequently, by Theorem~\ref{130}, $\NCQSym_q$ is a connected Hopf algebra.

\subsubsection{Symmetric functions in commuting variables}\label{194}

Let $V=[X^*,X^*]$, so that the quotient monoid $X^*/V$ is isomorphic to the free commutative monoid $X^\oplus$. Note that two words in $X^*$ are congruent modulo $V$ if and only if they share the same multiset of letters. In particular, each congruence class $[u]$ is finite, as it contains at most $\ell(u)!$ words. Let $U\subset X^*$ be the system of representatives for the classes of $X^*/V$ formed by words $x_{r_1}\cdots x_{r_k}$ such that $r_1\leq\cdots\leq r_k$ with $k\geq0$, that is, $U=\{\sort(u)\mid u\in X^*\}$, with canonical retraction $\ssf=\sort:X^*\to U$.

Our goal is to determine when the pair $(V,U)$ defines an invariant-retraction pair as in Subsection~\ref{176}. Since the pairs generating $V$ consist of words of the same length, $V$ is immediately bigraded. Moreover, for any $(u,u'),(v,v')\in V$, we have $\inv(u,v)=\inv(u',v')$ because $\sort(u)=\sort(u')$ and $\sort(v)=\sort(v')$. Thus, $\tau_q(u,v)=\tau_q(u',v')$, so $V$ is compatible with the twisting $\tau_q$. In addition, $V$ is stable under the symmetrizing action because if two words share the same multiset of letters, their images under the action also share the same multiset of letters.

For the $q$-weight function $\omega_q$ to be well-defined on the quotient, we must have $\omega_q(\sigma,u)=\omega_q(\sigma,u')$ for all $(u,u')\in V$; however, this condition depends on the parameter $q$. For instance, considering $u=x_1x_2$ and $u'=x_2x_1$, so that $(u,u')\in V$, and letting $\sigma=(1\,\,2)\in\sym_\Q$ be the transposition swapping $1$ and $2$, we obtain\[q=q^{1-0}=q^{\inv(\rho_\sigma(u))-\inv(u)}=\omega_q(\sigma,u)=\omega_q(\sigma,u')=q^{\inv(\rho_\sigma(u'))-\inv(u')}=q^{0-1}=q^{-1},\]which implies $q^2=1$, and thus $q\in\{1,-1\}$. In the remainder of this subsubsection, we assume that $q=\pm1$, which implies $\phi_q(n,m)\in\{\pm1\}$ for all $n,m\in\N$.

To prove that $V$ is compatible with the weight function $\omega_q$, it suffices to show that $\omega_q(\sigma,u)=\omega_q(\sigma,u')$ for any words $u,u'\in X^*$ that differ by a single swap of adjacent letters, as such pairs generate $V$. Let $u=wx_ix_jw'$ and $u'=wx_jx_iw'$ for some $w,w'\in X^*$ and $x_i,x_j\in X$. If $i=j$, the identity holds trivially. If $i\neq j$, swapping the adjacent distinct letters $x_i$ and $x_j$ changes the total number of inversions by one, that is, $\inv(u')=\inv(u)\pm1$, which implies $\inv(u')\equiv\inv(u)+1\pmod 2$. Applying the symmetrizing action yields $\rho_\sigma(u)=\rho_\sigma(w)x_{\sigma(i)}x_{\sigma(j)}\rho_\sigma(w')$ and $\rho_\sigma(u')=\rho_\sigma(w)x_{\sigma(j)}x_{\sigma(i)}\rho_\sigma(w')$. Since $\sigma$ is injective, $\sigma(i)\neq\sigma(j)$, which means that $\rho_\sigma(u)$ and $\rho_\sigma(u')$ also differ from each other by swapping two adjacent distinct letters. Thus, we similarly obtain $\inv(\rho_\sigma(u'))\equiv\inv(\rho_\sigma(u))+1\pmod 2$. Subtracting these congruences yields\[\inv(\rho_\sigma(u'))-\inv(u')\equiv\inv(\rho_\sigma(u))-\inv(u)\pmod2.\]Because $q^2=1$, the value of $q^k$ depends only on the parity of $k$, which implies\[\omega_q(\sigma,u')=q^{\inv(\rho_\sigma(u'))-\inv(u')}=q^{\inv(\rho_\sigma(u))-\inv(u)}=\omega_q(\sigma,u).\]

Finally, to verify compatibility with the coproduct structure, let $u=wx_ix_jw'$ and $u'=wx_jx_iw'$. Recall that $\psf=\qsf\circ\psi$. If the images of $x_i$ and $x_j$ under $\psi$ both belong to $X$ or both to $Y$, the projection under $\qsf$ to that component simply swaps two adjacent letters, leaving its multiset unchanged, while the other component is unaffected. If one letter belongs to $X$ and the other to $Y$, the canonical projection $\qsf$ separates them into distinct components, so $\qsf(\psi(u))=\qsf(\psi(u'))$. In all cases, we obtain $(\sort\times\sort)(\psf(u))=(\sort\times\sort)(\psf(u'))$. Therefore, $(V,U)$ is an invariant-retraction pair.

Together with these ingredients, we set $\Sym_q:=\funcs{\field}{X^\oplus}^{\sym_\Q}$, the \emph{algebra of $(\pm1)$-symmetric functions in commuting variables}. In particular, $\Sym:=\Sym_1$ is the classical algebra of symmetric functions in commuting variables~\cite{Mac99,St99}. Let $T_{q,V}=\{t\in T_q\mid[t]\text{ is }\omega_q^{_V}\text{-admissible}\}$. For any $t\in T_{q,V}$, Lemma~\ref{160} yields\[c_t=|[t]\cap\sym_\Q t|=\prod_{r\geq1}m_r!\,,\]where $m_r$ is the number of distinct letters occurring in $t$ with multiplicity $r$. Notice that $c_t$ depends only on the partition of letter multiplicities of $t$, meaning that $\sym_\Q[u]=\sym_\Q[v]$ in $X^\oplus$ implies $c_u=c_v$. Consequently, since $(V,U)$ is an invariant-retraction pair and $(\psf,\hsf)$ is a coproduct structure on $X^*$, Theorem~\ref{153} and Proposition~\ref{166} imply that the coproduct structure on $X^*$ descends to the free commutative monoid $X^\oplus$ so that $\Sym_q$ is a connected Hopf algebra, isomorphic to the quotient of $\NCSym_q$ by the graded bi-ideal $J$ generated by
\[\{m_u^q\mid u\in T_q\setminus T_{q,V}\}\cup\{m_u^q-\omega_q^{_V}(\sigma,[u])\,m_v^q\mid u,v\in T_{q,V}\text{ and }\sigma\in\sym_\Q\text{ such that }\rho_\sigma^{_V}([u])=[v]\}.\]
When $q=1$, every orbit is admissible, so $T_{1,V}=T_1$, the orbits in $X^\oplus$ are indexed by integer partitions, and the generators of $J$ simplify to the differences $m_u^1-m_v^1$ for all restricted growth words $u,v\in T_1$ having the same partition of letter multiplicities~\cite[Subsection~3.6]{HiNoTh08}. When $q=-1$, a transposition swapping two distinct letters of equal multiplicity $a$ in $[t]$ changes the number of inversions by an integer congruent to $a^2$ modulo $2$, so $[t]$ is $\omega_{-1}^{_V}$-admissible if and only if no odd multiplicity is repeated.

\subsubsection{Quasisymmetric functions in commuting variables}

Following the construction of Subsubsection~\ref{182} for $q\in\{1,-1\}$, we work over the free commutative monoid $X^\oplus=X^*/V$, equipped with the induced twisting $\tau_q^{_V}([u],[v])=\tau_q(\sort(u),\sort(v))$, the degree map $\deg([u])=\ell(u)$, the induced symmetrizing action $\rho^{_V}:\sym_\Q\to\Aut(X^\oplus)$, and the induced weight function $\omega_q^{_V}(\sigma,[u])=\omega_q(\sigma,\sort(u))$. In analogy with Subsubsection~\ref{182}, we apply the method of restricted actions of Subsection~\ref{171} to $X^\oplus$ using the generating set $S=\{s_{a,b}\mid a<b\text{ in }\Q\}$ of all transpositions in $\sym_\Q$ and the family of subsets\[X_s^\oplus=\{u\in X^\oplus\mid|\alph(u)\cap[x_a,x_b]|\leq1\}\quad\text{for each }s=s_{a,b}\in S.\]By the same arguments as in the noncommutative case, each $X_s^\oplus$ is stable and multiplicative under $\rho_s^{_V}$, and its complement $X^\oplus\setminus X_s^\oplus$ is an ideal of $X^\oplus$. Thus, Propositions~\ref{054} and~\ref{055} yield a bigraded quasi-multiplicative restricted action $\bar{\rho}^{_V}:\bar{\sym}_\Q\to\sym(X^\oplus)$ and a restricted weight function $\bar{\omega}_q^{_V}:\bar{\sym}_\Q\times X^\oplus\to\field^\times$ satisfying Condition~\eqref{018}.

The restricted orbit $\bar{\sym}_\Q u$ of a commutative monomial $u=x_{r_1}^{a_1}\cdots x_{r_k}^{a_k}\in X^\oplus$ with $r_1<\cdots<r_k$ and $a_i\geq1$ consists of all monomials $x_{s_1}^{a_1}\cdots x_{s_k}^{a_k}$ sharing the exact same exponent sequence $(a_1,\dots,a_k)$ with $s_1<\cdots<s_k$. Notice that for any $\sigma\in\bar{\sym}_\Q$, the sorted representatives in $U\subset X^*$ of both $u$ and $\bar{\rho}_\sigma^{_V}(u)$ have zero inversions, which implies that the restricted weight function is trivial, that is, $\bar{\omega}_q^{_V}(\sigma,u)=1$ for all $\sigma\in\bar{\sym}_\Q$ and $u\in X^\oplus$. Consequently, every $\bar{\sym}_\Q$-orbit in $X^\oplus$ is $\bar{\omega}_q^{_V}$-admissible, and a canonical system of representatives $\bar{T}^\oplus\subset X^\oplus$ is given by the monomials $x_1^{a_1}\cdots x_k^{a_k}$ with $k\geq0$ and $a_1,\dots,a_k\geq1$, which are naturally indexed by integer compositions $\alpha=(a_1,\dots,a_k)$.

We set $\QSym_q:=\funcs{\field}{X^\oplus}^{\bar{\sym}_\Q}$, the \emph{algebra of $(\pm1)$-quasisymmetric functions in commuting variables}; in particular, $\QSym:=\QSym_1$ is the classical algebra of quasisymmetric functions~\cite{Ges84,LuMyWi13}. By Proposition~\ref{170}, the coproduct structure $(\psf^{_V},\hsf)$ on $X^\oplus$ inherited from $\Sym_q$ transfers to $(\bar{\sym}_\Q,\bar{\rho}^{_V},\bar{\omega}_q^{_V})$, endowing $\QSym_q$ with the structure of a connected Hopf algebra such that the inclusion $\Sym_q\hookrightarrow\QSym_q$ is an injective Hopf algebra homomorphism. Equivalently, $(V,U)$ is an invariant-retraction pair on $X^*$ with respect to $(\bar{\sym}_\Q,\bar{\rho},\bar{\omega}_q)$, so Theorem~\ref{153} yields a surjective Hopf algebra homomorphism $[\ ]_{\bar{\sym}_\Q}:\NCQSym_q\to\QSym_q$, completing the commutative square of connected Hopf algebras:\[\begin{array}{ccc}\NCSym_q&\hookrightarrow&\NCQSym_q\\[4pt]
\downarrow&&\downarrow\\[4pt]\Sym_q&\hookrightarrow&\QSym_q\end{array}\]

\subsubsection{MacMahon symmetric functions}

The constructions of Subsection~\ref{183} allow us to recover symmetric functions in multiple alphabets both in noncommuting and commuting variables~\cite{Ro01,RoRoSt02}. Throughout this subsubsection, we specialize to $q=1$. Fix a positive integer $k\geq1$, and consider $k$ copies of the free monoid $M_1=\cdots=M_k=X^*$, each equipped with the trivial twisting $\tau_1\equiv1$, the word-length grading $\deg(u)=\ell(u)$, the symmetrizing action $\rho:\sym_\Q\to\Aut(X^*)$, and the trivial weight function $\omega_1\equiv1$ from Subsubsection~\ref{123}.

In the noncommutative setting, we consider the $k$-fold free product monoid $(X^*)^{*k}:=X^**\cdots*X^*$, which is naturally isomorphic to the free monoid $(X^{(1)}\sqcup\cdots\sqcup X^{(k)})^*$ generated by $k$ disjoint copies of the alphabet $X$. By Propositions~\ref{067}, \ref{184}, and~\ref{068}, $(X^*)^{*k}$ is equipped with the trivial free twisting $\tau_1^*\equiv1$, the total length degree map $\deg(u_1\cdots u_m)=\ell(u_1)+\cdots+\ell(u_m)$, the free diagonal action $\rho^*:\sym_\Q\to\Aut((X^*)^{*k})$ acting simultaneously on all $k$ copies of $X$, and the trivial free diagonal weight function $\omega_1^*\equiv1$, satisfying Conditions~\eqref{015} to~\eqref{051}. Every $\sym_\Q$-orbit is $\omega_1^*$-admissible, with representatives indexed by $k$-colored set partitions~\cite{NoTh10,Ro01,RoRoSt02,Strot26}. Duplicating each copy of $X$ via the order-isomorphism $\iota:\Q\to\Q\sqcup\Q$ endows $\MNCSym^{(k)}:=\funcs{\field}{(X^*)^{*k}}^{\sym_\Q}$ with the structure of a connected Hopf algebra, the algebra of \emph{MacMahon symmetric functions in noncommuting variables}.

In the commutative setting, we consider instead the $k$-fold direct product of free commutative monoids $(X^\oplus)^k:=X^\oplus\times\cdots\times X^\oplus$, which is a quotient of $(X^*)^{*k}$ by the commutator congruence. By Proposition~\ref{076}, $(X^\oplus)^k$ is equipped with the trivial twisting, the degree map $\deg(u_1,\dots,u_k)=\ell(u_1)+\cdots+\ell(u_k)$, the diagonal action $\rho^\times:\sym_\Q\to\Aut((X^\oplus)^k)$ given by $\rho_\sigma^\times(u_1,\dots,u_k)=(\rho_\sigma^{_V}(u_1),\dots,\rho_\sigma^{_V}(u_k))$, and the trivial diagonal weight function $\omega_1^\times\equiv1$. The resulting algebra of invariants $\MSym^{(k)}:=\funcs{\field}{(X^\oplus)^k}^{\sym_\Q}$ is the connected Hopf algebra of \emph{MacMahon symmetric functions in commuting variables}, whose monomial basis is indexed by vector partitions, that is, multisets of non-zero vectors in $\N^k$~\cite{Ro01,RoRoSt02}.

\subsection{Symmetric functions in superspace}\label{192}

Let $\Theta=(\theta_r\mid r\in\Q)$ be another infinite alphabet linearly ordered by the natural order of its indices, and let $\phi^s:\Z_2\times\Z_2\to\field^\times$ be defined by $\phi^s(x,y)=(-1)^{xy}$. Note that $\phi^s$ is a bicharacter because for all $x,y,z\in\Z_2$,\[\phi^s(x,y)=\phi^s(y,x)\quad\text{and}\quad\phi^s(x+y,z)=(-1)^{(x+y)z}=(-1)^{xz+yz}=(-1)^{xz}(-1)^{yz}=\phi^s(x,z)\,\phi^s(y,z).\]

The \emph{fermionic monoid} is the monoid $\Theta^\wedge$ presented by generators $\Theta$ subject to the relations $\theta_r\theta_{r'}=\theta_{r'}\theta_r$ and $\theta_r^2=1$ for all $r,r'\in\Q$, which is in fact an abelian group. Note that any element in $\Theta^\wedge$ can be uniquely written as a word $\theta_{r_1}\cdots\theta_{r_k}$ such that $k\geq0$ and $r_1<\cdots<r_k$. Furthermore, for $u,v\in\Theta^\wedge$, we have $u=v$ if and only if $\alph(u)=\alph(v)$. Thus, the elements of $\Theta^\wedge$ can be naturally identified with finite subsets of $\Q$ by writing $\theta_I=\theta_{r_1}\cdots\theta_{r_k}\in\Theta^\wedge$ and $\theta_I^*=\theta_{r_1}\cdots\theta_{r_k}\in\Theta^*$ with $I=\alph(\theta_I)=\{r_1<\cdots<r_k\}$, where $\theta_\emptyset=1$ and $\theta_\emptyset^*=1$. In the group $\Theta^\wedge$, we have $\theta_I\theta_J=\theta_{I\bigtriangleup J}$, which reduces to $\theta_I\theta_J=\theta_{I\cup J}$ whenever $I\cap J=\emptyset$. Given $I,J$ finite subsets of $\Q$, a \emph{cross-inversion} of $(\theta_I,\theta_J)$ is an inversion of the word $\theta_I^*\theta_J^*$, that is, a pair $(i,j)\in I\times J$ such that $i>j$. The number of such cross-inversions is denoted by $\inv(\theta_I,\theta_J)$. Given a finite subset $K\subset\Q$, we have $\sort(u)=\theta_K^*$ for any representative $u\in\Theta^*$ of $\theta_K$ of minimal length, whence\begin{gather}\label{177}\inv(\theta_I,\theta_J)=\inv(\theta_I^*,\theta_J^*)=\inv(\theta_I^*\theta_J^*).\end{gather}

Let $\tau_\wedge:\Theta^\wedge\times\Theta^\wedge\to\field$ be the map defined for all finite subsets $I,J\subset\Q$ by\[\tau_\wedge(\theta_I,\theta_J)=\begin{cases}(-1)^{\inv(\theta_I,\theta_J)}&\text{if }I\cap J=\emptyset,\\0&\text{otherwise}.\end{cases}\]By~\eqref{177}, if $I,J$ are disjoint finite subsets of $\Q$, we have\begin{gather}\label{178}\tau_\wedge(\theta_I,\theta_J)=(-1)^{\inv(\theta_I,\theta_J)}=(-1)^{\inv(\theta_I^*,\theta_J^*)}=\tau_{-1}(\theta_I^*,\theta_J^*).\end{gather}We claim that $\tau_\wedge$ defines a twisting on $\Theta^\wedge$; indeed, let $I,J,K$ be finite subsets of $\Q$. If $I,J,K$ are pairwise disjoint, we have $\theta_I\theta_J=\theta_{I\cup J}$ and $\theta_J\theta_K=\theta_{J\cup K}$, so by~\eqref{178} we obtain $\tau_\wedge(\theta_I,\theta_J)\tau_\wedge(\theta_I\theta_J,\theta_K)=\tau_\wedge(\theta_I,\theta_J\theta_K)\tau_\wedge(\theta_J,\theta_K)$ because $\tau_{-1}$ is a twisting on $\Theta^*$. If $I,J,K$ are not pairwise disjoint, at least one factor in each product vanishes, so both sides are equal to zero. Moreover, since every finite subset $I\subset\Q$ is disjoint from $\emptyset$, by~\eqref{178} we obtain $\tau_\wedge(\theta_I,1)=1=\tau_\wedge(1,\theta_I)$, and therefore $\tau_\wedge$ is a twisting on $\Theta^\wedge$. With respect to this twisting, the monoid $\Theta^\wedge$ is locally finite, since for any finite subset $K\subset\Q$ there are only finitely many pairs of finite subsets $(I,J)$ of $\Q$ such that $\tau_\wedge(\theta_I,\theta_J)\neq 0$ and $\theta_I\theta_J=\theta_K$, namely, all disjoint partitions $I\cup J=K$, of which there are $2^{|K|}$. Notice that if $I,J$ are disjoint finite subsets of $\Q$ with $|I|=h$ and $|J|=k$, then $\tau_\wedge(\theta_I,\theta_J)=\sgn(\pi_{I,J})$, where $\pi_{I,J}$ is the unique permutation of $[h+k]$ that sorts $\theta_I^*\theta_J^*$ into $\theta_{I\cup J}^*$\,; that is, if $I=\{r_1<\cdots<r_h\}$, $J=\{r'_1<\cdots<r'_k\}$, and $I\cup J=\{s_1<\cdots<s_{h+k}\}$, then $s_{\pi_{I,J}(m)}=r_m$ for $m\in[h]$, and $s_{\pi_{I,J}(m)}=r'_{m-h}$ for $h+1\le m\le h+k$.

For any finite subset $I\subset\Q$, we set $|\theta_I|:=\bar{|I|}\in\Z_2$, the canonical projection of $|I|$ modulo $2$. It is immediate that the map $\theta_I\mapsto|\theta_I|$ defines a monoid homomorphism from $\Theta^\wedge$ to $\Z_2$ because $|I\bigtriangleup J|\equiv|I|+|J|\pmod2$. Thus, the monoid $\Theta^\wedge$ is $\Z_2$-graded for the trivial twisting, with fibres $\ubar{\Theta}^\wedge_{\bar{0}}$ and $\ubar{\Theta}^\wedge_{\bar{1}}$. Furthermore, the map $\deg:\Theta^\wedge\to\N$ given by $\deg(\theta_I)=|I|$ defines a grading on $\Theta^\wedge$ with respect to $\tau_\wedge$, since $\tau_\wedge(\theta_I,\theta_J)=0$ whenever $I\cap J\neq\emptyset$, while $\deg(\theta_I\theta_J)=|I\cup J|=|I|+|J|=\deg(\theta_I)+\deg(\theta_J)$ whenever $I\cap J=\emptyset$, with $\deg(\theta_I)=0$ if and only if $I=\emptyset$. Note that each homogeneous component $(\Theta^\wedge)^d=\{\theta_I\in\Theta^\wedge\mid\deg(\theta_I)=d\}$ intersects exactly one fibre of the $\Z_2$-grading; in fact, since $|\theta_I|=\bar{d}$ for all $\theta_I\in(\Theta^\wedge)^d$, the component $(\Theta^\wedge)^d$ is contained in the single fibre $\ubar{\Theta}^\wedge_{\bar{d}}$. Consequently, $(\Theta^\wedge)^d_x=(\Theta^\wedge)^d\cap\ubar{\Theta}^\wedge_x$ is empty for $x\neq\bar{d}$, and equals $(\Theta^\wedge)^d$ for $x=\bar{d}$.

We now consider the symmetrizing action $\rho:\sym_\Q\to\Aut(\Theta^\wedge)$ on the monoid $\Theta^\wedge$, defined by\[\rho_\sigma(\theta_I)=\theta_{\sigma(I)}=\theta_{\sigma(r_1)}\cdots\theta_{\sigma(r_k)},\quad\text{for all}\quad\sigma\in\sym_\Q\quad\text{and}\quad I=\{r_1<\cdots<r_k\}\subset\Q.\]Since permutations preserve cardinalities, $|\sigma(I)|=|I|$, the symmetrizing action is immediately bigraded. Furthermore, let $\omega_s:\sym_\Q\times\Theta^\wedge\to\field^\times$ be the map defined by\[\omega_s(\sigma,\theta_I)=(-1)^{\inv(\rho_\sigma(\theta_I^*))}\]Notice that $\omega_s(\sigma,\theta_I)=\sgn(\pi_{\sigma,I})$, where $\pi_{\sigma,I}$ is the unique permutation of $[k]$ that sorts $\rho_\sigma(\theta_I^*)$ into $\theta_{\sigma(I)}^*$\,; that is, if $\sigma(I)=\{s_1<\cdots<s_k\}$, then $\sigma(r_i)=s_{\pi_{\sigma,I}(i)}$ for all $i\in[k]$. We claim that $\omega_s$ is a weight function for the action $\rho$, which makes $(\tau_\wedge,\omega_s)$ a multiplicative pair. Indeed, for any $\sigma,\varsigma\in\sym_\Q$, the permutation $\pi_{\varsigma,I}$ orders the indices obtained after applying $\varsigma$, while $\pi_{\sigma,\varsigma(I)}$ orders the indices obtained after subsequently applying $\sigma$, which implies $\pi_{\sigma\varsigma,I}=\pi_{\sigma,\varsigma(I)}\pi_{\varsigma,I}$. Hence\[\omega_s(\sigma\varsigma,\theta_I)=\sgn(\pi_{\sigma\varsigma,I})=\sgn(\pi_{\sigma,\varsigma(I)}\pi_{\varsigma,I})=\sgn(\pi_{\sigma,\varsigma(I)})\sgn(\pi_{\varsigma,I})=\omega_s(\sigma,\rho_\varsigma(\theta_I))\,\omega_s(\varsigma,\theta_I).\]By definition, if $I\cap J\neq\emptyset$, we have $\tau_\wedge(\theta_I,\theta_J)=0=\tau_\wedge(\rho_\sigma(\theta_I),\rho_\sigma(\theta_J))$, so that~\eqref{028} holds trivially. If $I\cap J=\emptyset$, let $\pi$ be the unique permutation of $[h+k]$ sorting $\rho_\sigma(\theta_I^*)\rho_\sigma(\theta_J^*)$ into $\theta_{\sigma(I\cup J)}^*$. The left-hand side of~\eqref{028} factors this sorting by first ordering $\theta_I^*\theta_J^*$ via $\pi_{I,J}$ and then sorting the image under $\sigma$ via $\pi_{\sigma,I\cup J}$. Similarly, the right-hand side factors the same sorting by first ordering each block $\rho_\sigma(\theta_I^*)$ and $\rho_\sigma(\theta_J^*)$ independently via $\pi_{\sigma,I}$ and $\pi_{\sigma,J}$, and then sorting their product via $\pi_{\sigma(I),\sigma(J)}$; specifically,\[\begin{array}{rclll}
\tau_\wedge(\theta_I,\theta_J)\,\omega_s(\sigma,\theta_I\theta_J)&=&\sgn(\pi_{I,J})\,\sgn(\pi_{\sigma,I\cup J})\\&=&\sgn(\pi)\\&=&\sgn(\pi_{\sigma(I),\sigma(J)})\,\sgn(\pi_{\sigma,I})\,\sgn(\pi_{\sigma,J})&=&\tau_\wedge(\rho_\sigma(\theta_I),\rho_\sigma(\theta_J))\,\omega_s(\sigma,\theta_I)\,\omega_s(\sigma,\theta_J).\end{array}\]

\begin{rem}
Let $I$ be a finite subset of $\Q$. Notice that the stabilizer of $\theta_I$ is $G_{(\theta_I)}=\{\sigma\in\sym_\Q\mid\sigma(I)=I\}$. If $|I|\le 1$, then $\inv(\rho_\sigma(\theta_I^*))=0$, so $\omega_s(\sigma,\theta_I)=(-1)^0=1$ for all $\sigma\in\sym_\Q$. If $|I|\ge 2$, taking distinct elements $a<b$ in $I$, the transposition $s_{a,b}=(a\;b)$ belongs to $G_{(\theta_I)}$ and satisfies $\omega_s(s_{a,b},\theta_I)=(-1)^1=-1\neq 1$. Consequently, the only $\omega_s$-admissible orbits in $\Theta^\wedge$ are those of degree $0$ and $1$. Since $\sym_\Q$ acts transitively on singletons, $\{1,\theta_{\{1\}}\}$ is a system of representatives for these $\omega_s$-admissible orbits. Therefore, Conditions~\eqref{015} to~\eqref{051} hold.
\end{rem}\medskip

Let $\phi:(\Z_2\times\N)\times(\Z_2\times\N)\to\field^\times$ be defined by $\phi((x,m),(y,n))=(-1)^{xy}\,q^{mn}$ for all $x,y\in\Z_2$ and $m,n\in\N$. Notice that $\phi$ defines a bicharacter; indeed, for all $(x,m),(x',m'),(y,n)\in\Z_2\times\N$, we have
\begin{gather*}
\phi((x+x',m+m'),(y,n))=(-1)^{(x+x')y}\,q^{(m+m')n}=(-1)^{xy}(-1)^{x'y}\,q^{mn}q^{m'n}=\phi((x,m),(y,n))\,\phi((x',m'),(y,n)),\\
\phi((x,m),(y+y',n+n'))=(-1)^{x(y+y')}\,q^{m(n+n')}=(-1)^{xy+xy'}\,q^{mn+mn'}=\phi((x,m),(y,n))\,\phi((x,m),(y',n')).
\end{gather*}
Since $\phi_q(m,n)=q^{mn}=\phi((\bar{0},m),(\bar{0},n))$ for all $m,n\in\N$, the constructions in Subsection~\ref{179} can be carried out using $\phi$ instead of $\phi_q$. Indeed, $X^*$ and $X^\oplus$ can be regarded as $(\Z_2\times\N)$-graded monoids by setting $|u|=(\bar{0},\ell(u))$ for every monomial $u$. Consequently,\[\ubar{X}_{(\bar{0},d)}^*=\ubar{X}_d^*\,,\qquad \ubar{X}_{(\bar{1},d)}^*=\emptyset\,,\qquad\ubar{X}_{(\bar{0},d)}^\oplus=\ubar{X}_d^\oplus\,,\quad\text{and}\quad\ubar{X}_{(\bar{1},d)}^\oplus=\emptyset\quad\text{for all}\quad d\in\N.\]Similarly, since $\phi^s(x,y)=(-1)^{xy}=\phi((x,0),(y,0))$ for all $x,y\in\Z_2$, the fermionic monoid $\Theta^\wedge$ can be regarded as a $(\Z_2\times\N)$-graded monoid by setting $|\theta_I|=(\bar{|I|},0)$ for all $\theta_I\in\Theta^\wedge$. Consequently,\[\ubar{\Theta}_{(x,0)}^\wedge=\ubar{\Theta}_x^\wedge\quad\text{and}\quad\ubar{\Theta}_{(x,d)}^\wedge=\emptyset\quad\text{for all}\quad x\in\Z_2\quad\text{and}\quad d\geq1.
\]From now on, we work with this unified bicharacter $\phi$, restricting to $\phi_q$ or $\phi^s$ whenever appropriate.

\subsubsection{Symmetric functions in noncommuting variables in superspace}\label{185}

Let $X^\theta=\Theta^\wedge\times X^*$ be the direct product monoid of the fermionic monoid $\Theta^\wedge$ and the free monoid $X^*$, equipped with the product twisting defined from $\tau_\wedge$ and $\tau_q$ as in Proposition~\ref{075}\eqref{009}; specifically, for $(\theta_I,u),(\theta_J,v)\in X^\theta$,\begin{gather}\label{180}
\tau^\times((\theta_I,u),(\theta_J,v))=\tau_\wedge(\theta_I,\theta_J)\,\tau_q(u,v)=\begin{cases}(-1)^{\inv(\theta_I^*,\theta_J^*)}\,q^{\inv(u,v)}&\text{if }I\cap J=\emptyset,\\[4pt]0&\text{otherwise}.\end{cases}
\end{gather}Notice that the twisting $\tau^\phi$ on $X^\theta$ defined from $\phi$ as in Proposition~\ref{083} is trivial; indeed,\[\tau^\phi((\theta_I,u),(\theta_J,v))=\phi(|u|,|\theta_J|)=\phi((\bar{0},\ell(u)),(\bar{|J|},0))=\phi^s(\bar{0},\bar{|J|})\,\phi_q(\ell(u),0)=(-1)^{\bar{0}\bar{|J|}}q^{\ell(u)0}=1.\]Consequently, the combined twisting $\tau_\phi=\tau^\times\tau^\phi$ coincides with the product twisting $\tau^\times$ in~\eqref{180}.

Since $\Theta^\wedge$ and $X^*$ are locally finite with respect to $\tau_\wedge$ and $\tau_q$, respectively, Proposition~\ref{092}\eqref{093} ensures that $X^\theta$ is locally finite with respect to $\tau^\times$.

As in Subsection~\ref{120}, we endow $X^\theta$ with the $(\Z_2\times\N)$-grading defined for all $(\theta_I,u)\in X^\theta$ by $|(\theta_I,u)|=|\theta_I|+|u|=(\bar{|I|},0)+(\bar{0},\ell(u))=(\bar{|I|},\ell(u))$, whose fibres with respect to the trivial twisting are given by\[\ubar{X}^\theta_{(x,d)}=\{(\theta_I,u)\in X^\theta\mid(\bar{|I|},\ell(u))=(x,d)\}=\ubar{\Theta}^\wedge_x\times\ubar{X}^*_d\quad\text{for all}\quad(x,d)\in\Z_2\times\N.\]Furthermore, the monoid $X^\theta$ inherits a grading with respect to $\tau^\times$ via the degree map $\deg:X^\theta\to\N$ defined by $\deg(\theta_I,u)=\deg(\theta_I)+\deg(u)=|I|+\ell(u)$, which satisfies $\deg(\theta_I,u)=0$ if and only if $(\theta_I,u)=(1,1)$. By~\eqref{091}, for each $d\in\N$, the corresponding homogeneous component is given by the finite disjoint union\[(X^\theta)^d=\bigsqcup_{h+k=d} (\Theta^\wedge)^h\times(X^*)^k.\]Observe that for any $(\theta_I,u)\in (X^\theta)^d$, the length $\ell(u)\in[d]_0$ uniquely determines the $(\Z_2\times\N)$-degree $(\bar{|I|},\ell(u))$, since $|I|=d-\ell(u)$. Consequently, each $(X^\theta)^d$ intersects at most $d+1$ fibres of the $(\Z_2\times\N)$-grading; specifically,\[(X^\theta)^d_{(x,n)}=(X^\theta)^d\cap\ubar{X}^\theta_{(x,n)}=\begin{cases}(\Theta^\wedge)^{d-n}\times(X^*)^n&\text{if }n\in[d]_0\text{ and }x=\bar{d-n},\\\emptyset&\text{otherwise}.\end{cases}\]By identifying $(\theta_r,1)$ and $(1,x_r)$ with $\theta_r$ and $x_r$, respectively, any simple monomial $(\theta_I,u)\in X^\theta$ can be uniquely represented as $\theta_{i_1}\cdots\theta_{i_m}x_{j_1}\cdots x_{j_n}$, where $I=\{i_1<\cdots<i_m\}\subset\Q$ and $u=x_{j_1}\cdots x_{j_n}\in X^*$. Notice that under the twisted convolution product in $\funcs{\field}{X^\theta}$, for any $r,s\in\Q$ with $r<s$, we have\[\theta_s\cdot\theta_r=\tau_\wedge(\theta_s,\theta_r)\,\theta_s\theta_r=-\theta_r\theta_s=-\theta_r\cdot\theta_s\quad\text{and}\quad\theta_r\cdot\theta_r=\tau_\wedge(\theta_r,\theta_r)\,\theta_r^2=0\cdot1=0.\]Furthermore, since $X^\theta$ is a direct product monoid, the elements of $\Theta$ and $X$ commute with each other; specifically, $x_r\cdot\theta_s=(\theta_s,x_r)=\theta_s\cdot x_r$ for all $r,s\in\Q$. Hence, the twisted algebra $\funcs{\field}{X^\theta}$ is naturally identified with the algebra of formal power series of bounded degree in noncommuting variables $x_r$ and anticommuting variables $\theta_r$.

By Proposition~\ref{076}, the monoid $X^\theta$ can be endowed with the diagonal action $\rho^\times:\sym_\Q\to\Aut(X^\theta)$ and the diagonal weight function $\omega^\times:\sym_\Q\times X^\theta\to\field^\times$, defined for all $\sigma\in\sym_\Q$ and $(\theta_I,u)\in X^\theta$ by\[\rho^\times_\sigma(\theta_I,u)=(\rho_\sigma(\theta_I),\rho_\sigma(u))\quad\text{and}\quad\omega^\times(\sigma,(\theta_I,u))=\omega_s(\sigma,\theta_I)\,\omega_q(\sigma,u).\]Moreover, $(\tau^\times,\omega^\times)$ is a multiplicative pair for the multiplicative bigraded action $\rho^\times$, which ensures that Condition~\eqref{018} holds. We call the twisted algebra of invariants $\sNCSym_q:=\funcs{\field}{X^\theta}^{\sym_\Q}$ the \emph{algebra of $q$-symmetric functions in noncommuting variables in superspace}.

We now determine the $\omega^\times$-admissible $\sym_\Q$-orbits in $X^\theta$. Let $(\theta_I,u)\in X^\theta$, where $I\subset\Q$ is finite and $u\in X^*$. Note that any permutation $\sigma\in(\sym_\Q)_{(\theta_I,u)}$ must fix every index in $\ind(u):=\{r\in\Q\mid x_r\in\alph(u)\}$. Consequently, $\sigma$ can only act nontrivially by permuting the fermionic indices that do not occur in $u$, that is, the elements of $I\setminus\ind(u)$. Since $\rho_\sigma$ fixes $\alph(u)$, we have $\omega_q(\sigma,u)=q^{\inv(u)-\inv(u)}=1$, and so the diagonal weight simplifies to\[\omega^\times(\sigma,(\theta_I,u))=\omega_s(\sigma,\theta_I)\,\omega_q(\sigma,u)=\omega_s(\sigma,\theta_I)=(-1)^{\inv(\rho_\sigma(\theta_I^*))}.\]If $|I\setminus\ind(u)|\geq2$, there exist distinct indices $a<b$ in $I\setminus\ind(u)$ such that the transposition $s_{a,b}$ fixes every element in $\ind(u)$ and preserves $I$, so $s_{a,b}\in(\sym_\Q)_{(\theta_I,u)}$. However, $\omega^\times(s_{a,b},(\theta_I,u))=\omega_s(s_{a,b},\theta_I)=-1\neq1$, which implies that the orbit $\sym_\Q(\theta_I,u)$ is not $\omega^\times$-admissible. Conversely, if $|I\setminus\ind(u)|\leq1$, the only permutation in $(\sym_\Q)_{(\theta_I,u)}$ acting on $I$ is the identity, which implies $\omega^\times(\sigma,(\theta_I,u))=1$ for all $\sigma\in(\sym_\Q)_{(\theta_I,u)}$. Therefore, an orbit $\sym_\Q(\theta_I,u)$ is $\omega^\times$-admissible if and only if $|I\setminus\ind(u)|\leq1$.

In analogy with the classical case, the admissible orbits can be canonicalized via restricted growth words. We define the set of \emph{super restricted growth pairs} $T^\theta \subset X^\theta$ as the collection of pairs $(\theta_I,u)$ such that $u\in T_q$ is a restricted growth word on the alphabet $\{x_1,\ldots,x_k\}$ for some $k\geq0$, and $I\subseteq[k]_0$ with $|I\setminus[k]|\leq1$, where $0\in I$ indicates the presence of a unique fermionic index absent from $u$. Every $\omega^\times$-admissible orbit contains a unique representative in $T^\theta$. Under this description, the elements of $T^\theta$ are in one-to-one correspondence with set superpartitions~\cite{ArGoMa25}. Since for each degree $d\in\N$ there are only finitely many such pairs, Condition~\eqref{051} holds. Consequently, by Theorem~\ref{050}, the algebra $\sNCSym_q$ admits the monomial linear basis $\{m_{(\theta_I,u)}^\times\mid(\theta_I,u)\in T^\theta\}$.

To equip $\sNCSym_q$ with a Hopf structure, we apply the duplication construction of Subsubsection~\ref{035} using the order-isomorphism $\iota:\Q\to\Q\sqcup\Q$ and the injective group homomorphism $\hsf:\sym_\Q\times\sym_\Q\to\sym_\Q$ defined in Subsubsection~\ref{123}. By definition, the monoid $X^\theta$ is presented by the generators $x_r$ and $\theta_s$ with $r,s\in\Q$, subject to the relations $\theta_r^2=1$, $x_r\theta_s=\theta_sx_r$, and $\theta_r\theta_s=\theta_s\theta_r$ for all $r,s\in\Q$. In the free product $X^\theta*X^\theta$, we denote the generators of the first copy by $x_r,\theta_r$ and those of the second copy by $y_r,\vartheta_r$. We define $R^\theta$ as the monoid congruence on $X^\theta*X^\theta$ generated by the pairs $(\theta_ry_s,y_s\theta_r)$, $(\vartheta_rx_s,x_s\vartheta_r)$, and $(\theta_r\vartheta_s,\vartheta_s\theta_r)$ for all $r,s\in\Q$. Denoting the two ordered copies in the codomain of $\iota$ by $\Q_1\sqcup\Q_2$, the assignment\[x_r\mapsto\begin{cases}x_{\iota(r)}&\text{if }\iota(r)\in\Q_1,\\y_{\iota(r)}&\text{if }\iota(r)\in\Q_2,\end{cases}\qquad\text{and}\qquad\theta_r\mapsto\begin{cases}\theta_{\iota(r)}&\text{if }\iota(r)\in\Q_1,\\\vartheta_{\iota(r)}&\text{if }\iota(r)\in\Q_2\end{cases}\]preserves the defining relations of $X^\theta$, thereby extending to a unique monoid homomorphism $\psi_\theta:X^\theta\to(X^\theta*X^\theta)/R^\theta$. Since $\iota$ is a bijection, the assignment $x_i,\theta_i\mapsto x_{\iota^{-1}(i)},\theta_{\iota^{-1}(i)}$ for $i\in\Q_1$ and $y_j,\vartheta_j\mapsto x_{\iota^{-1}(j)},\theta_{\iota^{-1}(j)}$ for $j\in\Q_2$ preserves the relations defining $(X^\theta*X^\theta)/R^\theta$, defining an inverse homomorphism. Thus, $\psi_\theta$ is a monoid isomorphism induced by the order-isomorphism $\iota$. Composing this isomorphism with the canonical projection $\qsf_\theta:(X^\theta*X^\theta)/R^\theta\to X^\theta\times X^\theta$ defines the monoid homomorphism $\psf_\theta:=\qsf_\theta\circ\psi_\theta$. Specifically, identifying $\Q_1$ and $\Q_2$ with $\Q$ via their canonical order-isomorphisms $\iota_1:\Q\to\Q_1$ and $\iota_2:\Q\to\Q_2$, for each $(\theta_I,u)\in X^\theta$ we have $\psf_\theta(\theta_I,u)=((\theta_{I_1},u_1),(\theta_{I_2},u_2))$, where $I_1=\iota_1^{-1}(\iota(I)\cap\Q_1)$ and $I_2=\iota_2^{-1}(\iota(I)\cap\Q_2)$ so that $I=\iota^{-1}(\iota_1(I_1)\sqcup\iota_2(I_2))$, and $u_1,u_2\in X^*$ are the subwords of $u$ formed by the letters mapped by $\iota$ to $\Q_1$ and $\Q_2$. Note that the fibres of $\psf_\theta$ are finite; indeed, for any $((\theta_{I_1},u_1),(\theta_{I_2},u_2))\in X^\theta\times X^\theta$, any element in its preimage under $\psf_\theta$ is of the form $(\theta_{\iota^{-1}(\iota_1(I_1)\sqcup\iota_2(I_2))},u)\in X^\theta$, where the word $\psi_\theta(1,u)$ in the generators $x$ and $y$ is a $(u_1,u_2)$-shuffle. Since the set of $(u_1,u_2)$-shuffles is finite as noted in Subsubsection~\ref{123}, $\psf_\theta$ has finite fibres. Moreover, since $\ell(u)=\ell(u_1)+\ell(u_2)$ and $|I|=|I_1|+|I_2|$, the map $\psf_\theta$ preserves both the $(\Z_2\times\N)$-degree and the $\N$-degree; thus, $\psf_\theta$ is bigraded.

We now verify that $\psf_\theta$ is twisting-preserving with respect to $\tau^\times$ and the bicharacter $\phi$; that is, for all $a,b\in X^\theta$ with $\psf_\theta(a)=(a_1,a_2)$ and $\psf_\theta(b)=(b_1,b_2)$, we must show that $\tau^\times(a,b)=\tau^\times(a_1,b_1)\,\tau^\times(a_2,b_2)\,\phi(|a_2|,|b_1|)$. Let $a=(\theta_I,u)$ and $b=(\theta_J,v)$ in $X^\theta$, where $a_i=(\theta_{I_i},u_i)$ and $b_i=(\theta_{J_i},v_i)$ for $i\in\{1,2\}$. By~\eqref{180}, if $I\cap J\neq\emptyset$, then $I_1\cap J_1\ne\emptyset$ or $I_2\cap J_2\ne\emptyset$, which implies $\tau^\times(a,b)=0$ and $\tau^\times(a_1,b_1)\,\tau^\times(a_2,b_2)=0$. Assume now that $I\cap J=\emptyset$. Since all indices in $\Q_1$ strictly precede those in $\Q_2$, the cross-inversions decompose as\[\inv(\theta_I^*,\theta_J^*)=\inv(\theta_{I_1}^*,\theta_{J_1}^*)+\inv(\theta_{I_2}^*,\theta_{J_2}^*)+|I_2||J_1|\quad\text{and}\quad\inv(u,v)=\inv(u_1,v_1)+\inv(u_2,v_2)+\ell(u_2)\ell(v_1).\]
Consequently,\[\begin{array}{rclll}\tau^\times(a,b)&=&\tau_\wedge(\theta_I,\theta_J)\,\tau_q(u,v)&=&(-1)^{\inv(\theta_{I_1}^*,\theta_{J_1}^*)+\inv(\theta_{I_2}^*,\theta_{J_2}^*)+|I_2||J_1|}\,q^{\,\inv(u_1,v_1)+\inv(u_2,v_2)+\ell(u_2)\ell(v_1)}\\[0.2cm]&&&=&\tau^\times(a_1,b_1)\,\tau^\times(a_2,b_2)\,(-1)^{|I_2||J_1|}\,q^{\ell(u_2)\ell(v_1)}.\end{array}\]On the other hand, since $|a_2|=(\bar{|I_2|},\ell(u_2))$ and $|b_1|=(\bar{|J_1|},\ell(v_1))$, evaluating the bicharacter $\phi$ yields\[\phi(|a_2|,|b_1|)=(-1)^{\bar{|I_2|}\,\bar{|J_1|}}\,q^{\ell(u_2)\ell(v_1)}=(-1)^{|I_2||J_1|}\,q^{\ell(u_2)\ell(v_1)}.\]Therefore, $\tau^\times(a,b)=\tau^\times(a_1,b_1)\,\tau^\times(a_2,b_2)\,\phi(|a_2|,|b_1|)$, which proves that $\psf_\theta$ is twisting-preserving.

Furthermore, as established in Subsubsection~\ref{123}, the permutations in $\hsf(\sigma,\varsigma)$ act independently on each copy and preserve the relative order between the two copies. More precisely, for any $a=(\theta_I,u)\in X^\theta$, the conjugate action of $\hsf(\sigma,\varsigma)=\iota^{-1}\circ(\sigma\sqcup\varsigma)\circ\iota$ maps the indices of $I_1$ and letters of $u_1$ via $\sigma$, and those of $I_2$ and $u_2$ via $\varsigma$. Hence,\[\psf_\theta(\rho_{\hsf(\sigma,\varsigma)}^\times(a))=(\rho_\sigma^\times(a_1),\rho_\varsigma^\times(a_2))=(\rho^\times)_{(\sigma,\varsigma)}^\otimes(\psf_\theta(a)),\]and for $\psf_\theta(a)=(a_1,a_2)$, the diagonal weight factors as $\omega^\times(\hsf(\sigma,\varsigma),a)=\omega^\times(\sigma,a_1)\,\omega^\times(\varsigma,a_2)=(\omega^\times)^\otimes((\sigma,\varsigma),\psf_\theta(a))$. Thus, the compatibility conditions in~\eqref{125} hold. Finally, the verification of coassociativity~\eqref{108} and counitality~\eqref{122} proceeds identically to that in Subsubsection~\ref{123}: the compositions $(p_\theta\times\id)\circ p_\theta$ and $(\id\times p_\theta)\circ p_\theta$ correspond to partitioning the indices into three consecutive ordered copies of $\Q$, yielding an order-preserving bijection between their fibres that preserves orbits and diagonal weights. Therefore, $(p_\theta,h)$ is a coproduct structure on $X^\theta$, and Condition~\eqref{127} holds. By Theorem~\ref{130}, the algebra $\sNCSym_q$ is a connected Hopf algebra endowed with the coproduct $\Delta=\Phi_2\circ\nabla$.

\subsubsection{Quasisymmetric functions in noncommuting variables in superspace}\label{187}

Here we apply the method of restricted actions from Subsection~\ref{171} to the monoid $X^\theta$ constructed in Subsubsection~\ref{185}. As in Subsubsection~\ref{182}, let $S$ be the set of all transpositions $s_{a,b}$ with $a,b\in\Q$ and $a<b$. For each $s=s_{a,b}\in S$, we consider the subset\[X_s^\theta=\{(\theta_I,u)\in X^\theta\mid|\ind(\theta_I,u)\cap[a,b]|\leq1\},\quad\text{where}\quad\ind(\theta_I,u)=I\cup\ind(u)\quad\text{and}\quad[a,b]=\{r\in\Q\mid a\leq r\leq b\}.\]Since $s([a,b])=[a,b]$ and $\ind(\rho_s^\times(\theta_I,u))=s(\ind(\theta_I,u))$, we have $|\ind(\rho_s^\times(\theta_I,u))\cap[a,b]|=|\ind(\theta_I,u)\cap[a,b]|$, which implies $\rho_s^\times(X_s^\theta)=X_s^\theta$. Now let $z=(\theta_I,u)$ and $z'=(\theta_J,v)$ be elements in $X^\theta$ such that $\tau^\times(z,z')\neq0$ and $zz'\in X_s^\theta$. By~\eqref{180}, the condition $\tau^\times(z,z')\neq0$ is equivalent to $I\cap J=\emptyset$, in which case $\theta_I\theta_J=\theta_{I\cup J}$ and $zz'=(\theta_{I\cup J},uv)$. Consequently,\[\ind(zz')=I\cup J\cup\ind(uv)=I\cup J\cup\ind(u)\cup\ind(v)=\ind(z)\cup\ind(z').\]Thus, the condition $|\ind(zz')\cap[a,b]|\leq1$ forces both $|\ind(z)\cap[a,b]|\leq1$ and $|\ind(z')\cap[a,b]|\leq1$, meaning $z,z'\in X_s^\theta$. This proves that the complement $X^\theta\setminus X_s^\theta$ is an ideal of $X^\theta$ with respect to $\tau^\times$. Furthermore, $\rho_s^\times$ is multiplicative on $X_s^\theta$ because the diagonal action $\rho^\times$ is multiplicative on $X^\theta$. Therefore, by Propositions~\ref{054} and~\ref{055}, the family $\{X_s^\theta\mid s\in S\}$ induces a bigraded quasi-multiplicative restricted action $\bar{\rho}^\times:\bar{\sym}_\Q\to\sym(X^\theta)$ and a restricted weight function $\bar{\omega}^\times:\bar{\sym}_\Q\times X^\theta\to\field^\times$ such that $(\tau^\times,\bar{\omega}^\times)$ is a multiplicative pair. Hence, Condition~\eqref{018} holds for $\bar{\sym}_\Q$.

We set $\sNCQSym_q:=\funcs{\field}{X^\theta}^{\bar{\sym}_\Q}$, the \emph{algebra of $q$-quasisymmetric functions in noncommuting variables in superspace}. In particular, $\sNCQSym:=\sNCQSym_1$ is the algebra of quasisymmetric functions in noncommuting variables in superspace~\cite{PreArGoMa26A}.

An element $(\theta_I,u)\in X^\theta$ is called \emph{packed} if either $(\theta_I,u)=(1,1)$, or $\ind(\theta_I,u)=[m]$ for some $m\geq1$~\cite{PreArGoMa26A}. The \emph{packed word} of $(\theta_I,u)$, denoted by $\pack(\theta_I,u)$, is the unique packed element in $X^\theta$ obtained by replacing the $i$th smallest index in $\ind(\theta_I,u)$ by $i$ in both $I$ and $u$; for instance, $\pack(\theta_{\{-1,4\}},x_{7/2}x_{-1}x_4)=(\theta_{\{1,3\}},x_2x_1x_3)$. Since each $s_{a,b}\in S$ acts nontrivially on $(\theta_I,u)$ under $\bar{\rho}^\times$ only when $\ind(\theta_I,u)$ contains exactly one index in $\{a,b\}$ and no indices in $(a,b)$, the restricted action $\bar{\rho}^\times$ preserves the relative order of $\ind(\theta_I,u)$. Moreover, any two elements in $X^\theta$ sharing the same packed word can be transformed into one another by a finite sequence of such order-preserving reindexings. Consequently, two elements belong to the same $\bar{\sym}_\Q$-orbit if and only if they have the same packed word.

Given $(\theta_I,u)\in X^\theta$ and $\sigma\in(\bar{\sym}_\Q)_{(\theta_I,u)}$, by Proposition~\ref{054}\eqref{060} there exists $\varsigma\in\sym_\Q$ such that $\rho_\varsigma^\times(\theta_I,u)=\bar{\rho}_\sigma^\times(\theta_I,u)=(\theta_I,u)$ and $\bar{\omega}^\times(\sigma,(\theta_I,u))=\omega^\times(\varsigma,(\theta_I,u))$ while $\varsigma$ preserves the relative order of $\ind(\theta_I,u)$. Since the only order-preserving permutation of the finite set $\ind(\theta_I,u)$ into itself is the identity, $\varsigma$ fixes every index in $I$ and in $\ind(u)$, yielding $\bar{\omega}^\times(\sigma,(\theta_I,u))=1$. Therefore, every $\bar{\sym}_\Q$-orbit in $X^\theta$ is $\bar{\omega}^\times$-admissible. We choose as our system of representatives the set $\bar{T}^\theta\subset X^\theta$ of all packed words. Since the restricted action preserves the bigrading, each homogeneous component $(X^\theta)^d$ intersects only finitely many $\bar{\sym}_\Q$-orbits. Thus, Condition~\eqref{051} holds, and by Theorem~\ref{057}, the algebra $\sNCQSym_q$ admits the monomial basis $\{\bar{m}_{(\theta_I,u)}^\times\mid (\theta_I,u)\in\bar{T}^\theta\}$, such that for every $(\theta_J,v)\in T^\theta$,\[m_{(\theta_J,v)}^\times=\sum_{(\theta_I,u)\in T_{(\theta_J,v)}}\omega^\times(\sigma_{(\theta_I,u)},(\theta_J,v))\,\bar{m}_{(\theta_I,u)}^\times,
\]where $T_{(\theta_J,v)}$ is the set of packed words in the orbit $\sym_\Q(\theta_J,v)$, and $\sigma_{(\theta_I,u)}\in\sym_\Q$ satisfies $\rho_{\sigma_{(\theta_I,u)}}^\times(\theta_J,v)=(\theta_I,u)$.

Finally, we verify that the coproduct structure $(\psf_\theta,\hsf)$ on $X^\theta$ defined in Subsubsection~\ref{185} transfers to $\sNCQSym_q$ via Proposition~\ref{170}. As shown in Subsubsection~\ref{182}, for every $s=s_{a,b}\in S$, both $\hsf(s,1)$ and $\hsf(1,s)$ belong to $S$. For any $z=(\theta_I,u)\in X^\theta$ with $\psf_\theta(z)=(z_1,z_2)=((\theta_{I_1},u_1),(\theta_{I_2},u_2))$, because $\iota:\Q\to\Q_1\sqcup\Q_2$ is an order-isomorphism with $\Q_1<\Q_2$, the closed interval $[\iota_1^{-1}(\iota(a)),\iota_1^{-1}(\iota(b))]$ in the first copy is entirely contained in $\iota^{-1}(\Q_1)$ and is disjoint from $\iota^{-1}(\iota_2(\ind(z_2)))$. Since $\iota_1^{-1}\circ\iota$ restricts to an order-isomorphism between $\ind(z)\cap[\iota^{-1}(\iota_1(a)),\iota^{-1}(\iota_1(b))]$ and $\ind(z_1)\cap[a,b]$, their cardinalities coincide, which implies $z\in X_{\hsf(s,1)}^\theta$ if and only if $z_1\in X_s^\theta$. Symmetrically, $z\in X_{\hsf(1,s)}^\theta$ if and only if $z_2\in X_s^\theta$. Furthermore, for any $u\in X^\theta\times X^\theta\times X^\theta$, the order-preserving reindexing between the two partitions of $\Q$ into three consecutive ordered copies $\Q_1<\Q_2<\Q_3$ induces a bijection $f_u:((\psf_\theta\times\id)\circ\psf_\theta)^{-1}(\{u\})\to((\id\times\psf_\theta)\circ\psf_\theta)^{-1}(\{u\})$ that preserves the relative order of all fermionic and bosonic indices; hence, for every $v$ in the fibre, $w=f_u(v)$ shares the exact same packed word as $v$, so $w=\bar{\rho}_\sigma^\times(v)$ for some $\sigma\in\bar{\sym}_\Q$ with $\bar{\omega}^\times(\sigma,v)=1$. Likewise, for any $t'\in\bar{T}^\theta$, if $\psf_\theta(v)=((1,1),t')$ or $\psf_\theta(v)=(t',(1,1))$, the order-isomorphism $\iota$ ensures that $\pack(v)=t'$, so each fibre is a singleton contained in $\bar{\sym}_\Q t'$ with restricted weight $1$. Therefore, by Proposition~\ref{170}, $\sNCQSym_q$ is a connected Hopf algebra, and the canonical inclusion $\sNCSym_q\hookrightarrow\sNCQSym_q$ is an injective homomorphism of Hopf algebras.

\subsubsection{Symmetric functions in commuting variables in superspace}\label{186}

Let $V_\theta=\id_{\Theta^\wedge}\times V$ be the congruence on $X^\theta=\Theta^\wedge\times X^*$ defined by $((\theta_I,u),(\theta_J,v))\in V_\theta$ if and only if $I=J$ and $(u,v)\in V=[X^*,X^*]$. The quotient monoid $X^\theta/V_\theta$ is naturally isomorphic to $X^{\theta,\oplus}:=\Theta^\wedge\times X^\oplus$, the direct product of the fermionic monoid $\Theta^\wedge$ and the free commutative monoid $X^\oplus$. Each congruence class $[(\theta_I,u)]$ is finite of size at most $\ell(u)!$, and $U_\theta=\{(\theta_I,\sort(u))\mid (\theta_I,u)\in X^\theta\}$ forms a canonical system of representatives with retraction $\ssf_\theta=\id_{\Theta^\wedge}\times\sort:X^\theta\to U_\theta$.

We now determine when $(V_\theta,U_\theta)$ defines an invariant-retraction pair on $X^\theta$ as in Subsection~\ref{176}. Since $V$ preserves word lengths and multisets of letters while the fermionic component $\theta_I$ remains fixed under $V_\theta$, the congruence $V_\theta$ is immediately bigraded and stable under the diagonal action $\rho^\times$ of $\sym_\Q$. Moreover, for any $((\theta_I,u),(\theta_I,u'))\in V_\theta$ and $((\theta_J,v),(\theta_J,v'))\in V_\theta$, we have $\sort(u)=\sort(u')$ and $\sort(v)=\sort(v')$, so $\inv(u,v)=\inv(u',v')$, which implies $\tau^\times((\theta_I,u),(\theta_J,v))=\tau^\times((\theta_I,u'),(\theta_J,v'))$. Thus, $V_\theta$ is compatible with the twisting $\tau^\times$.

For the diagonal weight $\omega^\times(\sigma,(\theta_I,u))=\omega_s(\sigma,\theta_I)\,\omega_q(\sigma,u)$ to be compatible with $V_\theta$, we must have $\omega_s(\sigma,\theta_I)\,\omega_q(\sigma,u)=\omega_s(\sigma,\theta_I)\,\omega_q(\sigma,u')$ for all $(u,u')\in V$ and $\sigma\in\sym_\Q$. Since $\omega_s(\sigma,\theta_I)\in\{\pm1\}$ is non-zero, this is equivalent to $\omega_q(\sigma,u)=\omega_q(\sigma,u')$ for all $(u,u')\in V$, which holds if and only if $q\in\{1,-1\}$ as shown in Subsubsection~\ref{194}. Assuming $q=\pm1$, the compatibility of $V_\theta$ with the coproduct map $\psf_\theta$ follows directly from that of $V$ with $\psf$ in Subsubsection~\ref{194}, because $\psf_\theta(\theta_I,u)=((\theta_{I_1},u_1),(\theta_{I_2},u_2))$ acts as $\psf$ on the bosonic component while partitioning $I$ deterministically into $\iota_1(I_1)\sqcup\iota_2(I_2)$. Thus, $(V_\theta,U_\theta)$ is an invariant-retraction pair on $X^\theta$ for $q=\pm1$.

We set $\sSym_q:=\funcs{\field}{X^{\theta,\oplus}}^{\sym_\Q}$, the \emph{algebra of $(\pm1)$-symmetric functions in commuting variables in superspace}; in particular, $\sSym:=\sSym_1$ is the classical algebra of symmetric functions in superspace~\cite{DeLaMa03,DeLaMa04,DeLaMa06,FiLaPi19}. To describe its basis and the graded bi-ideal $J_\theta$, let $T_{V_\theta}^\theta=\{t\in T^\theta\mid[t]\text{ is }(\omega^\times)^{V_\theta}\text{-admissible}\}$. Any class in $X^{\theta,\oplus}=\Theta^\wedge\times X^\oplus$ can be written uniquely as $(\theta_I,x_{r_1}^{a_1}\cdots x_{r_k}^{a_k})$ with $r_1<\cdots<r_k$ and $a_i\geq1$; for convenience, set $a_r=0$ for all $r\in\Q\setminus\{r_1,\dots,r_k\}$. A permutation $\sigma\in\sym_\Q$ belongs to the stabilizer of $(\theta_I,x_{r_1}^{a_1}\cdots x_{r_k}^{a_k})$ in $X^{\theta,\oplus}$ if and only if $\sigma(I)=I$ and $a_{\sigma(r)}=a_r$ for all $r\in\Q$. When $q=1$, we have $\omega_1\equiv1$, so $(\omega^\times)^{V_\theta}(\sigma,[(\theta_I,u)])=\omega_s(\sigma,\theta_I)$. If there exist two distinct indices $r\neq s$ in $I$ with $a_r=a_s$, the transposition $(r\;s)$ belongs to the stabilizer of $[(\theta_I,u)]$ in $X^{\theta,\oplus}$ and satisfies $\omega_s((r\;s),\theta_I)=-1\neq1$. Conversely, if the multiplicities $\{a_r\mid r\in I\}$ are pairwise distinct, every permutation in the stabilizer of $[(\theta_I,u)]$ fixes $I$ pointwise, so $\omega_s(\sigma,\theta_I)=1$. Therefore, $[(\theta_I,u)]$ is admissible in $\sSym_1$ if and only if $a_r\neq a_s$ for all distinct $r,s\in I$. The resulting admissible orbits in $X^{\theta,\oplus}$ are naturally indexed by \emph{superpartitions}~\cite{DeLaMa03,FiLaPi19}. On the other hand, when $q=-1$, a transposition $(r\;s)$ swapping two indices with $a_r=a_s=a$ changes the inversions by an integer congruent to $a^2$ modulo $2$, so $[(\theta_I,u)]$ is admissible in $\sSym_{-1}$ if and only if no odd multiplicity $a_r$ is repeated among indices outside $I$, and no even multiplicity $a_r$ is repeated among indices inside $I$. By Theorem~\ref{153} and Proposition~\ref{166}, $\sSym_q$ is a connected Hopf algebra, isomorphic to the quotient of $\sNCSym_q$ by the graded bi-ideal $J_\theta$ generated by\[\{m_t^\times\mid t\in T^\theta\setminus T_{V_\theta}^\theta\}\cup\{c_{t'}\,m_t^\times-c_t\,(\omega^\times)^{V_\theta}(\sigma,[t])\,m_{t'}^\times\mid t,t'\in T_{V_\theta}^\theta\text{ and }\sigma\in\sym_\Q\text{ such that }(\rho^\times)_\sigma^{V_\theta}([t])=[t']\},\]where $c_t=|[t]\cap\sym_\Q t|\geq1$ for all $t\in T_{V_\theta}^\theta$.

\subsubsection{Quasisymmetric functions in commuting variables in superspace}

Combining the quotient construction of Subsubsection~\ref{186} for $q\in\{1,-1\}$ with the restricted action method of Subsubsection~\ref{187}, we work over the monoid $X^{\theta,\oplus}=\Theta^\wedge\times X^\oplus=X^\theta/V_\theta$, equipped with the induced twisting $(\tau^\times)^{V_\theta}([(\theta_I,u)],[(\theta_J,v)])=\tau_\wedge(\theta_I,\theta_J)\,\tau_q(\sort(u),\sort(v))$, the bigrading $|[(\theta_I,u)]|=(\bar{|I|},\ell(u))$ and $\deg([(\theta_I,u)])=|I|+\ell(u)$, the induced diagonal action $(\rho^\times)^{V_\theta}:\sym_\Q\to\Aut(X^{\theta,\oplus})$, and the induced diagonal weight $(\omega^\times)^{V_\theta}(\sigma,[(\theta_I,u)])=\omega_s(\sigma,\theta_I)\,\omega_q(\sigma,\sort(u))$. For each transposition $s=s_{a,b}\in S$ with $a<b$ in $\Q$, we define the subset\[X_s^{\theta,\oplus}=\{(\theta_I,u)\in X^{\theta,\oplus}\mid |(I\cup\ind(u))\cap[a,b]|\leq1\}.\]As in Subsubsection~\ref{187}, each $X_s^{\theta,\oplus}$ is stable and multiplicative with respect to $(\tau^\times)^{V_\theta}$ under $(\rho^\times)_s^{V_\theta}$, and its complement $X^{\theta,\oplus}\setminus X_s^{\theta,\oplus}$ is an ideal with respect to $(\tau^\times)^{V_\theta}$. Propositions~\ref{054} and~\ref{055} thus yield a bigraded quasi-multiplicative restricted action $\bar{\rho}^{\theta,\oplus}:\bar{\sym}_\Q\to\sym(X^{\theta,\oplus})$ and a restricted weight function $\bar{\omega}^{\theta,\oplus}:\bar{\sym}_\Q\times X^{\theta,\oplus}\to\field^\times$. Because the restricted action preserves the strict relative order of $I\cup\ind(u)$, both sorted representatives experience zero relative inversions under $\bar{\sym}_\Q$, which implies that $\bar{\omega}^{\theta,\oplus}\equiv1$ is trivial and every $\bar{\sym}_\Q$-orbit in $X^{\theta,\oplus}$ is admissible. A canonical system of representatives $\bar{T}^{\theta,\oplus}\subset X^{\theta,\oplus}$ consists of all elements $(\theta_I,x_1^{a_1}\cdots x_m^{a_m})$ with $m\geq0$, $I\subseteq[m]$, and $a_1,\dots,a_m\geq0$ such that $a_i\geq1$ whenever $i\notin I$, so that $I\cup\ind(x_1^{a_1}\cdots x_m^{a_m})=[m]$, which are naturally indexed by \emph{supercompositions} or \emph{dotted compositions}~\cite{FiLaPi19,PreArGoMa26A}.

We set $\sQSym_q:=\funcs{\field}{X^{\theta,\oplus}}^{\bar{\sym}_\Q}$, the \emph{algebra of $(\pm1)$-quasisymmetric functions in commuting variables in superspace}; in particular, $\sQSym:=\sQSym_1$ is the algebra of quasisymmetric functions in superspace~\cite{FiLaPi19}. By Proposition~\ref{170} and Theorem~\ref{153}, the coproduct structure on $X^\theta$ descends and restricts to $X^{\theta,\oplus}$, making $\sQSym_q$ a connected Hopf algebra that fits into the commutative square of connected Hopf algebras:\[\begin{array}{ccc}
\sNCSym_q&\hookrightarrow&\sNCQSym_q\\[4pt]
\downarrow&&\downarrow\\[4pt]
\sSym_q&\hookrightarrow&\sQSym_q
\end{array}\]

\subsection{Additional examples}\label{193}

\subsubsection{Symmetric functions in partially commutative variables}

Following~\cite{Dol21,PreDau24}, let $A$ be a finite totally ordered alphabet of elements called \emph{colors}, and consider the partially commutative monoid\[M_A=\langle x_{a,r}\mid a\in A,\;r\in\Q,\;x_{a,r}x_{b,s}=x_{b,s}x_{a,r}\text{ whenever }r\neq s\rangle.\]Equivalently, $M_A$ is the quotient of the free monoid generated by $X_A=\{x_{a,r}\mid a\in A,\;r\in\Q\}$ by the congruence generated by the pairs $(x_{a,r}x_{b,s},x_{b,s}x_{a,r})$ with $r\neq s$. Thus, generators with distinct second indices commute, whereas generators sharing the same second index do not commute. For a word $w=a_1\cdots a_m\in A^*\setminus\{1\}$ and $r\in\Q$, we write $x_{w,r}=x_{a_1,r}\cdots x_{a_m,r}$. Every non-identity element $u\in M_A$ admits a unique normal form\[u=x_{w_1,r_1}\cdots x_{w_k,r_k},\quad\text{with }r_1<\cdots<r_k\text{ in }\Q\text{ and }w_1,\dots,w_k\in A^*\setminus\{1\}.\]The sequence of words $I=(w_1,\dots,w_k)$ is called the \emph{sentence} associated with $u$, and we set $\ind(u)=\{r_1,\dots,r_k\}\subset\Q$. We equip $M_A$ with the trivial twisting $\tau\equiv1$, the trivial color group $\G=\{1\}$ with bicharacter $\phi\equiv1$, and the degree map $\deg:M_A\to\N$ given by the total number of colors, $\deg(x_{w_1,r_1}\cdots x_{w_k,r_k})=\ell(w_1)+\cdots+\ell(w_k)$. Any factorization of an element in normal form is obtained by splitting each word $w_j$ into a prefix and a suffix, possibly empty, yielding only finitely many factorizations. Hence, $M_A$ is locally finite and graded with $M_A^0=\{1\}$, so Conditions~\eqref{015} to~\eqref{017} hold.

The finitary symmetric group $\sym_\Q$ acts on $M_A$ by bigraded monoid automorphisms via the symmetrizing action $\rho:\sym_\Q\to\Aut(M_A)$ defined on generators by $\rho_\sigma(x_{a,r})=x_{a,\sigma(r)}$. Taking the trivial weight function $\omega\equiv1$, Condition~\eqref{018} is immediately satisfied. Two elements in $M_A$ belong to the same $\sym_\Q$-orbit if and only if their associated sentences share the same multiset of words $\{w_1,\dots,w_k\}$. Ordering $A^*$ by the \emph{graded lexicographic order} $\preceq_{\text{glex}}$, where $v\preceq_{\text{glex}}w$ if $\ell(v)<\ell(w)$, or if $\ell(v)=\ell(w)$ and $v$ precedes $w$ lexicographically, each orbit contains a unique representative $u_P=x_{w_1,1}\cdots x_{w_k,k}$ whose associated sentence $P=(w_1,\dots,w_k)$ is sorted, that is, $w_1\succeq_{\text{glex}}\cdots\succeq_{\text{glex}}w_k$; such a sorted sentence is called a \emph{p-sentence}~\cite{PreDau24}. Since $A$ is finite, there are only finitely many p-sentences of each degree $d\in\N$. Thus, Condition~\eqref{051} holds for the system of representatives $T_A=\{u_P\mid P\text{ is a p-sentence}\}$. We set $\Sym_A:=\funcs{\field}{M_A}^{\sym_\Q}$, the algebra of \emph{colored symmetric functions in partially commutative variables}~\cite{PreDau24}, with monomial basis $\{m_P:=m_{u_P}\mid u_P\in T_A\}$. To endow $\Sym_A$ with a Hopf structure via Subsubsection~4.4.1, let $\iota:\Q\to\Q_1\sqcup\Q_2$ with $\Q_1<\Q_2$ and $\hsf:\sym_\Q\times\sym_\Q\to\sym_\Q$ be the maps from Subsubsection~5.1.1. Since elements in $M_A$ supported on $\Q_1$ and $\Q_2$ commute with each other, $\iota$ induces a bigraded monoid isomorphism $\psi_A:M_A\to(M_A*M_A)/[M_A,M_A]\simeq M_A\times M_A$. Taking $\psf_A=\psi_A$, the fibres of $\psf_A$ are singletons, and $(\psf_A,\hsf)$ forms a coproduct structure on $M_A$ satisfying Condition~\eqref{127}. By Proposition~\ref{052} and Theorem~\ref{130}, $\Sym_A$ is a connected Hopf algebra with product and coproduct on the monomial basis given by\[m_P\cdot m_Q=\sum_R c_{P,Q}^R\,m_R\quad\text{and}\quad\Delta(m_P)=\sum_{Q\sqcup R=P}m_Q\otimes m_R,\]where $c_{P,Q}^R=c_{u_P,u_Q}^{u_R}$ is given by~\eqref{053} and $Q\sqcup R=P$ runs over all ordered pairs of p-sentences $(Q,R)$ whose disjoint union as multisets of words equals $P$, recovering the Hopf algebra structure of~\cite[Theorem~4.6]{PreDau24}.

Furthermore, applying the restricted action method of Subsection~\ref{171} with the transpositions $S=\{s_{a,b}\mid a<b\text{ in }\Q\}$ and the subsets $M_{A,s}=\{u\in M_A\mid|\ind(u)\cap[a,b]|\leq1\}$ yields the restricted action $\bar{\rho}:\bar{\sym}_\Q\to\sym(M_A)$ with trivial restricted weight $\bar{\omega}\equiv1$. Two elements of $M_A$ lie in the same $\bar{\sym}_\Q$-orbit if and only if they have the same associated sentence $I=(w_1,\dots,w_k)$, with canonical representative $u_I=x_{w_1,1}\cdots x_{w_k,k}$. By Theorem~\ref{057} and Proposition~\ref{170}, the resulting invariant space $\QSym_A:=\funcs{\field}{M_A}^{\bar{\sym}_\Q}$ is the connected \emph{Hopf algebra of colored quasisymmetric functions in partially commutative variables}~\cite{Dol21,PreDau24}, and the canonical inclusion $\Sym_A\hookrightarrow\QSym_A$ is an injective homomorphism of Hopf algebras.

\subsubsection{Signed-symmetric functions}

Consider the signed alphabet $X_\pm=\{x_r\mid r\in\Q^\times\}$ together with the free monoid $X_\pm^*$ generated by $X_\pm$. We equip $X_\pm^*$ with the trivial twisting $\tau\equiv1$, the trivial color group $\G=\{1\}$ with bicharacter $\phi\equiv1$, and the word-length grading $\deg(u)=\ell(u)$, so that Conditions~\eqref{015} to~\eqref{017} are satisfied.

Let $\sym_\Q^B$ denote the \emph{finitary hyperoctahedral group} or \emph{finitary signed symmetric group} on $\Q^\times$, that is, the group of all bijections $\sigma:\Q^\times\to\Q^\times$ that fix all but finitely many elements and commute with sign inversion, meaning $\sigma(-r)=-\sigma(r)$ for all $r\in\Q^\times$. The \emph{signed symmetrizing action} $\rho^B:\sym_\Q^B\to\Aut(X_\pm^*)$ is the bigraded multiplicative action defined on words $u=x_{r_1}\cdots x_{r_d}\in X_\pm^*$ by\[\rho_\sigma^B(x_{r_1}\cdots x_{r_d})=x_{\sigma(r_1)}\cdots x_{\sigma(r_d)}\quad\text{for all }\sigma\in\sym_\Q^B.\]Taking the trivial weight function $\omega\equiv1$, Condition~\eqref{018} holds, and we define $\NCSym_B:=\funcs{\field}{X_\pm^*}^{\sym_\Q^B}$, the \emph{algebra of signed-symmetric functions in noncommuting variables}.

Since $\omega\equiv1$, every $\sym_\Q^B$-orbit in $X_\pm^*$ is admissible. Observe that two words $u=x_{r_1}\cdots x_{r_d}$ and $v=x_{s_1}\cdots x_{s_d}$ belong to the same $\sym_\Q^B$-orbit if and only if for all $i,j\in[d]$, we have $|r_i|=|r_j|\Leftrightarrow |s_i|=|s_j|$ and $r_i=r_j\Leftrightarrow s_i=s_j$. Consequently, each orbit contains a unique representative $t=x_{r_1}\cdots x_{r_d}$ such that the word of absolute values $x_{|r_1|}\cdots x_{|r_d|}\in T_1$ is a restricted growth word on $\{x_1,x_2,\ldots\}$ and the first occurrence of each index $|r_i|$ in $t$ carries a positive sign $r_i>0$. Denoting by $T_B\subset\{x_{\pm1},x_{\pm2},\dots\}^*$ the set of all such \emph{signed restricted growth words}, the elements of $T_B$ of length $d$ are in canonical bijection with signed-index set partitions of $[d]$, of which there are $\sum_{k=0}^d2^{d-k}S(d,k)$, where $S(d,k)$ denotes the Stirling number of the second kind. Thus, Condition~\eqref{051} holds, and $\NCSym_B$ admits the monomial basis $\{m_t^B\mid t\in T_B\}$.

To equip $\NCSym_B$ with a Hopf structure, let $\Q_+=\{r\in\Q\mid r>0\}$ and choose an order-isomorphism $\iota_+:\Q_+\to\Q_{+,1}\sqcup\Q_{+,2}$ with $\Q_{+,1}<\Q_{+,2}$. Extending $\iota_+$ by sign symmetry via $\iota(-r)=-\iota_+(r)$ yields a bijection $\iota:\Q^\times\to\Q_1^\times\sqcup\Q_2^\times$, where $\Q_s^\times=\Q_{+,s}\sqcup(-\Q_{+,s})$ for $s\in\{1,2\}$. As in Subsubsection~\ref{123}, $\iota$ induces a monoid isomorphism $\psi_B:X_\pm^*\to X_\pm^**X_\pm^*$ and an injective group homomorphism $\hsf_B:\sym_\Q^B\times\sym_\Q^B\to\sym_\Q^B$ given by $\hsf_B(\sigma,\varsigma)=\iota^{-1}\circ(\sigma\sqcup\varsigma)\circ\iota$. Composing $\psi_B$ with the canonical projection $\qsf_B:X_\pm^**X_\pm^*\to X_\pm^*\times X_\pm^*$ defines $\psf_B=\qsf_B\circ\psi_B$. Because $\omega\equiv1$ and the reindexing between three copies of $\Q^\times$ commutes with sign inversion, $(\psf_B,\hsf_B)$ satisfies Condition~\eqref{127}, endowing $\NCSym_B$ with the structure of a connected Hopf algebra.

Finally, our general constructions in Section~\ref{168} immediately produce the remaining members of this signed family:
\begin{itemize}
\item Restricted actions ($\NCQSym_B$): Let $S_B$ be the symmetric generating set of $\sym_\Q^B$ consisting of the signed transpositions $s_{a,b}=(a\;b)(-a\;-b)$ for $0<a<b$ in $\Q_+$ and the sign changes $s_a=(a\;-a)$ for $a\in\Q_+$. Setting $(X_\pm^*)_{s_{a,b}}=\{u\in X_\pm^*\mid |\{|r|\mid x_r\in\alph(u)\}\cap[a,b]|\leq1\}$ and $(X_\pm^*)_{s_a}=X_\pm^*$, Subsection~\ref{171} and Proposition~\ref{170} yield a connected Hopf algebra of \emph{signed-quasisymmetric functions in noncommuting variables} $\NCQSym_B:=\funcs{\field}{X_\pm^*}^{\bar{\sym}_\Q^B}$.
\item Commutative quotients ($\Sym_B$ and $\QSym_B$): Let $V_B=[X_\pm^*,X_\pm^*]$ be the commutation congruence on $X_\pm^*$, with quotient $X_\pm^*/V_B\simeq X_\pm^\oplus$ and sorted representatives $U_B=\{\sort(u)\mid u\in X_\pm^*\}$. Since $(V_B,U_B)$ is an invariant-retraction pair with respect to $\sym_\Q^B$, Theorem~\ref{153} yields the quotient connected Hopf algebra $\Sym_B:=\funcs{\field}{X_\pm^\oplus}^{\sym_\Q^B}$ of \emph{signed-symmetric functions in commuting variables}, whose monomial basis is indexed by finite multisets of unordered pairs $\{a,b\}$ of nonnegative integers with $a+b>0$ (recording the exponents of $x_r$ and $x_{-r}$ for each active index $r>0$), together with its quasisymmetric counterpart $\QSym_B:=\funcs{\field}{X_\pm^\oplus}^{\bar{\sym}_\Q^B}$, whose monomial basis is indexed by finite sequences of such unordered pairs $\{a,b\}$.
\end{itemize}

%\begin{itemize}
%\item[\eqref{015}]\conditions{0}
%\item[\eqref{016}]\conditions{1}
%\item[\eqref{017}]\conditions{2}
%\item[\eqref{018}]\conditions{3}
%\item[\eqref{051}]\conditions{4}
%\item[\eqref{127}]\conditions{5}
%\item[\eqref{161}]\conditions{6}
%\end{itemize}

\subsubsection*{Acknowledgements}
The first named author acknowledges the financial support of DIDULS/ULS, through
the project PR2553853. The second named author was partially supported by the grant ANID-FONDECYT Iniciaci\'on No. 11241418. The third named author acknowledges the financial support of Fondo de Apoyo a la Investigaci\'on DIUA309-2025.

\bibliographystyle{plainurl}
\bibliography{../../../../Projects/LaTeX/bibtex.bib}

\end{document}